\documentclass[a4paper,12pt,reqno]{amsart}

\usepackage{amsmath,amssymb,amsthm,bm}
\usepackage{caption,graphicx}
\usepackage[english]{babel}
\usepackage{amscd}
\usepackage{amsgen}
\usepackage[final]{epsfig}
\usepackage{latexsym}
\usepackage{amsfonts}
\usepackage{url}
\usepackage{slashed}
\usepackage[all]{xy}
\usepackage{here}
\usepackage{comment}
\usepackage{fixmath}
\usepackage{scalerel,stackengine}
\stackMath
\newcommand\reallywidehat[1]{%
\savestack{\tmpbox}{\stretchto{%
  \scaleto{%
    \scalerel*[\widthof{\ensuremath{#1}}]{\kern-.6pt\bigwedge\kern-.6pt}%
    {\rule[-\textheight/2]{1ex}{\textheight}}
  }{\textheight}%
}{0.5ex}}%
\stackon[1pt]{#1}{\tmpbox}%
}
\newtheorem{thm}{Theorem}
\newtheorem{lemm}{Lemma}
\newtheorem{prop}{Proposition}
\newtheorem{cor}{Corollary}

\newtheorem{defn}{Definition}
\newtheorem{probl}{Problem}

\newtheorem{ex}{Example}

\allowdisplaybreaks

\newcommand{\toitself}{\mathbin{\scalebox{.85}{%
    \lefteqn{\scalebox{.5}{$\blacktriangleleft$}}\raisebox{.34ex}{$\supset$}}}}

\begin{document}

\title{Gromov's Oka principle, fiber bundles and the conformal module}

\author{Burglind J\"oricke }

\address{IHES, 35 Route de Chartres\\ 91440 Bures-sur-Yvette\\ France}

\email{joericke@googlemail.com}

\date{Received: date / Accepted: date}

\dedicatory{To the memory of my colleague Mikael Passare
who died tragically in 2011}



\begin{abstract} The conformal module of conjugacy classes of braids is an invariant that appeared earlier than the entropy of conjugacy classes of braids, and is inverse proportional to the entropy. Using the relation between the two invariants we give a short conceptional proof of an earlier result on the conformal module.
Mainly, we consider situations, when the conformal module of conjugacy classes of braids
serves as obstruction for the existence of homotopies (or isotopies) of smooth objects involving braids to the respective holomorphic objects, and present theorems on the restricted validity of Gromov's Oka principle in these situations.
\end{abstract}

\maketitle

\centerline \today


\section{Introduction}
\label{section:8.1}
\setcounter{equation}{0}

The conformal module of conjugacy classes of braids
appeared first (without name) in the paper \cite{GL} in connection with the interest of the authors in Hilbert's Thirteen's Problem. This invariant of conjugacy classes of braids
is undeservedly almost forgotten, although it appeared again in \cite{G} (for conjugacy classes of elements of the findamental group of the twice punctured complex plane instead of conjugacy classes of braids). It defines obstructions for the validity of Gromov's Oka principle
concerning homotopies of continuous or smooth objects involving braids to the respective holomorphic objects. The interesting point is that this invariant is inverse proportional to a popular dynamical invariant of conjugacy classes of braids, the entropy (see \cite{Jo})). Starting with the paper \cite{FLP} the entropy has been studied extensively. In some cases it has been computed explicitly. In the present paper we give applications of the concept of the conformal module of conjugacy classes of braids. We will study situations when this invariant gives obstructions for the validity of
Gromov's Oka principle and present theorems on the restricted validity of this principle. We will also give a short conceptional proof of (a slight improvement of) the first result \cite{GL} related to the conformal module of conjugacy classes of braids using the relation of this invariant to the entropy.

We will consider braids as elements of the fundamental group 
of the symmetrized
configuration space. More, detailed,
let $C_n(\mathbb{C})=\{(z_1,\ldots,z_n)\in \mathbb{C}^n: z_j\neq z_k \;\mbox{for}\; j\neq k\}$ be the $n$-dimensional configuration space.
The symmetrized configuration space is its quotient $C_n(\mathbb{C})\diagup \mathcal{S}_n$ by the diagonal action of the symmetric group $\mathcal{S}_n$.
It carries a complex structure as a quotient of a complex manifold by a free and discontinuous action of a group.
A point in $C_n(\mathbb{C})\diagup \mathcal{S}_n$ can also be considered as subset of $\mathbb{C}$ consisting of $n$ points.
A geometric braid on $n$ strands (an $n$-braid for short) with base point $E_n$ is a collection of $n$ disjoint arcs in $[0,1]\times \mathbb{C}$ (the strands) which join the set $\{0\}\times E_n$ with the set $\{1\}\times E_n$, so that the projection of each strand to the first factor $[0,1]$ is a homeomorphism.
A braid on $n$ strands ($n$-braids for short)
with base point $E_n\in C_n(\mathbb{C})\diagup \mathcal{S}_n$
is an isotopy class of geometric braids with base point $E_n$, equivalently,
it is a homotopy class of loops with
base point $E_n$ in the symmetrized configuration space,
i.e. an element of the fundamental group $\pi_1(C_n (\mathbb {C})
\diagup
{\mathcal S}_n, E_n)$ of the symmetrized configuration space
with base point $E_n$.
Alternatively, braids on $n$ strands are given by words in the
Artin braid group \index{Artin group} ${\mathcal B}_n$.

The conjugacy classes $\widehat{ \mathcal{B}}_n$
\index{$\widehat{\mathcal{B}}_n$} of elements of the braid group, equivalently, of
the fundamental group $\pi_1( C_n(\mathbb{C})\diagup \mathcal{S}_n,E_n)$ of the symmetrized configuration space, can be
interpreted as free homotopy classes of loops in $C_n(\mathbb{C})\diagup \mathcal{S}_n$.
We say that a continuous mapping $f$ of a round
annulus $A= \{z \in \mathbb{C}: \, r<|z|<R\},\;$ $0\leq r < R \leq \infty,\;$
into $ C_n(\mathbb{C})\diagup \mathcal{S}_n  $ represents an element $\hat b\in \widehat{ \mathcal{B}}_n$ if for some (and hence for
any) circle $\,\{|z|=\rho \} \subset A\,$ the loop $\,f:\{|z|=\rho
\} \rightarrow  C_n(\mathbb{C})\diagup \mathcal{S}_n\,$ represents $\,\hat b$.

More generally,
let $\mathcal{X}$ be a topological space, and $\hat e$ a conjugacy class of elements of the fundamental group of $\mathcal{X}$. $\hat e$ can be interpreted as a free homotopy class of loops in $\mathcal{X}$.
A continuous map $g$ from an
annulus $A = \{ z \in {\mathbb C} : r < \vert z \vert < R \}$ into
${\mathcal X}$ is said to represent the conjugacy class $\hat e$
if for some (and hence for any) $\rho \in (r,R)$ the map $g :
\{\vert z \vert = \rho \} \to {\mathcal X}$ represents $\hat e$.
Let $\mathcal{X}$ be a smooth oriented manifold. If the mapping is a homeomorphism we will speak about an annulus representing $\hat e$ contained in $\mathcal{X}$.

Ahlfors \cite{A1} defined the conformal module \index{conformal module ! of an
annulus} of an annulus $\,A= \{z \in \mathbb{C}:\; r<|z|<R\}\,\,$ as
$\,m(A)= \frac{1}{2\pi}\,\log(\frac{R}{r})\,.\, $ \index{$m(A)$} Two
annuli of finite conformal module are conformally equivalent iff
they have equal conformal module. If a manifold $\Omega$ is
conformally equivalent to a round annulus $A$, its conformal module is
defined to be $m(A)$. Recall that any domain in the complex plane
with fundamental group isomorphic to the group of integer numbers
$\mathbb{Z}$ is conformally equivalent to an annulus.

Associate with each conjugacy class of the fundamental group of
$C_n(\mathbb{C})\diagup \mathcal{S}_n  $, or, equivalently, with each conjugacy class of
$n$-braids, its conformal module defined as follows.

\index{conformal module ! of a conjugacy class of braids}

\begin{defn}\label{defn8.1} Let $\hat b$ be a
conjugacy class of $n$-braids, $n \geq 2$. The conformal module
$M(\hat b)$ of $\hat b$ is defined as $ M(\hat b)=
sup_{\mathcal{A}(\hat b)}\, m(A),$ where $\mathcal{A}(\hat b)$
denotes the set of all annuli which admit a holomorphic mapping into
$ C_n(\mathbb{C})\diagup \mathcal{S}_n$ which represents $\hat b$.
\end{defn}

More generally, for any complex manifold $\mathcal{X}$ and any conjugacy class $\hat e$ of elements of the fundamental group of $\mathcal{X}$ the conformal module $\mathcal{M}(\hat{e})$ of $\hat e$ is defined as the supremum of the conformal modules of annuli that admit a holomorphic mapping into $\mathcal{X}$ representing $\hat e$.

There is another useful point of view on the symmetrized configuration space.
A monic polynomial of degree $n$
can be
parameterized either by the unordered
tuple of its zeros or by its coefficients. Denote by $\mathfrak{P}_n$ the set of monic polynomials of degree $n$ without multiple zeros, and by $\overline{\mathfrak{P}_n}$ the set of all monic polynomials of degree $n$.
The set $\mathfrak{P}_n$ \index{$\mathfrak{P}_n$}
can be identified with
the symmetrized configuration space $C_n(\mathbb{C})\diagup \mathcal{S}_n$ by assigning to each point in $C_n(\mathbb{C})\diagup \mathcal{S}_n$ (an unordered tuple of points in $\mathbb{C}$) the monic polynomial whose zero set coincides with this unordered tuple.
Equivalently, assigning to each polynomial in $\mathfrak{P}_n$ the ordered $n$-tuple of its coefficients, we may identify  $\mathfrak{P}_n$ with the complement of the algebraic
hypersurface $\{\textsf{D}_n =0\}$ in complex Euclidean space
$\mathbb{C}^n$. \index{$\mathbb{C}^n$} Here $\textsf{D}_n$
\index{$\textsf{D}_n$} denotes the discriminant of polynomials in
$\overline{\mathfrak{P}_n}$. The function $\textsf{D}_n$ is a polynomial in
the coefficients of the elements of $\overline{\mathfrak{P}_n}$, that vanishes exactly on the elements of $\overline{\mathfrak{P}_n}$ with multiple zeros.
The complex structure on $\mathbb{C}^n\setminus \{\textsf{D}_n =0\}$ inherited from
$\mathbb{C}^n$ coincides with the complex structure of $C_n(\mathbb{C})\diagup \mathcal{S}_n$ under the identification of the two sets.

In connection with his interest in the Thirteen's Hilbert Problem Arnold studied
the topological invariants of the space $\mathfrak{P}_n$ (\cite{Ar}). On the other hand, the collection of conformal modules of conjugacy classes of elements of the fundamental group of the space of monic polynomials of degree $n$ is a collection of conformal invariants of $\mathfrak{P}_n\cong \mathbb{C}_n\setminus \textsf{D}_n$.

We will consider mappings from topological spaces to the space $\mathfrak{P}_n$.
A continuous mapping from a
topological space $X$ into the set of monic polynomials $\overline{\mathfrak{P}}_n$ of fixed
degree $n$ (maybe, with multiple zeros) is called a quasipolynomial.
\index{quasipolynomial} It can be written as a function in two
variables $\,f(x,\zeta)= a_0(x)
+ a_1(x)\zeta+ ... + a_{n-1}(x)\zeta^{n-1}
+ \zeta^n, \,$
 $x \in X,\; \zeta \in \mathbb{C}\;,$
for continuous functions $a_j, \; j=1,...,n\,,$ on $X$. If $X$ is a
complex manifold and the mapping is holomorphic it is called an
algebroid function. If the image of the map is contained in the space
$\mathfrak{P}_n$ of monic polynomials of degree $n$ without multiple zeros, it is called separable. A separable quasipolynomial
\index{quasipolynomial ! separable} is called solvable
\index{quasipolynomial ! solvable} if it can be globally written as
a product of quasipolynomials of degree $1$.
It is called
irreducible \index{quasipolynomial ! irreducible} if it can not be
written as product of two quasipolynomials of positive degree, and reducible otherwise.
We also call a solvable quasipolynomial a globally solvable family of polynomials, and an irreducible quasipolynomial a globally irreducible family of polynomials.
Two
separable quasipolynomials are called isotopic (or homotopic through separable quasipolynomials) if there is a continuous
family of separable quasipolynomials \index{isotopy ! of
quasipolynomials} joining them. An algebroid function
\index{algebroid function} on $\mathbb C ^n$ whose coefficients are
polynomials is called an algebraic function.

The first result related to the conformal module of braids was a theorem on reducibility of holomorphic quasipolynomials due to Gorin, Lin, Petunin, and Zjuzin.
Let $f$ be an algebroid function on an annulus $A=\{z\in \mathbb{C}: r<|z|<R\}$, $0\leq r<1 < R\leq \infty$. For each $z\in A$ the mapping $f_z(\zeta)= f(z,\zeta),\, \zeta \in \mathbb{C},$ is a polynomial.
We denote by ${\sf D}_n(f_z)$ its discriminant. The index of the mapping
$$
A \ni z \to {\sf D}_n (f_z) \in {\mathbb C} \backslash \{0\} \, ,
\quad z \in A \, ,
$$
is the degree of the map
$$
z \to \frac{{\sf D}_n (f_z)}{\vert {\sf D}_n (f_z) \vert}
$$
from $\{ \vert z \vert = 1\}$ to itself.
We give a short conceptional proof of the following theorem, which is the theorem of these authors with an improved constant.
\begin{thm}\label{thm8.3}
Suppose
$n$ is a prime number. If $A=\{z\in \mathbb{C}: r<|z|<R\}$, $0\leq r<1 < R\leq \infty$, is an annulus of conformal module
$m(A) > \frac{2\pi}{\log 2} \, n$, and $f$ is a separable algebroid function on $A$ of degree $n$ such that the index of its discriminant is divisible
by $n$, then $f$ is reducible.
\end{thm}

For simplicity we
restricted ourselves to the case when the degree of the
quasipolynomial is a prime number. (The general case has been
considered by Lin and later by Zjuzin. An alternative proof of the
general case can be given as well.) In \cite{GL} Gorin and Lin
proved the existence of lower bounds $r_n$ of the conformal module $m(A)$ for which the statement of the theorem is true.
Zjuzin \cite{Z} proved that one can take $r_n = n \cdot \rho_0$ for an
absolute constant $\rho_0$, and
Petunin showed that $\rho_0=10^7$ works.

Our main goal is to apply the concept of conformal module to problems concerning the restricted validity of Gromov's Oka principle. We will
consider a continuous or smooth mapping $f$ from an open Riemann surface
(i.e. from a Riemann surface without closed connected components)
to the complex manifold $C_n(\mathbb{C})\diagup \mathcal{S}_n\cong \mathfrak{P}_n$ and ask whether the mapping is homotopic to a holomorphic mapping. This is
the question about the restricted validity of Gromov's Oka principle in this situation.
We will also ask the respective question for fiber bundles, namely, consider a smooth bundle over a Riemann surface whose fibers are closed surfaces, is the bundle isotopic to a holomorphic one.

Gromov  \cite{G} formulated his Oka principle as "an expression of an optimistic expectation with regard to the validity of the $h$-principle for holomorphic maps in the situation when the source manifold is Stein". Holomorphic maps $X\to Y$ from a
complex manifold $X$ to a complex manifold $Y$ are said to satisfy the $h$-principle if each continuous map from $X$ to $Y$ is homotopic to a holomorphic map. We call a target manifold $Y$ a Gromov-Oka manifold if the $h$-principle holds for holomorphic maps from any Stein manifold to $Y$.
Gromov \cite{G} gave a sufficient condition on a complex manifold $Y$ to be a Gromov-Oka manifold.

The question of understanding Gromov-Oka manifolds received a lot of attention. It turned out to be fruitful to strengthen the requirement on the target $Y$ by combining the 
$h$-principle for holomorphic maps
with a Runge type approximation property. Manifolds $Y$ satisfying the stronger requirement are called Oka manifolds. It has been proved by Forstneric that a manifold is an Oka manifold iff Runge approximation holds for holomorphic mappings from neighbourhoods of compact convex subsets of $\mathbb{C}^n$ into $Y$ (see \cite{Forst1}).
For more details, examples of Oka manifolds and an account on modern development of Oka theory based on Oka manifolds see \cite{Forst}.

Open Riemann surfaces are Stein manifolds (see e.g. \cite{Fo}), but the space $\mathfrak{P}_n$ is not a Gromov-Oka manifold. The results of this paper can be interpreted as a beginning of an investigation of the restricted validity of Gromov's Oka principle and of obstructions for its validity in the case when the target is a complex manifold that is not a Gromov-Oka manifold but is of interest in other aspects.

In \cite{G} Gromov mentioned maps from annuli to the twice punctured complex plane as simplest interesting example where the $h$-principle for holomorphic mappings fails (see Section \ref{sec:9.4a} for the case when the target is the twice punctured complex plane). Further, Gromov remarks that for a complex manifold $Y$ which is not a Gromov-Oka manifold the mappings from annuli to $Y$ play a special role for understanding the "conformal rigidity" of the manifold $Y$. He mentions the conformal invariant of $Y$ which is called here the conformal module of a conjugacy class of elements of the fundamental group of $Y$. For $Y=\mathfrak{P}_n$ this invariant is the conformal module of a conjugacy class of $n$-braids.
The problem is interesting for other target spaces $Y$, for instance for complements of arbitrary algebraic hypersurfaces of $\mathbb{C}^n$. The collection of conformal
modules of all conjugacy classes of elements of its fundamental group $\pi_1(Y,y_0)$ is a biholomorphic invariant of the complex
manifold $Y$, which is expected to play a special role for understanding the problem.
Obstructions for Gromov's Oka principle are based
on the relation between conformal invariants of the source and the target.

Runge's approximation theorem on Riemann surfaces implies immediately the following fact.
A surface is called finite if its fundamental group is finitely generated.  Let$f$ be a continuous mapping from a finite open Riemann surface $X$ to $\mathfrak{P}_n$.
There exists a domain $ X' \subset X$ (depending on $f$) which is a deformation retract of $X$ and is
diffeomorphic to $X$ such that the restriction $f \mid  X'$
is homotopic
to a holomorphic mapping from $ X'$ to $\mathfrak{P}_n$.
Slightly reformulated, for each continuous mapping from a
finite open Riemann surface $X$ to $\mathfrak{P}_n$ there exists a homotopic conformal
structure on $X$ for which the mapping is homotopic
to a holomorphic one.
This fact is known in more general situations as ''soft Oka principle" (see \cite{Forst1}).

The following considerations give more precise information in the case
when $X$ is an open Riemann surface of the ''simplest non-trivial topological type'', namely $X$ is equal to an annulus $A$,
and motivate further questions.
For an annulus $A$ the question which continuous mappings from $A$ to $\mathfrak{P}_n$ are homotopic to a holomorphic mapping has an almost complete answer.
Indeed, a continuous mapping $f:A\to \mathfrak{P}_n$ represents a conjugacy class
$\hat{b}_f$ of braids. By the definition of the conformal module, $f$ is homotopic to a holomorphic mapping if $m(A)<\mathcal{M}(\hat{b}_f)$ and is not homotopic to a holomorphic mapping if  $m(A)>\mathcal{M}(\hat{b}_f)$.
Thus there is a measure of "complexity" of the mapping $f:A\to \mathfrak{P}_n$, namely the conformal module $\mathcal{M}(\hat{b}_f)$, whose relation to the conformal invariant $m(A)$ of $A$ "almost" determines whether $f$ is homotopic to a holomorphic mapping.
Recall that $\mathcal{M}(\hat{b}_f)$ is inverse proportional to the entropy of $\hat{b}_f$. It is reasonable to interpret the entropy as "complexity" which grows with increasing entropy, while the chance for the mapping to be homotopic to a holomorphic one decreases. For $n=3$ there are effective upper and lower bounds for the conformal module (equivalently, for the entropy) of $3$-braids that differ by multiplicative constants (see \cite{Jo2},\cite{Jo3}).

The facts known for mappings from annuli to $\mathfrak{P}_n$ motivate the following two problems related to the restricted validity of Gromov's Oka principle.
\begin{probl}\label{probl1}
Fix a connected open (i.e. non-compact) complex manifold $X$ and a complex manifold $Y$. Obtain information about the set of continuous or smooth mappings $X\to Y$ that are homotopic to holomorphic mappings.
\end{probl}
\begin{probl}\label{probl2}

Consider a connected smooth oriented open manifold $X$, a complex manifold $Y$, and a smooth orientation preserving mapping $f:X\to Y$. Obtain information on the set of complex manifolds $X'$ which admit a smooth orientation preserving isomorphism $X'\to X$ for which the pull back of $f$ is homotopic to a holomorphic mapping.
\end{probl}
In the sequel all manifolds will be oriented, and all manifolds (including Riemann surfaces) will be connected, if not said otherwise.

In a forthcoming paper (see \cite{Jo4}) we address Problem \ref{probl1}. For instance, we give an upper bound for the number of homotopy classes of irreducible mappings from a finite open Riemann surface
to the thrice punctured Riemann sphere,
that contain a holomorphic mapping. A non-contractible mapping to the trice punctured Riemann sphere is called irreducible if it is not homotopic to a mapping into a punctured disc.
The estimates can be regarded as a quantitative information concerning the restricted validity of Gromov's Oka principle for some cases when the target manifold is not a Gromov-Oka manifold.

Here we will address Problem \ref{probl2}.
We restrict ourselves to surfaces $X$. Recall that
a surface is called finite if its fundamental group is finitely generated.
Each finite open Riemann surface $X$ is conformally equivalent
to a domain (denoted again by $X$) on a closed Riemann surface $X^c$ such that each connected component of the complement $X^c \setminus X$ is either a point or a closed topological disc with smooth boundary \cite{Sto}.
The connected components of the complement will be called holes.
A Riemann surface is called of first kind,
if it is closed, or it is obtained from a closed Riemann surface by removing finitely many points (called punctures). Otherwise the connected Riemann surface is called of second kind. If all holes of a finite open Riemann surface are closed topological discs, the Riemann surface is said to have only thick ends.

For a smooth oriented surface $X$ an orientation preserving homeomorphism $\omega:X\to \omega(X)$ from $X$ to a Riemann surface $\omega(X)$ is called a conformal structure on $X$. A mapping $f$ from $X$ to a complex manifold $Y$ is said to be holomorphic for the complex structure $\omega$ if the mapping $f\circ \omega^{-1}$ is holomorphic on $\omega(X)$.
\begin{defn}\label{defn8.2} Let $X$ be an oriented finite open smooth surface and $Y$ a
complex manifold.
A continuous mapping $f:X\to Y$ is said to be homotopic to a holomorphic mapping
for a conformal structure $\omega:X\to \omega(X)$ on $X$, if $f \circ w^{-1}$ is homotopic to a holomorphic mapping from $\omega(X)$ to $Y$.
If this holds for any conformal structure $\omega$ on $X$ with only thick ends (i.e. with $\omega(X)$ having only thick ends),
the mapping $f:X\to Y$ is said to have the Gromov-Oka property.
\end{defn}
 \index{conformal structure}

Let $n=2$. For an annulus $A$ each continuous mapping $f:A\to\mathfrak{P}_2$
has the Gromov-Oka property.
The following lemma holds for $n=3$.
\begin{lemm}\label{lem0}
A continuous mapping $f$ from an annulus $A$ to $\mathfrak{P}_3$ is homotopic to a holomorphic mapping for each conformal structure of second kind on $A$ if and only if it has the Gromov-Oka property. This happens if and only if the inequality $\mathcal{M}(\hat{b}_f)> \frac{\pi}{2}(\log\frac{3+\sqrt{5}}{2})^{-1}$ holds.
\end{lemm}

Notice that $\log\frac{3+\sqrt{5}}{2}$ is the smallest non-vanishing
entropy among $3$-braids.

On the other hand only very few mappings $f$ with infinite $\mathcal{M}(\hat{b}_f)$ are homotopic to a holomorphic mapping on an annulus of first kind (i.e. on the punctured plane). This fact motivated the definition of mappings with the Gromov-Oka property from a finite open smooth surface to a complex manifold. Further, the following question arises, which is part of Problem \ref{probl2}.

\smallskip

\noindent {\bf Problem 2a.}
{\it Given a finite open smooth surface $X$ and a complex manifold $Y$
with a set of generators whose conjugacy classes have infinite conformal module, which mappings $f:X\to Y$ have the Gromov-Oka property? In particular, which mappings $f:X\to  \mathfrak{P}^n$ have the Gromov-Oka property?}

\smallskip

We consider now a smooth surface $X$ of genus one with a hole and address the question which mappings $X\to\mathfrak{P}^n$ have the Gromov-Oka property.
A Riemann surface of genus $1$ will be called a torus
and a Riemann surface of genus $1$ with a hole will be called a torus with a hole.

Note first, that any continuous map from an oriented finite open smooth surface $X$ into $\mathfrak{P}^2$ has the Gromov-Oka property, moreover, any such map is homotopic to a holomorphic map for any conformal structure on $X$. (See the proofs of Lemma \ref{lem*} and Theorem \ref{thmEl.0} below).

Let ${\mathcal X}$ be a topological space and let $\hat e$
be a conjugacy class
of  elements of the fundamental group $\pi_1 ({\mathcal X} , x_0)$ of
${\mathcal X}$.
Choose a loop $\gamma$ in ${\mathcal X}$ that represents $\hat e$.
For a mapping $f:\mathcal{X}\to \mathfrak{P}^n$ the conjugacy class of braids represented by the restriction $f|\gamma$ depends only on $f$ and on $\hat e$ and is denoted by $\hat{f_e}$.
Notice that $\mathcal{M}(\hat{f}_e)=\infty$ iff the restriction of $f$ to an annulus representing $\hat e$ has the Gromov-Oka property.

We consider now mappings from a smooth orienetd manifold of genus one with a hole to $\mathfrak{P}^3$.
The following lemma states that the Gromov-Oka property of the mapping on the manifold
implies the Gromov-Oka property of its restrictions to certain annuli.
\begin{lemm}\label{lem0'}
Let $X$ be a smooth oriented manifold of genus one with a hole.
If $f:X\to \mathfrak{P}^3$ has the Gromov-Oka property then $\mathcal{M}(\hat{f_e})=\infty$ if $e$ is a primitive element of $\pi_1(X,x_0)$ or $e$ is the commutator of two elements that form a basis of
$\pi_1(X,x_0)$.
\end{lemm}
The lemma is an immediate consequence of the following more general lemma.
\begin{lemm}\label{lem0''}
Let $X$ be an oriented smooth finite open surface and let $e\in \pi_1(X,x_0)$ be an element whose conjugacy class can be represented
by a simple closed curve in $X$. Then for any $a>0$ there exists a conformal structure $\omega:X\to \omega(X)$ on $X$, such that the Riemann surface $\omega(X)$ has only thick ends and contains an annulus $A$ representing $\hat e$ of conformal module $m(A)>a$.
\end{lemm}

Vice versa, Theorem \ref{thm8.0} below
states the existence of a finite number of annuli in a smooth surface of genus $1$ with a hole $X$ with the following property. If for a separable quasipolynomial of degree $3$ the restriction to each of these annuli
has the Gromov-Oka property, then the quasipolynomial on $X$ has the Gromov-Oka property.
Notice that the conjugacy classes of the restriction of a mapping to curves representing a set of generators of the fundamental group of $X$ do not determine the mapping.

\begin{thm}\label{thm8.0}
Let $X$ be a smooth surface of genus one with a hole, and let $\mathcal{E}=\{e_1,e_2\}$ be a chosen set of generators of the fundamental group of $X$ with base point $x_0$. Denote by $\mathcal{E}_0$ the set $\{e_1,e_2, e_2 e_1^{-1}, e_2 e_1^{-2}, e_1 e_2 e_1^{-1} e_2^{-1}\}$.

Let $f$ be a separable quasipolynomial of degree $3$ on $X$. Suppose that for each $e\in \mathcal{E}_0$ the inequality $\mathcal{M}(\hat{f_e})>\frac{\pi}{2}(\log\frac{3+\sqrt{5}}{2})^{-1}$ holds.
Then the quasipolynomial is isotopic to a separable algebroid function for any conformal structure of second kind on $X$.
\end{thm}

Notice that a Riemann surface with one hole is of second kind iff it has only thick ends.

The condition that a mapping $X\to \mathfrak{P}^3$ has the Gromov-Oka property is rather restrictive (see Theorem \ref{thm8.1} below).

We do not know whether the respective result holds for quasipolynomials on arbitrary smooth oriented finite open surfaces. The following theorem says that it holds for orientation preserving mappings from arbitrary smooth oriented finite open surfaces to the twice punctured complex plane (and, hence, for orientation preserving mappings from arbitrary smooth oriented finite open surfaces into $\mathfrak{P}_3$ whose
monodromies correspond to pure $3$-braids, equivalently, for mappings into $C_3(\mathbb{C})$).

\begin{thm}\label{prop0}
Let $X$ be a smooth oriented surface of genus $g\geq 0$ with $m\geq 1$ holes.
There exists a subset $\mathcal{E}'$ of the fundamental group $\pi_1(X,x_0)$ of $X$ consisting of at most $(2g+m-1)^3$ elements such that the conjugacy class $\hat e$ of each element of $\mathcal{E}'$ can be represented by
a simple closed curve and the following holds.

A smooth orientation preserving mapping $f:X\to \mathbb{C}\setminus \{-1,1\}$ has the Gromov-Oka property if and only if for each $e\in \mathcal{E}'$ the restriction of the mapping to an annulus in X that represents $\hat e$ has the Gromov-Oka property.
\end{thm}

For $g>0$ the statement holds for all smooth mappings (not only for orientation preserving ones).
Notice that $\pi_1(X,x_0)$ is a free group in $2g+m-1$ generators.
The mappings $X\to \mathbb{C}\setminus \{-1,1\}$ that enjoy the Gromov-Oka property
are completely described in Section \ref{sec:9.4a}.
The case of mappings from open Riemann surfaces to $\mathfrak{P}^n$ with $n>3$ is widely unknown even in the case when $X$ is a surface of genus one with a hole (see the discussion at the end of Section \ref{sec:8.3}).

We consider now an analog of Theorem \ref{thm8.0} for fiber bundles.
Note that each smooth bundle over a smooth surface of genus one with a hole with each fiber being a closed oriented surface of genus one is isotopic to a smooth bundle with fibers equipped with complex structures depending smoothly on the point in the base manifold (see Section \ref{sec:9.7}). We call the latter bundles torus bundles. Theorem \ref{thm8.0} has an analog concerning isotopic deformations of smooth torus bundles
over a Riemann surface of genus one with a hole to holomorphic bundles.
A smooth torus bundle $\mathfrak{F}$ over a smooth finite open surface $X$ is said to have the Gromov-Oka property, if for each smooth conformal structure $\omega:X\to \omega(X)$ on $X$ with only thick ends the push-forward $\mathfrak{F}_{\omega}$ to $\omega(X)$ is isotopic to a holomorphic torus bundle.

\begin{thm}\label{thmEl.9a}
Let $X$ be a smooth surface of genus one with a hole with base point $x_0$, and with a chosen set
$\mathcal{E}=\{e_1,e_2\}$ of generators of $\pi_1(X,x_0)$. Define the set
$\mathcal{E}_0=
\{e_1,e_2, e_1 e_2^{-1},  e_1 e_2^{-2}, e_1e_2 e_1^{-1} e_2^{-1}\}$ as in Theorem {\rm \ref{thm8.0}}. Consider a smooth bundle $\mathfrak{F}$
over $X$ whose fibers are closed Riemann surfaces of genus one. Suppose for each $e\in \mathcal{E}_0$ the restriction of the bundle $\mathfrak{F}$
to an annulus
in $X$ representing $\hat e$ has the Gromov-Oka property.
Then the bundle over $X$ has the Gromov-Oka property.
\end{thm}

A similar statement (see Theorem \ref{thmEl.0} below) is proved for fibers bundles each fiber of which is the complex plane with a set of three distinguished points. This statement is closely related to Theorem
\ref{thm8.0} on quasipolynomials. Theorem \ref{thmEl.9a} is proved by representing the torus bundles as double branched coverings of bundles with fibers equal to $\mathbb{C}$ with a set of three distinguished points.
In section \ref{sec:9.7} we prove a more comprehensive theorem on torus bundles.

\medskip

\noindent \textbf{Acknowledgement}.
The author is grateful to the IHES, the CRM Barcelona, and the Indiana University Bloomington, where a big part of the work on the paper has been done.
A talk given at the Algebraic
Geometry seminar at Courant Institute and a stimulating discussion
with F.Bogomolov and M.Gromov had a special impact. The author is also
indebted to many other mathematicians for interesting and helpful
discussions and for information concerning references to the
literature, among them R.Bryant, D.Calegari, and
M.Zaidenberg. The figures were drawn by M.Vergne and F.Dufour.

\section{ Families of polynomials, solvability and reducibility.}
\label{sec:8.2}

In this section we will give the proof of Theorem \ref{thm8.3} on reducibility of holomorphic quasipolynomials, and present some preliminary results on braids and mapping classes, that will be needed later.

Notice first, that the isomorphism between $\pi_1(\,C_n(\mathbb {C}) \diagup {\mathcal S}_n\,,\; E^0_n\,)$ with $E^0_n=\{0,\ldots,\frac{n-1}{n}\}$ and the Artin braid group with generators $\sigma_j, \, j=1,\ldots,n-1,$ is given as follows. Take any arc $\gamma:[0,1]\to C_n(\mathbb {C}) \diagup {\mathcal S}_n$ with endpoints equal to $E^0_n$ that represents an element of the fundamental group $\pi_1(\,C_n(\mathbb {C}) \diagup {\mathcal S}_n\,,\; E^0_n\,)$. Associate to the arc $\gamma$ its lift  $\tilde {\gamma}=(\tilde {\gamma}_1,\ldots,\tilde {\gamma}_n): [0,1]\to C_n(\mathbb {C})$ for which $\tilde {\gamma}_j(0)= \frac{j-1}{n},\,j=1,\ldots,n$. For $j=1,\ldots,n-1,$ the isomorphism maps the generator $\sigma_j$ to  the homotopy class of  arcs in  $C_n(\mathbb {C})\diagup \mathcal{S}_n$ with endpoints equal to $E_n$, whose associated lift is represented by the following mapping. The mapping moves the points $\tilde {\gamma}_{j}(0)$ and $\tilde {\gamma}_{j+1}(0)$ with uniform speed in positive direction along a half-circle with center at the point $\frac{1}{2}(\tilde {\gamma}_{j}(0)+\tilde {\gamma}_{j+1}(0))$, and
fixes $\tilde {\gamma}_{j'}(0)$ for all $j'\neq j,j+1$.
For more details on braids see \cite{Bi} and \cite{KaTu}.
\index{$\mathcal{B}_n$}

In the sequel
we will need the relation between braids and mapping classes which we recall now.
Let $A$
be an oriented manifold.
Let $A_1$ and $A_2$ be disjoint closed subsets of $A$. Denote by
${\rm Hom} ^+(A ; A_1 , A_2)$ \index{${\rm Hom}^+ (A ; A_1 , A_2)$} the
set of orientation preserving self-homeomorphisms of $A$ that fix $A_1$ pointwise and $A_2$
setwise.
Equip the
set ${\rm Hom}^+ (A ; A_1 , A_2)$ with the compact open topology. The set of connected components
of ${\rm Hom}^+ (A;A_1,A_2)$ is denoted by ${\mathfrak M}
(A;A_1,A_2)$. \index{$\mathfrak{M} (A;A_1,A_2)$} It carries the structure of a group under composition, and is called a mapping class group.

Let $\mathbb{D}$ be the unit disc in the complex plane
and let
$E_n\subset \mathbb{D}$ be a set containing $n$ points.
The braid group ${\mathcal
B}_n$ is isomorphic to $\mathfrak{M}(\overline{\mathbb{D}};\partial \mathbb{D}, E_n)$. The set $\mathfrak{M}(\overline{\mathbb{D}};\partial \mathbb{D}, E_n)$ is called the mapping class group of the closed disc with set $E_n$ of $n$ distinguished points. For each generator $\sigma_j$ of the braid group we consider the element $\mathfrak{m}_{\sigma_j}\in \mathfrak{M}(\overline{\mathbb{D}};\partial \mathbb{D}, E^0_n)$
which is represented by a homeomorphism $\varphi_j\in {\rm Hom} ^+(\overline{\mathbb{D}};\partial \mathbb{D}, E^0_n )$ that is the identity outside the disc with center $\frac{1}{2}(\frac{j-1}{n}+ \frac{j}{n})$ and radius $\frac{1}{n}$, equals rotation by the angle $\pi$ on the disc with center $\frac{1}{2}(\frac{j-1}{n}+ \frac{j}{n})$ and radius $\frac{1}{2n}$, and rotates each point in the remaining annulus by an angle between $0$ and $\pi$. The mapping ${\mathcal
B}_n \ni \sigma_j\to \mathfrak{m}_{\sigma_j}^{-1}$ extends to a homomorphism from  ${\mathcal
B}_n $ onto $\mathfrak{M}(\overline{\mathbb{D}};\partial \mathbb{D}, E^0_n)$.

Each self-homeomorphism $\varphi$ of $\overline{\mathbb{D}}$ that represents the mapping class $\mathfrak{m}_b$ of a braid $b$
extends to a self-homeomorphism of the Riemann sphere $\mathbb{P}^1$ by taking it equal to
the identity outside $\overline{\mathbb{D}}$. The extended homeomorphism represents an element of $\mathfrak{M}(\mathbb{P}^1;\infty, E_n)$, which is denoted by $\mathcal{H}_{\infty}(m_b)$.
The group $\mathfrak{M}(\mathbb{P}^1;\infty, E_n)$ is isomorphic to the braid group modulo its center $\mathcal{B}_n\diagup \mathcal{Z}_n$. A braid $b$ is called periodic if the mapping class $\mathcal{H}_{\infty}(\mathfrak{m}_b)$ contains a periodic mapping.

We explain now Thurston's notion of irreducible surface homeomorphisms and irreducible braids.
Let $S$ be a finite smooth (oriented) surface. It is either closed or homeomorphic to a closed surface with a finite number of punctures. We will assume from the beginning that $S$ is either closed or punctured.

A finite non-empty set of mutually disjoint simple closed curves $\{ C_1 ,
\ldots , C_{\alpha}\}$ on a closed or punctured surface $S$ is called
admissible if no $C_i$ is homotopic
to a point in $X$, or to a puncture,
or
to a $C_j$ with $i \ne j$. Thurston calls an isotopy class $\mathfrak{m}$ of self-homeomorphisms
of $S$ (in other words, a mapping class on $S$) reducible if there is an admissible system
of curves $\{ C_1 ,\ldots , C_{\alpha}\}$ on $S$ such that
some (and, hence, each) element in $\mathfrak{m}$ maps the system to an isotopic system. In
this case we say that the system  $\{ C_1 ,
\ldots , C_{\alpha}\}$ reduces $\mathfrak{m}$. A mapping class which is not reducible is
called irreducible. Mappings are called reducible if they belong to a reducible mapping class, and are called irreducible otherwise. The mapping class of the identity is reducible.

Let $S$ be a closed or punctured surface with set $E$ of distinguished points. We say that $\varphi$ is a self-homeomorphism of $S$ with distinguished points $E$, if $\varphi$ is a self-homeomorphism of $S$ that maps the set of distinguished points $E$ to itself.
Notice that each self-homeomorphism of the punctured surface $S\setminus E$ extends to a self-homeomorphism of the surface $S$ with set of distinguished points $E$.
We will sometimes identify
self-homeomorphisms of $S\setminus E$ and self-homeomorphism of $S$ with set $E$ of distinguished points.

For a finite (closed or punctured) surface $S$ and a finite subset $E$ of $S$ a finite non-empty set of mutually disjoint simple closed curves $\{ C_1 ,
\ldots , C_{\alpha}\}$ in $S\setminus E$ is called
admissible for $S$ with set of distinguished points $E$ if
it is admissible for $S \setminus E$. An admissible system of
curves for $S$ with set of distinguished points $E$ is said to reduce a mapping class $\mathfrak{m}$ on $S$ with set of distinguished
points $E$, if the induced mapping class on $S\setminus E$ is reduced by this system
of curves.

Conjugacy classes of reducible mapping classes can be decomposed in some sense into
irreducible components, and conjugacy classes of reducible mapping classes can be recovered
from the irreducible components up to products of commuting Dehn twists. Conjugacy classes
of irreducible mapping classes are classified and studied.

A Dehn twist \index{Dehn twist} about a simple closed curve $\gamma$
in an oriented surface $S$ is a mapping that is isotopic to the following one. Take a tubular neighbourhood of
$\gamma$ and parameterize it as a round annulus  $A=\{ e^{-\varepsilon} < \vert z \vert < 1\}$  so that $\gamma$ corresponds to $|z|=e^{-\frac{\varepsilon}{2}}$. The mapping is
an orientation preserving self-homeomorphism of $S$
which is the identity outside $A$ and is equal to the mapping $e^{-\varepsilon s +2\pi i t}\to  e^{-\varepsilon s +2\pi i (t+s) }$ for  $e^{-\varepsilon s +2\pi i t}\in A$, i.e. $s\in (0,1)$. Here $\varepsilon$ is a small positive number.

A braid $b\in \mathcal{B}_n$ is called irreducible if the associated mapping class $\mathfrak{m}_b\in \mathfrak{M}(\overline{\mathbb{D}}; \partial \mathbb{D},E_n)$ is irreducible, and reducible otherwise.

The following simple Lemma concerns global solvability of holomorphic families of polynomials.


\begin{lemm}\label{lem8.0}  Let $X$ be a closed Riemann surface
of positive genus with a closed smoothly bounded topological disc removed. Suppose $f$ is an
irreducible separable algebroid function of degree $3$ on $X$.
Suppose $X$ contains a domain $A$, one of whose boundary components
coincides with the boundary circle of $X$, such that $A$ is
conformally equivalent to an annulus of conformal module strictly
larger than $\frac{\pi}{2}\, (\log( \frac{3 +\sqrt{5}}{2}))^{-1}$.
Then $f$ is solvable over $A$.
\end{lemm}

We postpone the proof of Lemma \ref{lem8.0}, and state some simple lemmas which will be needed in the sequel. Let $A$ be an annulus.
For a continuous mapping $f$ from an annulus $A$ to the symmetrized configuration space $C_n(\mathbb{C})\diagup \mathcal{S}_n$
we denote as before by $\hat b_{f}$ the conjugacy class of braids
that is represented by the mapping $f$.

\begin{lemm}\label{lem8.1} Suppose the separable quasipolynomial $f$ of degree
$n$ on an annulus $A$ is irreducible and $n$ is prime. Then the
represented conjugacy class of braids $\hat b_{f}$ is irreducible.
\end{lemm}
Notice that, on the other hand, conjugacy classes of irreducible
pure braids define solvable, hence reducible, quasipolynomials on
the circle.

\smallskip

\noindent {\bf Proof.} Recall that $\tau_n:\mathcal{B}_n\to \mathcal{S}_n$ is the natural projection from the braid group to the symmetric group. $\tau_n$ maps conjugacy classes of braids to conjugacy classes of permutations. If $f$ is an irreducible separable quasipolynomial, the conjugacy class of
braids $\hat b_{f}$ projects to a conjugacy class $\tau_n(\hat b_{f})$ of
$n$-cycles. The lemma is now a consequence of the following known
lemma (see e.g. \cite{FaMa}). \hfill $\Box$
\begin{lemm}\label{lemm8.2} If $n$ is prime then any braid $b \in {\mathcal
B}_n$, for which $\tau_n (b)$ is an $n$-cycle, is irreducible.
\end{lemm}
For convenience of the reader we give the short argument.\\

\noindent {\bf Proof of Lemma \ref{lemm8.2}.} If $b$ was reducible
then some homeomorphism $\varphi$ which represents the mapping class
corresponding to $b$ would fix setwise an admissible system of
curves ${\mathcal C}$. Let $C_1$ be one of the curves in ${\mathcal
C}$ and let $\delta_1$ be the topological disc contained in
${\mathbb D}$, bounded by $C_1$. $\delta_1$ contains at least two
distinguished points $z_1$ and $z_2$. Since $\varphi$ permutes the
distinguished points along an $n$-cycle there is a power
${\varphi}^k$ of $\varphi$ which maps $z_1$ to $z_2$. Since $n$ is
prime, $\varphi^k$ also permutes the distinguished points along an
$n$-cycle. Hence it maps some distinguished point in $\delta_1$ to
the complement of $\delta_1$. We obtained that $\varphi^k
(\delta_1)$ intersects both, $\delta_1$ and its complement. Hence
$\varphi^k (C_1) \ne C_1$ but $\varphi^k (C_1)$ intersects $C_1$.
This contradicts the fact that the system of curves ${\mathcal C}$
was admissible and invariant under $\varphi$. \hfill $\Box$

\medskip

Let $\Delta_3=\sigma_1\sigma_2\sigma_1= \sigma_2\sigma_1\sigma_2$. The group $\langle \Delta^2\rangle$ generated by $\Delta^2$ coincides with the center $\mathcal{Z}_3$ of the braid group $\mathcal{B}_3$.  We need the following lemma.
\begin{lemm}\label{lemEl.0}
A mapping class in $\mathfrak{M}(\overline{\mathbb{D}}; \partial \mathbb{D}, E_3)$ is reducible if and only if the braid corresponding to it
is conjugate to $\sigma_1^k\Delta_3^{2\ell}$ for integers $k$ and $\ell$.

A mapping class
in $\mathfrak{M}(\mathbb{P}^1; \{\infty\}, E_3)$ is reducible if and only if it corresponds to a conjugate of $\sigma_1^k\diagup \mathcal{Z}_3$.
\end{lemm}
\noindent {\bf Proof.} Any admissible system of curves on $\overline{{\mathbb{D}}}$ with three distinguished points consists of one curve that divides $\overline{\mathbb{D}}$ into two connected components.
One of these components contains $\partial \mathbb{D}$ and one of the distinguished points, the other one contains two distinguished points.
For a reducible mapping class $\mathfrak{m}$ in $\mathfrak{M}(\overline{\mathbb{D}}; \partial \mathbb{D}, E_3)$ we take an admissible curve $\gamma$ that reduces  $\mathfrak{m}$ and a representing homeomorphism $\varphi\in \mathfrak{m}$ that fixes $\gamma$ pointwise. Suppose the points in $E_3$ are labelled according to a homomorphism from
$\pi_1(C_3(\mathbb{C})\diagup \mathcal{S}_3, E_n)$ to $\pi_1(C_3(\mathbb{C})\diagup \mathcal{S}_3, E^0_n)$ (which is defined up to conjugation)
so that one of the connected components $\mathcal{C}_1$ of the complement of $\gamma$ contains $\zeta_1$ and $\zeta_2$, the other connected component $\mathcal{C}_2$ contains $\zeta_3$ and $\partial \mathbb{D}$.
Since $\varphi$ fixes $\partial \mathbb{D}$ and $\gamma$ pointwise, it also fixes the annulus $\overline{\mathcal{C}}_2$ setwise and fixes the point $\zeta_3$.
A self-homeomorphism of an annulus that fixes the boundary pointwise is a twist.
After an isotopy of $\varphi$, that fixes $\gamma$ setwise and fixes each distinguished point, we may assume that $\varphi$ fixes pointwise a simple arc that joins $\zeta_3$ with $\gamma$. Then
the restriction $\varphi\mid \overline{\mathcal{C}}_2$ represents a power of a Dehn twist in $\mathbb{D}$ about a curve in $\mathcal{C}_2$ that is homotopic in ${\overline{\mathcal{C}}}_2\setminus \{\zeta_3\}$ to $\partial \mathbb{D}$.
Such a Dehn twist corresponds to the braid $\Delta_3^2$.
The restriction $\varphi\mid \overline{\mathcal{C}}_1$ represents the class in $\mathfrak{M}(\overline{\mathcal{C}}_1;\partial \mathcal{C}_1, \{\zeta_1,\zeta_2\})$, corresponding to $\sigma^k$ for some integer $k$.
Hence the class represented by $\varphi$ in
$\mathfrak{M}(\overline{\mathbb{D}}; \partial \mathbb{D}, E_3)$ corresponds to a braid that is conjugate to $\sigma_1^k \Delta_3^{2\ell}$ for some integer $\ell$. Vice versa, it is clear that the mapping class $\mathfrak{m}_{\sigma_1^{2k} \Delta_3^{2\ell}}$ is reducible.

The statement concerning the reducible elements of $\mathfrak{M}(\mathbb{P}^1; \{\infty\}, E_3)$ is proved in the same way.
\hfill $\Box$

\medskip

It is known that the smallest non-vanishing entropy among
irreducible $3$-braids equals $\log \frac{3+\sqrt 5}2$ (see e.g.
\cite{Son}).  By the following Theorem of \cite{Jo} the number $\frac\pi2 \left(\log \frac{3+\sqrt 5}2 \right)^{-1}$ is the largest finite conformal module
among irreducible $3$-braids.

\begin{thm}\label{thm*}
For each braid $b\in \mathcal{B}_n$, $n\geq 2$, the following relation holds
\begin{equation}\label{eq22}
\mathcal{M}(\hat{ b})=\frac{\pi}{2}\frac{1}{h(\hat{ b})}
\end{equation}
for the entropy $h(\hat{ b})$ and the conformal module $\mathcal{M}(\hat{b})$ of its conjugacy class $\hat b$.
\end{thm}

The following lemma holds.

\begin{lemm}\label{lemm8.3} Suppose for a conjugacy class of braids $\hat b
\in \hat{\mathcal B}_3$ the conformal module ${\mathcal M}
(\hat b)$ satisfies the inequality ${\mathcal M} (\hat b) > \frac\pi2 \left(\log \frac{3+\sqrt 5}2 \right)^{-1}$. Then the following holds.
\begin{enumerate}
\item[(a)] ${\mathcal M} (\hat b) = \infty$.
\item[(b)] {\it If $\hat b$ is irreducible then $\hat b$ is the conjugacy class of a
periodic braid, i.e. either of $(\sigma_1 \, \sigma_2)^{\ell}$ for
an integer $\ell$ not divisible by $3$, or of $(\sigma_1 \, \sigma_2
\, \sigma_1)^{\ell}$ for an integer $\ell$ not divisible by $2$.}
\item[(c)] {\it If $\hat b$ is reducible then $\hat b$ is the conjugacy class of $\sigma_1^k \, \Delta_3^{2\ell}$ for integers $k$ and $\ell$.}
    \item[(d)] {\it If $\hat b$ is the conjugacy class of a commutator, then it is represented by a pure braid $b$.}
\end{enumerate}
\end{lemm}

\noindent {\bf Proof.} Suppose $\hat b$ is irreducible. Since $\frac\pi2 \left(\log \frac{3+\sqrt 5}2 \right)^{-1}$ is the largest finite conformal module
among irreducible $3$-braids the equality ${\mathcal M} (\hat b) = \infty$ holds. (Hence, also  $h(\hat b) = 0$.) This proves (a) in
the irreducible case.

Let $b \in {\mathcal B}_3$ represent the irreducible class $\hat
b$, let ${\mathfrak m}_b\in \mathfrak{M}(\overline{\mathbb{D}};\partial\mathbb{D},E_3)$ be the mapping class corresponding to $b$, and let ${\mathcal H}_{\infty} ({\mathfrak
m}_b)\in \mathfrak{M}(\mathbb{P}^1;\infty,E_3)$ be the corresponding mapping class on $\mathbb{P}^1$ with distinguished points. Since $\mathcal{M}(\hat{b})=\infty$, the class ${\mathcal H}_{\infty} ({\mathfrak
m}_b)$ is represented by a periodic mapping of the
complex plane ${\mathbb C}$ with three distinguished points (see \cite{Be1}). By a conjugation we may assume that the periodic mapping fixes zero and
hence is a rotation by a root of unity. If zero is not a
distinguished point, the three distinguished points are
equidistributed on a circle with center zero and the mapping is
rotation by a power of $e^{\frac{2\pi i}3}$. This mapping
corresponds to
$$
\reallywidehat{(\sigma_1 \, \sigma_2)^{\ell} / \langle \Delta_3^2 \rangle}
\, .
$$
If zero is a distinguished point, the other two distinguished points
are equidistributed on a circle with center zero. The mapping is a
rotation, precisely, a multiplication by a power of $-1$, and
corresponds to $\reallywidehat{(\sigma_1 \, \sigma_2 \, \sigma_1)^{\ell}
/ \langle \Delta_3^2 \rangle}$.
(Note that powers of $\Delta_3^2$ are reducible and the mapping class in $\mathfrak{M}(\mathbb{P}^1;\{\infty\}, E_3)$ corresponding to $\Delta_3^2$ is the identity.) We proved (b).

Suppose $\hat b$ is reducible, i.e. the mapping class
${\mathfrak m}_b$ is reducible for $b \in \hat b$.
Then by Lemma \ref{lemEl.0} the class $\hat b$ is the conjugacy class of
$\sigma_1^k\Delta_3^{2\ell}$ for integers $k$ and $\ell$. We obtained (c). The braid $b =
\sigma_1^k \, \Delta_3^{2\ell}$ has infinite conformal module. This
gives (a) in the reducible case.

Suppose $\hat b$ is the conjugacy class of a commutator $b \in
{\mathcal B}_3$. Then $b$ has exponent sum zero (more detailed,
the sum of exponents of the generators in a representing word equals zero). If ${\mathcal M}
(\hat b) = \infty$ then this is possible only if $b$ is conjugate to $\sigma_1^k
\, \Delta_3^{2\ell}$
with $k+6\ell=0$. Hence, $k$ is even and the
braid $b$ is pure. This proves (d). The lemma is proved. \hfill
$\Box$

\medskip
\noindent{\bf Proof of Lemma \ref{lem0}.} If $f$ is homotopic to a holomorphic function for any conformal structure of second kind on $A$ or for any conformal structure with only thick ends, then $\mathcal{M}(\hat{b_f})=\infty$.
Suppose on the other hand that $\mathcal{M}(\hat{b_f})> \frac{\pi}{2}\, (\log( \frac{3 +\sqrt{5}}{2}))^{-1}$.
Then Lemma \ref{lemEl.0} implies that
$\mathcal{M}(\hat{b_f})=\infty$. Hence $f$ is homotopic to a holomorphic function for any conformal structure on $A$ of finite conformal module. There are two conformal structures of infinite conformal module, one conformally equivalent to the punctured disc $\mathbb{D}\setminus \{0\}$, it is of second kind, and the other conformally equivalent to the punctured complex plane, it is of first kind. By Lemma \ref{lemEl.0} a conjugacy class of braids of infinite conformal module is either a periodic braid, or it is conjugate to $\sigma_1^k$. We claim that in the first case $f$ is isotopic on the punctured plane to the quasipolynomial
$f_{1,k}(z,\zeta)= z^{k}-\zeta^3, \, z\in \mathbb{C}\setminus\{0\}, \,\zeta\in \mathbb{C}$,
or to $f_{2,k}(z,\zeta)= \zeta (z^k-\zeta^2), \, z\in \mathbb{C}\setminus\{0\}, \,\zeta\in \mathbb{C}$, respectively.
In the second case we take any constant $R>1$. We claim that the quasipolynomial $f$ is isotopic on the punctured disc to $f_{3,k}(z,\zeta)= (z^k-\zeta^2) (\zeta -R)   , \, z\in \mathbb{D}\setminus\{0\}, \,\zeta\in \mathbb{C}$.
Indeed, for $z = \frac{1}{2} e^{2\pi i
t}$, $t \in [0,1]$, the set of solutions of the equation $f_{1,k}
(z,\zeta) = 0$ is
$$
E_1^{k}(t) \stackrel{def}= 2^{-\frac{1}{3}}
e^{\frac{ 2\pi ik}{3} t}\left\{1 , e^{\frac{2\pi ik}{3} } , e^{\frac{4\pi ik}{3} }\right\}  \, , \quad t \in [0,1] \, .
$$
The path $t\to E_1^{k}(t),\,t \in [0,1],$ in ${\mathfrak P}_3$ defines a geometric braid in the
conjugacy class of $(\sigma_1 \cdot \sigma_2)^{k} $.
The set of solutions of the equation $f_{2,k}
(z,\zeta) = 0$ is
$$
E_2^k(t) \stackrel{def}= 
2^{-\frac{1}{2}}   e^{\frac{ 2\pi ik}{2} t}\left\{-1,0,1\emph{}\right\}  \, , \quad t \in [0,1] \, .
$$
The respective path represents $\Delta_3^k$.
Finally, the set of solutions of the equation $f_{3,k}
(z,\zeta) = 0$ is
$$
E_3^k(t) \stackrel{def}=
  \left\{2^{-\frac{1}{2}} e^{\frac{ 2\pi ik}{2} t}, - 2^{-\frac{1}{2}}e^{\frac{ 2\pi ik}{2} t} ,R\right\}  \, , \quad t \in [0,1] \, .
$$
The respective path represents $\sigma_1^k$.

The quasipolynomial $f_3$ is not
isotopic to a holomorphic one on the punctured plane. This follows from statement $4$
of Theorem \ref{thmEl.0} below.
\hfill $\Box$

\medskip
\noindent{\bf Proof of Lemma \ref{lem8.0}.} Let  $\mathfrak{S}_f= \{(z,\zeta)\in X\times \mathbb{C}: f(z,\zeta)=0\}$ be the zero set of the algebroid function $f$. Since $f$ is a separable irreducible quasipolynomial, this set is connected, and we obtain
an unramified covering $\mathfrak{S}_f \rightarrow X$ of degree $3$. Hence, for the Euler characteristics $\chi(X)$ and $\chi(\mathfrak{S}_f)$ the relation
$\chi(\mathfrak{S}_f)=3\chi(X)$ holds. Let $m(X)$ and $m(\mathfrak{S}_f)$ be the number of boundary components of $X$ and $\mathfrak{S}_f$. Then
$2-2g(\mathfrak{S}_f)-m(\mathfrak{S}_f)= 3(2-2g(X)-m(X))$. Hence, since $m(X)=1$, $m(\mathfrak{S}_f)$ equals either $1$ or $3$. This means that
$f\mid A$ is either irreducible or solvable. If $f\mid A$ is irreducible, it represents an irreducible conjugacy class of braids  $\hat b_{f}$ \index{$\hat
b_{f}$}. By Lemma \ref{lemm8.3}  $\hat b_{f}$ must be the conjugacy class of a  periodic braid which corresponds to a $3$-cycle. This is impossible for conjugacy classes of products of braid commutators. \hfill $\Box$

\medskip

We will prove now Theorem \ref{thm8.3} using the relation between the conformal module and the entropy of conjugacy classes of braids (see Theorem 1 of \cite{Jo})
and the following result of Penner \cite{P} on the smallest non-vanishing
entropy among irreducible $n$-braids.

\begin{thm}\label{thm8.4}{\rm (Penner)} For $3{\sf{g}}-3+{\sf{m}} > 0$ we denote by $h_{\sf{g}}^{\sf{m}}$ the
smallest non-vanishing entropy among irreducible self-homeomorphisms
of Riemann surfaces of genus $\sf{g}$ with $\sf{m}$ distinguished points. Then
$$
h_{\sf g}^{\sf m} \geq \frac{\log 2}{12{\sf g} - 12+4 {\sf m}} \, .
$$
\end{thm}
The entropy is a conjugacy invariant, i.e. the entropy of conjugate mappings is equal. The entropy of a braid $b\in \mathcal{B}_n$ is equal to the entropy of the mapping class $\mathcal{H}(\mathfrak{m}_b)$ corresponding to $b$.
Hence, the smallest
non-vanishing entropy among irreducible $n$-braids, $n \geq 3$, is
bounded from below by
$$
\frac{\log 2}{4(n+1)-12} = \frac{\log 2}{4n-8} \geq \frac{\log 2}{4}
\cdot \frac 1n \, , \quad n \geq 3 \, .
$$
By equality \eqref{eq22}
the largest finite conformal module among irreducible conjugacy classes of $n$-braids does not exceed $\frac{\pi}{2}\frac{4}{\log 2}
\cdot n$  for $ n \geq 3$.

Recall that the discriminant ${\sf D}_n$ is a function on
the space of all polynomials $\overline{\mathfrak P}_n$ of degree
$n$ which vanishes exactly on the set of polynomials with multiple
zeros. Explicitly, for a polynomial ${\sf p} \in \overline{\mathfrak
P}_n$, ${\sf p} (\zeta) = \underset{j=1}{\overset{n}{\prod}} (\zeta
- \zeta_j)$, its discriminant equals
$$
{\sf D}_n ({\sf p}) = \prod_{i < j} (\zeta_i - \zeta_j)^2 \, .
$$

\medskip

\noindent {\bf Proof of Theorem \ref{thm8.3}.} Suppose $f$ is an irreducible
separable algebroid function on the annulus $A$ and
\begin{equation}\label{eq8.19'}
m(A) > \frac{\pi}{2} \ \frac{4}{\log2} \cdot n = \frac{2\pi}{\log 2}
\cdot n \, .
\end{equation}
Since $n$ is prime, by Lemma \ref{lem8.1} the conjugacy class $\hat b_{f} \in \hat{\mathcal
B}_3$ induced by $f$ is irreducible.
By inequality \eqref{eq8.19'} ${\mathcal M} (\hat b_{f}) = \infty$.

This implies that $\hat b_{f}$ is the conjugacy class of a
periodic braid corresponding to an $n$-cycle, hence it is the
conjugacy class of $(\sigma_1 \, \sigma_2 \ldots \sigma_n)^k$ for an
integer $k$ which is not divisible by $n$. The isotopy class of the
algebroid function $\tilde f (z,\zeta) = \zeta^n - z^k$, $z \in A$,
induces this conjugacy class of braids. In other words, $\tilde f$ is
isotopic to $f$.  Indeed, for $z = e^{2\pi i
t}$, $t \in [0,1]$, the set of solutions of the equation $\tilde f
(z,\zeta) = 0$ is
$$
E_n(t) \stackrel{def}=\left\{ e^{\frac{2\pi ik}{n} t} , e^{\frac{2\pi ik}n + \frac{2\pi i
k}n t} , \ldots , e^{\frac{2\pi i k (n-1)}n + \frac{2\pi i k}n t}
\right\} \, , \quad t \in [0,1] \, .
$$
The path $t\to E_n(t),\,t \in [0,1],$ in ${\mathfrak P}_n$ defines a geometric braid in the
conjugacy class of $(\sigma_1 \cdot \sigma_2 \cdot \ldots \cdot
\sigma_n)^k$.

Compute the discriminant ${\sf D}_n (\tilde f_z)$:
\begin{equation}
{\sf D}_n (\tilde f_z) \,= \,\prod_{0 \leq m < \ell < n} \left(
e^{\frac{2\pi ik}{n} \cdot m + \frac{2\pi i k}n \cdot t} -
e^{\frac{2\pi ik}{n} \cdot \ell + \frac{2\pi i k}n \cdot t}
\right)^2 \nonumber
\end{equation}
\begin{equation}\label{eq8.19} = e^{\frac{2\pi ik}{n} \cdot t \cdot n \cdot (n-1)}
\cdot \prod_{0 \leq m < \ell < n} \left( e^{\frac{2\pi ik}{n} m}
- e^{\frac{2\pi ik}{n} \ell} \right)^2 \, . 
\end{equation}
Hence
$$
{\sf D}_n (\tilde f_z) = e^{2\pi i k \cdot (n-1) \cdot t} \cdot c_n
\, , \quad z = e^{2\pi i t} \, , \quad t \in [0,1] \, ,
$$
where $c_n$ is a non-zero constant depending only on $n$. (It is
equal to the product in the last expression of \eqref{eq8.19}.) The index of
the discriminant equals $k \cdot (n-1)$ which is not divisible by
$n$. In other words, the condition for the discriminant excludes the
only possibility for the isotopy class of a separable algebroid function on annuli of the
given conformal module to be irreducible. Hence under the conditions
of Theorem \ref{thm8.3} the algebroid function $f$ must be reducible. \hfill $\Box$

\section{ The conformal module and isotopy of quasipolynomials to algebroid functions. }
\label{sec:8.3}

Two homomorphisms $h_1$ and $h_2$ from a group $G_1$ to a group $G_2$ are called conjugate if for an element $x_2\in G_2$ the equality $h_2(x)=x_2^{-1}h_1(x)x_2$ holds for all $x\in G_1$.

For two smooth open surfaces $\mathcal{X}$ and $\mathcal{Y}$ with base points $x_0 \in \mathcal{X}$ and $y_0 \in \mathcal{Y}$ and a continuous mapping $F:\mathcal{X} \to \mathcal{Y}$ with $F(x_0)=y_0$ we denote by $F_*: \pi_1(\mathcal{X},x_0) \to \pi_1(\mathcal{Y},y_0) $ the induced map on fundamental groups. For each element $e_0\in \pi_1(\mathcal{X},x_0)$ the image $F_*(e_0)$ is called the monodromy along $e_0$, and the homomorphism $F_*$ is called the monodromy homomorphism corresponding to $F$.
The homomorphism $F_*$ defines a one to one correspondence from
the set of homotopy classes of continuous mappings $\mathcal{X}\to \mathcal{Y}$ with fixed base point $x_0$ in the source and fixed value $y_0$ at the base point to the set of homomorphisms $\pi_1(\mathcal{X},x_0) \to \pi_1(\mathcal{Y},y_0) $.

Consider a free homotopy $F_t, \, t \in [0,1]$, of continuous mappings from $\mathcal{X}$ to $\mathcal{Y}$. Consider the curve $\alpha(t)\stackrel{def}=F_t(x_0),\, t\in [0,1]$. Suppose $F_0(x_0)=y_0$ and $F_1(x_0)=y_1$. Each curve $\beta$ in $\mathcal{Y}$ with initial point $y_0$ and terminating point $y_1$ defines an isomorphism  $\pi_1(\mathcal{Y},y_0) \to \pi_1(\mathcal{Y},y_1) $ by considering compositions of curves $\beta^{-1} \gamma \beta$ for loops $\gamma$ representing elements of $\pi_1(\mathcal{Y},y_0)$. For different curves $\beta$ the isomorphisms differ by conjugation. Identifying fundamental groups with different base point by a chosen isomorphism of the described kind
we obtain a mapping that assigns to each free homotopy class of maps  $\mathcal{X}\to \mathcal{Y}$ a conjugacy class of homomorphisms $\pi_1(\mathcal{X},x_0) \to \pi_1(\mathcal{Y},y_0) $.
The mapping is surjective. It is also injective by the following reason. Let $F_j:\mathcal{X}\to \mathcal{Y}$ be continuous mappings, $F_j(x_0)=y_0$, such that for an element $e\in \pi_1(\mathcal{Y},y_0)$ the equality $(F_1)_*=e^{-1}(F_0)_*\,e$ holds. Let $\alpha$ be a smooth curve in $\mathcal{Y}$ that represents $e$. There exists a free homotopy $F^t,\, t\in [0,1],$ of mappings such that $F^0=F_0$ and $F^t(x_0)=\alpha(t),\, t\in[0,1]$. Then $(F^1)_*=e^{-1}(F^0)_*\,e$. Hence, $F_1$ and $F^1$ define the same homorphism $\pi_1(\mathcal{X},x_0) \to \pi_1(\mathcal{Y},y_0) $, and hence, they are homotopic.
We obtain the following theorem (see also \cite{Ha},\cite{St}, \cite{Sp}.)

\begin{thm}\label{thmEl-1}
The free homotopy classes of continuous mappings from $\mathcal{X}$ to $\mathcal{Y}$
are in one-to-one correspondence to the set of conjugacy classes of homomorphisms between the fundamental groups of $\mathcal{X}$ and $\mathcal{Y}$.
\end{thm}

The crucial step for the proof of Theorem \ref{thm8.0} is the following Theorem \ref{thm8.1}.

\begin{thm}\label{thm8.1} Suppose $X$ is a smooth surface of genus one with a hole.
Let $\mathcal{E}=\{e_1,e_2\}$ be a chosen set of generators of the fundamental group of $X$ with base point $x_0$, and $\mathcal{E}_0= \{e_1,e_2, e_2 e_1^{-1}, e_2 e_1^{-2}, e_1 e_2 e_1^{-1} e_2^{-1}\}$. If $f$ is a separable quasipolynomial of degree $3$ on $X$ for which $\mathcal{M}(\hat{f_e})>\frac{\pi}{2}(\log\frac{3+\sqrt{5}}{2})^{-1}$ for $e\in \mathcal{E}_0$,
then the isotopy class of $f$
corresponds to the conjugacy class of a homomorphism
$$
\Phi : \pi_1 (X , x_0) \to \Gamma \subset {\mathcal B}_3
$$
for a subgroup $\Gamma$ of ${\mathcal B}_3$ \index{$\Phi$}
which is generated either by $\sigma_1 \, \sigma_2$, or by $\sigma_1 \, \sigma_2 \, \sigma_1$,
or by $\sigma_1$ and $\Delta_3^2$.
In particular, in all cases the image $\Phi ([e_1 , e_2])$ of the
commutator of the generators is the identity.
\end{thm}

The Theorems \ref{thm8.1} and \ref{thm8.0} describe the set of mappings from a smooth surface of genus one with a hole to $\mathfrak{P}_3$ that have the Gromov-Oka property.

The rest of this section is devoted to the proof of Theorem \ref{thm8.1}. Theorem \ref{thm8.0} and  results related to bundles will be proved in section \ref{sec:9.3}.

We start with the proof of Lemma \ref{lem0''}, and then provide several lemmas that are needed for the proof of Theorem \ref{thm8.1}.

\medskip

\noindent {\bf Proof of lemma \ref{lem0''}.} Take any conformal structure $\omega':X\to \omega'(X)$ on $X$ with only thick ends. Let $\gamma$ be a
smooth simple closed curve that represents the conjugacy class $\hat e$ of $e$.
Cut $\omega'(X)$ along $\omega'(\gamma)$. A tubular neighbourhood of $\omega'(\gamma)$ is cut into connected components $X_+$ and $X_-$, each being conformally equivalent to an annulus. Take an annulus $A$ of conformal module $m(A)>a$, and glue the $X_{\pm}$ conformally to annuli $A_{\pm}$ in $A$ that
are adjacent to different boundary components of $A$. If $a$ is large then $A_{\pm}$ are disjoint, and
we obtain a Riemann surface and a homeomorphism $\omega$ from $X$ onto this Riemann surface. The mapping $\omega:X\to \omega(X)$ gives the require complex structure. \hfill $\Box$

\medskip

The proof of the following lemma on permutations can be extracted for
instance from \cite{W}. For convenience of the reader we give the
short argument.

\begin{lemm}\label{lemm8.4} Let $n$ be a prime number. Any Abelian subgroup of
the symmetric group ${\mathcal S}_n$ which acts transitively on a
set consisting of $n$ points is generated by an $n$-cycle.
\end{lemm}

\noindent {\bf Proof.} Let ${\mathcal S}_n$ be the group of
permutations of elements of the set $\{ 1,\ldots , n \}$. Suppose
the elements $s_j \in {\mathcal S}_n$, $j=1,\ldots , m$, commute and
the subgroup $\langle s_1 , \ldots , s_m \rangle$ of ${\mathcal
S}_n$ generated by the $s_j$, $j = 1,\ldots , m$, acts transitively
on the set $\{ 1,\ldots , n\}$. Let $A_{s_1} \subset \{ 1,\ldots ,
n\}$ be a minimal $s_1$-invariant subset. Then $s_1 \mid A_{s_1}$ is
a cycle of length $k_1 = \vert A_{s_1} \vert$. (The order $\vert A
\vert$ of a set $A$ is the number of elements of this set.) For any
integer $\ell$ the set $s_2^{\ell} (A_{s_1})$ is minimal
$s_1$-invariant. Hence, two such sets are either disjoint or equal.
Take the minimal union $A_{s_1 s_2}$ of sets of the form $s_2^{\ell}
(A_{s_1})$ for some integers $\ell$ which contains $A_{s_1}$ and is
invariant under $s_2$. Then $s_2$ moves the $s_2^{\ell} (A_{s_1})$
along a cycle of length $k_2$ such that $\vert A_{s_1 s_2} \vert =
k_1 \cdot k_2$. $A_{s_1 s_2}$ is a minimal subset of $\{ 1,\ldots ,
n\}$ which is invariant for both, $s_1$ and $s_2$. Continue in this
way. We obtain a minimal set $A = A_{s_1 \ldots s_m} \subset
\{1,\ldots , n\}$ which is invariant under all $s_j$. By the
transitivity condition $A = \{1,\ldots , n\}$. The order $\vert A
\vert$ equals $k_1 \cdot \ldots \cdot k_m$. Since $n$ is prime,
exactly one of the factors, $k_{j_0}$ equals $n$, the other factors
equal $1$. Then $\vert A_{s_1 \ldots s_{j_0 - 1}} \vert = 1$, and
$s_{j_0}$ is a cycle of length $n$. Since each $s_j$ commutes with
$s_{j_0}$, each $s_j$ is a power of $s_{j_0}$. Indeed, consider a
bijection of $ \{1,\ldots , n\}$ onto the set of $n$-th roots of
unity so that $s_{j_0}$ corresponds to rotation by the angle
$\frac{2\pi}n$. In other words, put $\zeta = e^{\frac{2\pi i}n}$.
The permutation $s_{j_0}$ acts on $\{ 1 , \zeta , \zeta^2 , \ldots ,
\zeta^{n-1}\}$ by multiplication by $\zeta$. Consider an arbitrary
$s_j$. Then for some element $\ell_j$ of $\mathbb{Z}_n=\mathbb{Z}\diagup n\mathbb{Z}$, the equality $s_j (\zeta) =
\zeta^{\ell_j}$ holds, i.e. $s_j (\zeta) = \zeta^{\ell_j-1} \cdot \zeta =
(s_{j_0})^{\ell_j-1} (\zeta)$. Then for any other $n$-th root of
unity $\zeta^m$
\begin{align*}
s_j (\zeta^m) = & s_j ((s_{j_0}) \ ^{m-1}(\zeta))\\
= & (s_{j_0}) \
^{m-1}(s_j(\zeta))
=  \zeta^{m-1} \cdot \zeta^{\ell_j}
=  \zeta^{\ell_j-1} \cdot \zeta^m = (s_{j_0})^{\ell_j-1} (\zeta^m) \, .
\end{align*}
Hence $s_j = (s_{j_0})^{\ell_j-1}$. \hfill $\Box$

\begin{lemm}\label{lemm8.5} Let $e_1 , e_2$ be generators of a free group $F_2$.
Let $E \subset F_2$ be the finite subset $E = \{ e_1 , e_2 , e_2 \,
e_1^{-1} , e_2 \, e_1^{-2} \}$. Put $E_- \overset{\rm def}{=}
\{e^{-1} : e \in E \}$. Suppose $\Psi : F_2 \to {\mathcal S}_3$ is a
homomorphism from $F_2$ into the symmetric group ${\mathcal S}_3$
whose image is an Abelian subgroup of ${\mathcal S}_3$ which acts
transitively on the set of three elements. Then there are generators
${\sf e}_1 , {\sf e}_2 \in E \cup E_-$ of $F_2$ such that $\Psi ({\sf e}_1)$
is a $3$-cycle, $\Psi ({\sf e}_2) = {\rm id}$ and the commutator
$[{\sf e}_1 , {\sf e}_2]$ is conjugate to $[e_1 , e_2]$.
\end{lemm}

\noindent {\bf Proof.} By (the proof of) Lemma \ref{lemm8.4} the
image of one of the original generators is a $3$-cycle. If $\,\Psi
(e_1)$ is a $3$-cycle, then $\Psi (e_2 \, e_1^{-q})$ is the identity
for $q$ being either $0$ or $1$ or $2$. (Recall that by Lemma \ref{lemm8.4} $\Psi (e_2)$ is
a power of $\Psi (e_1)$ and $\Psi (e_1)^3 = {\rm id}$.) Also,
$$
[e_1 , e_2 \, e_1^{-q}] \, =\,  e_1 \, e_2 \, e_1^{-q} \, e_1^{-1}
\, e_1^q \, e_2^{-1}\, = \, [e_1 , e_2] \, .
$$
If $\Psi (e_1)$ is the identity then $\Psi (e_2)$ is a $3$-cycle.
For the pair $\,{\sf e}_1 = e_2\,$, ${\sf e}_2 = e_1^{-1}\,$ of
generators the commutator $\,[{\sf e}_1 , {\sf e}_2]\, =\, e_2 \,
e_1^{-1} \, e_2^{-1} \, e_1\,$ is conjugate to the commutator
$\,[e_1 , e_2]\,.\,$ \hfill $\Box$

\medskip

We will now consider a torus with a hole $X$. Recall that its fundamental group $\pi_1
(X,x_0)$ with base point $x_0$ is isomorphic to the free group $F_2$
with two generators.

For the proof of Theorem \ref{thm8.1} we need two further lemmas which we state and prove before proving the theorem.

\begin{lemm}\label{lemm8.7} Let ${\sf b}_1 = (\sigma_1 \, \sigma_2)^{\pm 1}$ and
${\sf b}_2 = w^{-1} \, \sigma_1^{2k} \, w$ for an integer $k$ and a
braid $w \in {\mathcal B}_3$. If the commutator $[{\sf b}_1 , {\sf
b}_2]$ is conjugate to $\sigma_1^{2k^*} \Delta_3^{2\ell^*}$ for
integers $k^*$ and $\ell^*$ then ${\sf b}_2$ is the identity and,
hence, also the commutator is the identity.
\end{lemm}

\noindent {\bf Proof of Lemma \ref{lemm8.7}.} Let $b \in {\mathcal B}_n$ be a
pure braid (considered as an isotopy class of geometric braids with
strands labeled). The linking number $\ell_{ij}$ of the
$i$-th and $j$-th strand is defined as follows. Discard all strands except the
$i$-th and the $j$-th strand. We obtain a pure braid $\sigma^{2m}$.
We call the integer number $m$ the linking number of the two
strands and denote it by $\ell_{ij}$.

Note that the linking numbers ${\ell}_{ij}^*$ of the braid
$\sigma_1^{2k^*} \Delta_3^{2\ell^*}$ are equal to ${\ell}_{12}^* = k^* + \ell^*$, ${\ell}_{23}^* =  \ell^*$ ,
${\ell}_{13}^* = \ell^*$. Since the braid is conjugate to a
commutator, the sum of exponents of generators in a word representing it must be zero. This means that $2k^* + 6\ell^* = 0$. Hence,
the (unordered) collection of linking numbers of $\sigma_1^{2k^*}
\Delta_3^{2\ell^*}$ is $\{ -2\ell^* , \ell^* , \ell^*\}$. Since
conjugation only permutes linking numbers between pairs of strands
of a pure braid, the unordered collection of linking numbers of
pairs of strands of $[{\sf b}_1 , {\sf b}_2]$ equals $\{ -2\ell^* ,
\ell^* , \ell^* \}$ for the integer $\ell^*$.

For the pure braid
$\sigma_1^{2k}$ the linking number between the first and the second
strand equals $k$, the linking numbers of the remaining pairs of
strands equal zero.

For a braid $w \in \mathcal{B}_3$ let $S_w$ be the permutation $S_w
= \tau_3 (w)$. We identify $\mathcal{B}_3$ with the fundamental group of $C_3(\mathbb{C})\diagup \mathcal{S}_3$ with base point $E_3$. $E_3$ can be considered as an unordered subset of $\mathbb{C}$ consisting of three points.
$S_w$ permutes this three points. We label the points of $E_3$ by $(x_1,x_2,x_3)$. The permutation $S_w$ defines
a permutation $s_w$ acting on the set
$(1,2,3)$, so that $S_w((x_1,x_2,x_3))=(x_{s_w(1)} , x_{s_w(2)} ,
x_{s_w(3)})$. Order the linking numbers of pairs of strands of a pure braid
$b \in \mathcal{B}_3$ as $(\ell_{23} ,
\ell_{13} , \ell_{12})$.
Since
$$
(s_w(1),s_w(2),s_w(3))
\overset{s_w^{-1}}{-\!\!\!-\!\!\!\longrightarrow} (1,2,3)
\overset{s_w}{-\!\!\!-\!\!\!\longrightarrow} (s_w(1),s_w(2),s_w(3)),
$$
the ordered tuple of linking numbers of pairs of strands of $w^{-1}
\, b \, w$ equals
$$
S_w((\ell_{23} , \ell_{13} , \ell_{12}))\, = \,
  (\ell_{s_w(2) s_w(3)} , \ell_{s_w(1) s_w(3)} ,
  \ell_{s_w(1)s_w(2)}).
$$

Hence, the ordered tuple of linking numbers between pairs of strands
of ${\sf b}_2 = w^{-1} \, \sigma_1^{2k} \, w$ equals $S_w (0,0,k)$.
Similarly, the ordered tuple of linking numbers between pairs of
strands of ${\sf b}_2^{-1}$ is equal to $S_w (0,0,-k)$.

Consider the commutator ${\sf b}_1 \, {\sf b}_2 \, {\sf b}_1^{-1} \,
{\sf b}_2^{-1}$.

\begin{figure}[h]
\begin{center}
\includegraphics[width=55mm]{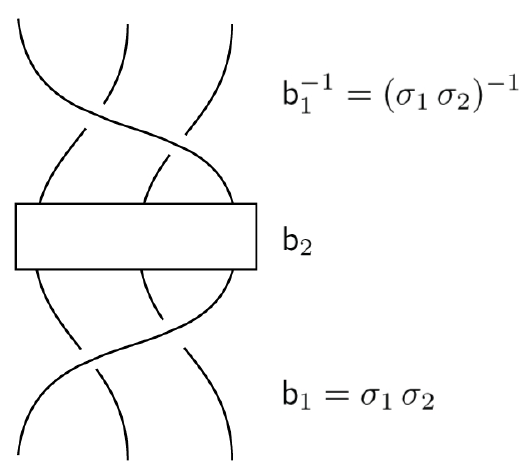}
\end{center}
\caption{The braid $b_1 b_2 b_1^{-1}= (\sigma_1 \, \sigma_2) b_2 (\sigma_1 \, \sigma_2)^{-1}$}\label{fig8.1}
\end{figure}

The ordered tuple of linking numbers of pairs of strands of the pure
braid ${\sf b}_1 \, {\sf b}_2 \, {\sf b}_1^{-1}$ equals $(S')^{-1}
\circ S_w (0,0,k)$ for the permutation $S' = \tau_3 ({\sf b}_1)$
(see Figure \ref{fig8.1}).

\noindent Hence the ordered tuple of linking numbers of pairs of
strands of the commutator ${\sf b}_1 \, {\sf b}_2 \, {\sf b}_1^{-1}
\circ {\sf b}_2^{-1}$ equals $(S')^{-1} \circ S_w
(0,0,k) + S_w (0,0,-k)$. Since $S'$ is a 3-cycle, the mapping  $S' (x_1 , x_2 ,
x_3) = (x_{s'(1)} , x_{s'(2)} , x_{s'(3)})$ does not fix (setwise) any pair
of points among the $x_1, x_2, x_3$. Hence the unordered
$3$-tuple of linking numbers of ${\sf b}_1\, {\sf b}_2 \, {\sf
b}_1^{-1} \, {\sf b}_2^{-1}$ is $\{ k,-k,0 \}$. It can coincide with
an unordered $3$-tuple of the form $\{-2 \ell^* , \ell^* , \ell^*\}$
only if $k=\ell^*=0$. Hence $\textsf{b}_2$ is the identity and the
commutator is the identity. \hfill $\Box$

\begin{lemm}\label{lemm8.8} The centralizer of $\sigma_1^k$ in $\mathcal{B}_3$, $k \ne 0$
an integer number, equals $\{\sigma_1^{k'} \, \Delta_3^{2\ell'} :
k',\ell' \in {\mathbb Z}\}$.
\end{lemm}

\noindent {\bf Proof.} We use the standard homomorphism $\vartheta :
{\mathcal B}_3 \to {\rm SL} (2,{\mathbb Z})$ for which
\index{$\vartheta$}
$$
\vartheta (\sigma_1) = {\sf A} \overset{\rm def}{=} \begin{pmatrix}
1&1 \\ 0&1 \end{pmatrix} \, , \quad \vartheta (\sigma_2) = {\sf B}
\overset{\rm def}{=} \begin{pmatrix}{\;\;\;\, 1}&0 \\\ -1&{1} \end{pmatrix} \, ,
$$
(see e.g. \cite{KaTu}). Here ${\rm SL} (2,{\mathbb Z})$ denotes the group of $2 \times 2$
matrices with determinant equal to $1$ and entries being integer numbers.  The kernel of $\vartheta$ is
generated by $(\sigma_1 \, \sigma_2 \, \sigma_1)^4$. For $b \in
{\mathcal B}_3$ we denote by $\langle b \rangle$ the subgroup of
${\mathcal B}_3$ generated by $b$. With this notation
$$
{\mathcal B}_3 \diagup \langle (\sigma_1 \, \sigma_2 \, \sigma_1)^4
\rangle
$$
is isomorphic to ${\rm SL} (2,{\mathbb Z})$. Notice that $\vartheta
((\sigma_1 \, \sigma_2 \, \sigma_1)^2) = - \begin{pmatrix} 1&0 \\
0&1 \end{pmatrix}$.

Let $b \in {\mathcal B}_3$ be in the centralizer of $\sigma_1^k$.
Then ${\sf V} \overset{\rm def}{=} \vartheta (b)$ is in the
centralizer of ${\sf A}^k$, in other words,
$$
\begin{pmatrix} 1&k \\ 0&1 \end{pmatrix} \begin{pmatrix} v_{11}&v_{12} \\ v_{21}&v_{22} \end{pmatrix} = \begin{pmatrix} v_{11}&v_{12} \\ v_{21}&v_{22} \end{pmatrix} \begin{pmatrix} 1&k \\ 0&1 \end{pmatrix} \quad \mbox{with} \quad {\sf V} =  \begin{pmatrix} v_{11}&v_{12} \\ v_{21}&v_{22} \end{pmatrix} \, ,
$$
i.e.
$$
\begin{pmatrix} v_{11} + k \, v_{21}&\;\;v_{12} + k \, v_{22} \\ v_{21}&\; \;v_{22} \end{pmatrix} =  \begin{pmatrix} v_{11}&\;\;k \, v_{11} + v_{12} \\ v_{21}&\;\;k \, v_{21} + v_{22} \end{pmatrix} \, .
$$
Since $k \ne 0$ we have $v_{21} = 0$ and $v_{11} = v_{22}$. Since
$\det {\sf V} = 1$ we obtain for an integer $m$
$$
{\sf V} = \begin{pmatrix} 1&m \\ 0&1 \end{pmatrix} \quad \mbox{or}
\quad {\sf V} = \begin{pmatrix} -1&-m \\ 0&-1 \end{pmatrix} = -
\begin{pmatrix} 1&m \\ 0&1 \end{pmatrix} \, .
$$
Hence, either $\vartheta (b \, \sigma_1^{-m}) = {\rm id}$ or
$\vartheta (b \, \sigma_1^{-m} \, \Delta_3^2) = {\rm id}$. Thus $b
\, \sigma_1^{-m} = \Delta_3^{2\ell}$ for some integer $\ell$. The
lemma is proved. \hfill $\Box$

\medskip

\noindent {\bf Proof of Theorem \ref{thm8.1}.} Let $\Phi : \pi_1 (X,x_0) \to {\mathcal B}_3$
be a homomorphism whose conjugacy class corresponds to the isotopy
class of the quasipolynomial $f$ considered as mapping $f:X\to \mathfrak{P}_3$. Then $\hat{f}_e=\widehat{\Phi(e)}$ for each $e\in \pi_1 (X,x_0)$. Denote by $\Psi = \tau_3 \circ \Phi : \pi_1 (X,x_0)
\to {\mathcal S}_3$ the related homomorphism into the symmetric
group.
Since $\mathcal{M}(\hat{f}_{[e_1,e_2]})>\frac{\pi}{2}(\log\frac{3+\sqrt{5}}{2})^{-1}$,
by Lemma \ref{lemm8.3} the conformal module of the commutator class $\widehat{[\Phi(e_1),\Phi(e_2)]}$ equals infinity,
and $\,\Psi ([e_1 , e_2]) \,= \,{\rm id}\,.\,$
\index{$\Psi$}

Consider first the case when the subgroup $\Psi
(\pi_1 (X,x_0))$ of ${\mathcal S}_3$ acts transitively on the set of
three points. Apply Lemma \ref{lemm8.5} to the free group $\pi_1
(X,x_0)$ with generators $e_1$ and $e_2$ and to the homomorphism
$\Psi$. We obtain new generators ${\sf e}_1 , {\sf e}_2 \in E \cup
E_-$ (with $E$ and $E_-$ being the sets of the lemma) such that
$\Psi ({\sf e}_1)$ is a $3$-cycle and $\Psi ({\sf e}_2) = {\rm id}$.
Put $s_1 = \Psi ({\sf e}_1)$, $s_2 = \Psi ({\sf e}_2)$, $b_1 = \Phi
({\sf e}_1)$, $b_2 = \Phi ({\sf e}_2)$.
Then $b_2$ and $[b_1,b_2]$
are pure braids (the latter holds since $[{\sf e}_1 , {\sf e}_2]$ is
conjugate to $[e_1 , e_2]$).
Since $\mathcal{M}(\hat{f}_{{\sf e}_1})>\frac{\pi}{2}(\log\frac{3+\sqrt{5}}{2})^{-1}$,
Lemma \ref{lemm8.3} implies that ${\mathcal M} (\widehat{b_1}) = \infty$.
Since $\tau_3 (b_1)$ is a $3$-cycle, by the same lemma
the braid $b_1$ must be conjugate to an integer power of $\sigma_1
\, \sigma_2$.

After conjugating $\Phi$ we may assume that $b_1 = (\sigma_1 \,
\sigma_2)^{3\ell \pm 1} = (\sigma_1 \, \sigma_2)^{\pm 1} \cdot
\Delta_3^{2\ell}$ for an integer $\ell$. Since $\mathcal{M}(\hat{f}_{{\sf e}_2})>\frac{\pi}{2}(\log\frac{3+\sqrt{5}}{2})^{-1}$,
Lemma \ref{lemm8.3} implies as above the equality
${\mathcal M} (\widehat{b_2}) = \infty$ for the pure braid $b_2$. Hence, $b_2
= w^{-1} \, \sigma_1^{2k} \, \Delta_3^{2\ell'} \, w$ for integers
$k$ and $\ell'$ and a conjugating braid $w \in {\mathcal B}_3$.
Since ${\mathcal M} (\widehat{[b_1 , b_2]}) = \infty$ the commutator $[b_1 , b_2]$ is
conjugate to $\sigma_1^{2k^*} \Delta_3^{2\ell^*}$ for integers $k^*$
and $\ell^*$.

Take into account that $\Delta_3^{2\ell}$ commutes with each $3$-braid
and apply Lemma \ref{lemm8.7} to $B_1\stackrel{def}= (\sigma_1 \, \sigma_2)^{\pm 1}$ and $B_2\stackrel{def}= w^{-1} \sigma_1^{2k^*} w$.
This gives the statement of Theorem \ref{thm8.1} for the case
when the subgroup $\Psi(\pi_1(X,x_0))$ acts transitively on the set of three points.

\smallskip

Consider now the case when the subgroup $\Psi
(\pi_1 (X,x_0))$ of ${\mathcal S}_3$ does not act transitively on the set of
three points.
Then $\Psi (\pi_1 (X,x_0))$ is generated either by a
transposition (we may assume the transposition to be $(12)$ by
choosing the mapping  $\Phi$ in its conjugacy class) or $\Psi (\pi_1 (X,x_0))$ is equal to
the identity in ${\mathcal S}_3$.

There are generators ${\sf e}_1 , {\sf e}_2 \in E \cup E_-$ such
that $[{\sf e}_1 , {\sf e}_2]$ is conjugate to $[e_1 , e_2]$ and
$\Psi ({\sf e}_2) = {\rm id}$. Indeed,
we need to show this for the case, when $\Psi(\pi_1(X;x_0))$
is generated by the transposition $(12)$ and neither $\Psi(e_1)$ nor
$\Psi(e_2)$ is the identity. Then $\Psi(e_1)=\Psi(e_2)=(12),\,$ and
$\Psi(e_2 e_1^{-1})=\rm {id}$. Hence as in the proof of Lemma
\ref{lemm8.5} the generators ${\sf e}_1 $ and $ {\sf e}_2$ can be
chosen as required.

Use again the notation $b_1 = \Phi ({\sf e}_1)$, $b_2 = \Phi ({\sf e}_2)$.
Since $\mathcal{M}(\hat{f}_{{\sf e}_1})$ and $\mathcal{M}(\hat{f}_{{\sf e}_2})$ are both bigger than $\frac{\pi}{2}(\log\frac{3+\sqrt{5}}{2})^{-1}$
we may assume by Lemma
\ref{lemm8.3}, that $b_1$ is conjugate to either
$\sigma_1^{ k_1} \, \Delta_3^{2\ell_1}$ for integers $\ell_1$ and $k_1$, or to $(\sigma_1 \sigma_2
\sigma_1)^{k_1}= \Delta_3^{k_1}\,$ for an integer $k_1$,
and $b_2$ is conjugate to $\sigma_1^{2 k_2} \,
\Delta_3^{2\ell_2}$ for integers $k_2$ and $\ell_2$.
Conjugating $\Phi$, we may assume that $b_2 =B_2\,\Delta_3^{2 \ell_2}$ with $B_2= \sigma_1^{2 k_2}$ for integers $\ell_2$ and $k_2$,
$b_1 = w^{-1} B_1 \, \Delta_3^{2 \ell_1} w$ for a braid $w\in\mathcal{B}_3$
and an integer $\ell_1$ with either $B_1 =\sigma_1 ^{k_1}$ for an integer $k_1$, or
$B_1 = \Delta_3$.
Since $\mathcal{M}(\hat{f}_{[{\sf e}_1,{\sf e}_1]}) > \frac{\pi}{2}(\log\frac{3+\sqrt{5}}{2})^{-1} $,
Lemma \ref{lemm8.3} implies as above that the commutator
$[B_1 , B_2]= [b_1 , b_2]$ is a pure braid with infinite conformal module,
and, hence,  $[B_1 , B_2]= [b_1 , b_2]$ is conjugate to
$\sigma_1^{2k^*} \, \Delta_3^{2\ell^*}$. Since the commutator has
exponent sum zero, the equality $2k^* + 6\ell^* = 0$ holds, and the
unordered tuple of linking numbers of the commutator equals $\{ \ell^* , \ell^*
,-2\ell^* \}$.

This implies that the commutator $[B_1 , B_2]$
is equal to the identity.
Indeed, the ordered tuple of linking numbers of pairs of strands of $B_2$ equals $(0,0,k_2)$, the ordered tuple of linking numbers of pairs of strands of $B_2^{-1}$ equals $(0,0,-k_2)$.
The ordered tuple of
linking numbers of pairs of strands of $B_1 \, B_2 \, B_1^{-1}$
equals $ (S')^{-1}(0,0,k_2)$, where $S' =
\tau_3 (B_1)$.
The unordered tuple of linking numbers of pairs of strands of
$[B_1 , B_2]$ equals either $\{ k_2 , -k_2 , 0\}$ or $\{0,0,0\}$
(in dependence on $S'$). Either of these unordered tuples can be
equal to $\{\ell^* , \ell^* , -2\ell^*\}$ for some integer $\ell^*$
only if $\ell^* = 0$. We obtained that $[b_1,b_2] =[B_1,B_2]= \rm{id}$.
If $k_2\neq0$, then by Lemma \ref{lemm8.8} the group $\Phi(\pi_1(X,x_0))$ is generated
by (a power of) $\sigma_1$ and $\Delta_3^2$. If $k_2=0$ then $\Phi(\pi_1(X,x_0))$ is generated by $b_1$, and hence, a conjugate group is generated either by $\sigma_1$ and $\Delta_3^{2}$, or by $\Delta_3$.
The theorem is proved. \hfill $\Box$

\medskip

We do not know the answer to the following problem.
\begin{probl}\label{probl8.1}  Can two braids in ${\mathcal B}_3$ of zero entropy
have a non-trivial commutator of zero entropy?
\end{probl}
If the answer is negative, then Theorem \ref{thm8.1} holds with $\mathcal{E}_0=\mathcal{E}$.
The proof of Theorem \ref{thm8.1}
gives the following corollary which is weaker than a negative answer
to Problem \ref{probl8.1}.
\begin{cor}\label{cor8.1}
Let ${b}_1$ and ${b}_2$ be braids in ${\mathcal{B}}_3$ with
$h(b_1)=h(b_2)=h([b_1,b_2])=0$. If one of the braids is pure or
$$h(b_2\, b_1^{-1})=h(b_2\,b_1^{-2})=0$$
then $[b_1,b_2]= \rm{id}$.
\end{cor}

On the other hand D.Calegari and A.Walker suggested the following example of two elements
of the braid group $\mathcal{B}_3$ with non-trivial commutator of
zero entropy.
\begin{ex}\label{ex8.1}
The commutator of the non-commuting braids $b_1= (\sigma_2)^{-1} \, \sigma_1\;$ and $ b_2 =
\sigma_2 \, (\sigma_1)^{-1}\,$  has entropy zero.
\end{ex}
It is enough to show that the commutator of the two braids equals $[b_1,b_2] =
(\sigma_2)^{-6} \, \Delta_3^2\,.\,$  Since
\begin{align*}
b_1^{-1} \, b_2^{-1} \,=\,&(\sigma_1)^{-1} \, \sigma_2 \, \sigma_1 \;
(\sigma_2)^{-1}\;\;\;\;\;  =\,(\sigma_1)^{-1} \, \sigma_2 \, \sigma_1
\,\sigma_2 \; (\sigma_2)^{-2} \, \\
=\,&(\sigma_1)^{-1} \, \sigma_1 \,
\sigma_2 \,\sigma_1 \; (\sigma_2)^{-2}
\, = \, \sigma_2 \, \sigma_1 \,(\sigma_2)^{-2}\, ,\\
b_1 \, b_2 \; \, =\,&(\sigma_2)^{-1} \, \sigma_1 \, \sigma_2 \;
(\sigma_1)^{-1}\;\; \;\;\; =\, (\sigma_2)^{-2} \, \sigma_2 \, \sigma_1 \,
\sigma_2 \;(\sigma_1)^{-1}\,\\
 =\,&(\sigma_2)^{-2} \, \sigma_1 \,
\sigma_2 \, \sigma_1 (\sigma_1)^{-1}
\,\; = \, (\sigma_2)^{-2} \, \sigma_1 \,\sigma_2\, ,
\end{align*}
we obtain
\begin{align*}
[b_1,b_2] = \, b_1 b_2 b_1^{-1}b_2^{-1} \, = \,&(\sigma_2)^{-3} \,
\sigma_2 \, \sigma_1 \, \sigma_2 \, \cdot \sigma_2 \, \sigma_1 \,
\sigma_2 \; (\sigma_2)^{-3} \,\\
=  \,&(\sigma_2)^{-3} \; \Delta_3^2 \;
(\sigma_2)^{-3}
= \, (\sigma_2)^{-6} \;\Delta_3^2\,.
\end{align*}

\begin{figure}[h]
\begin{center}
\includegraphics[width=7cm]{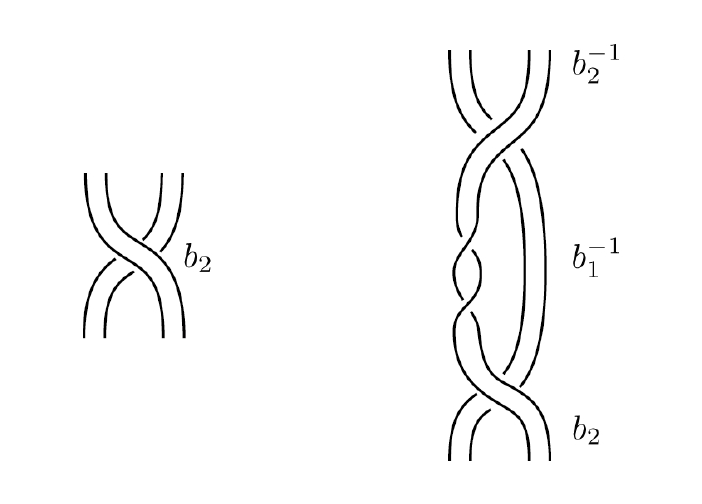}
\caption{Two $4$-braids of vanishing entropy with non-trivial commutator}\label{Fig8.2}
\end{center}
\end{figure}

Problem \ref{probl8.1} has a positive answer for braid groups on more than $3$ strands.
Consider the braids $ b_1 =\sigma_1^{-2}$ and $b_2 = (\sigma_2)^{-1} \,
(\sigma_1)^{-1} \, (\sigma_3)^{-1} \, (\sigma_2)^{-1}$ in $\mathcal{B}_4$ of zero entropy.
Their commutator equals $[b_1
, b_2] = \sigma_1^{-2} \cdot \sigma_3^{2} \ne {\rm id}$ and is a non-trivial braid of zero entropy (see Figure \ref{Fig8.2}).
Notice the following fact. Let $X$ be a surface of genus one with a disc removed.
Consider the separable quasipolynomial of degree four on $X$ whose
isotopy class corresponds to the conjugacy class of the homomorphism
$\Phi : \pi_1 (X,x_0) \to {\mathcal B}_4$ with $\Phi (e_1) = b_1$ and $\Phi (e_2) = b_2$.
By a similar argument as in the proof of Theorem \ref{thm8.0} (see Section \ref{sec:9.1})
the quasipolynomial is isotopic to a separable algebroid function for each
conformal structure of second kind on $X$.

\begin{probl}\label{probl8.2.} What is the analogue of Corollary
\ref{cor8.1} for braid groups on more than three strands? What is the analog of Theorem \ref{thm8.0} for quasipolynomials of degree $4$?
\end{probl}

\section{$(\sf{g,m})$-fiber bundles}
\label{sec:9.1}

We will consider bundles whose fibers are
closed surfaces of genus $\textsf{g}\geq 0$ with $\textsf{m}\geq 0$ distinguished points.
For brevity we call a smooth closed surface of genus $\sf{g}$ with $\sf{m}$ distinguished
points a surface of type $(\sf{g},\sf{m})$. In case the surface is equipped with a complex
structure we call it a Riemann surface of type $(\sf{g},\sf{m})$.

\begin{defn}\label{def9.1}{\rm (Smooth oriented $(\textsf{g},\textsf{m})$ fiber bundles.)}
Let $X$ be a smooth oriented manifold of dimension $k$, let
${\mathcal X}$ be a smooth (oriented)
manifold of dimension $k+2$ and ${\mathcal P} : {\mathcal X} \to X$
an orientation preserving smooth proper submersion such that for each point $x \in X$ the
fiber ${\mathcal P}^{-1} (x)$ is a smooth closed oriented surface of genus $\sf{g}$.
Let $\mathbold{E}$ be a smooth submanifold of $\mathcal{X}$ that intersects each fiber $\mathcal{P}^{-1}(x)$ along a set $E_x$ of $\sf{m}$ distinguished points.
Then the tuple ${\mathfrak F}_{\textsf{g},\textsf{m}} = ({\mathcal
X} , {\mathcal P} , \mathbold{E}, X)$ is called
a smooth (oriented) fiber bundle over $X$ with fibers being smooth closed oriented surfaces of
genus $\textsf{g}$ with $\textsf{m}$ distinguished points
(for short, a smooth oriented $(\textsf{g},\textsf{m})$-bundle).
\end{defn}
If ${\sf m}=0$ the set $\mathbold E$ is the empty set and we will often denote the bundle
by $({\mathcal X} , {\mathcal P} , X)$.
If ${\sf{m}}>0$ the mapping $x\to E_x$ locally defines $\sf m$ smooth sections of the $({\sf g},0)$-bundle $({\mathcal X} , {\mathcal P} , X)$ that is obtained by forgetting all distinguished points. A section of the bundle $({\mathcal X} , {\mathcal P} , X)$ is a continuous mapping $s:X\to {\mathcal X}$ with $\mathcal{P}\circ s$ being the identity on $X$.
$(\textsf{g},0)$-bundles will also be called genus $\textsf{g}$ fiber bundles.
For $\textsf{g}=1$ and $\textsf{m}=0$ the bundle is also called an elliptic fiber
bundle.  \index{fiber bundle !  elliptic}
We will consider mostly the case when $2{\sf g}-2 + {\sf m}>0$.

\index{$\mathfrak{F}_{{\sf g},{\sf m}}$} \index{fiber-bundle ! ${\sf g},{\sf m}$-fiber bundle}
Let $S$ be a smooth reference surface of genus $\sf g$ and $E\subset S$ a set of $\sf{m}$ distinguished points.
For an open subset $U$ of $X$ we consider the trivial bundle $\big(U\times S, {\rm {pr}}_1, U\times E,U\big)$
with set $\{x\}\times E$ of distinguished points in the fiber $\{x\}\times S$ over $x$. Here ${\rm pr}_1: U\times S\to U $ is the projection onto the first factor.
By Ehresmann's Fibration Theorem each smooth $(\sf{g},\sf{m})$-bundle ${\mathfrak
F}_{\textsf{g},\textsf{m}} = ({\mathcal
X} , {\mathcal P} , \mathbold{E}, X)$ with set of distinguished points $E_x \stackrel{def}=\mathbold{E}\cap{\mathcal P}^{-1}(x) $ in the fiber over $x$ is
locally smoothly trivial, i.e. each point in $X$ has a neighbourhood $U\subset X$
such that for a
surjective diffeomorphism $\varphi_U:\mathcal{P}^{-1}(U) \to U \times S$
the diagram

\[
\xymatrix{
{\mathcal{P}^{-1}(U)}      \ar[r]^{\varphi_U     } \ar[d]_{\mathcal{P}         } & {U \times
S}   \ar[dl]^{{\rm  pr}_1} \\
  U
}
\]
\noindent is commutative and $\varphi_U$ maps $\mathbold{E}\cap \mathcal{P}^{-1}(U)$ onto $U\times E$, equivalently,
for all $x\in U$ the diffeomorphism
$\varphi_U$ maps
the set of distinguished points $E_x=\mathbold{E}\cap \mathcal{P}^{-1}(x)$ in the fiber $\mathcal{P}^{-1}(x)$ to the set of
distinguished points $\{x\}\times E$ in the fiber $\{x\}\times S$.

The idea of the proof of Ehresmann's Theorem is the following. Choose smooth coordinates on $U$
by a mapping from a rectangular box in $\mathbb{R}^n$ to $U$.
Consider
smooth vector fields $v_j$ on $U$, which form a basis of the tangent
space of $U$ at each point of $U$. Take smooth vector fields $V_j$ on ${\mathcal
P}^{-1} (U)$ that are tangent to $\mathbold{E}$ at points of this set and are mapped to
$v_j$ by the differential of
${\mathcal P}$. Such vector fields can easily be obtained locally.
To obtain the globally defined vector fields $V_j$ on ${\mathcal P}^{-1} (U)$
one uses partitions of unity. The required diffeomorphism $\varphi_U$ is
obtained by composing the flows of the vector fields $V_j$
(in any fixed order).

In this way a trivialization of the bundle can be obtained over any simply connected smooth manifold.

\smallskip

In the case when the base manifold is a Riemann surface,  a
holomorphic $(\textsf{g}$,$\textsf{m})$ fiber bundle over $X$
is defined as follows.
\begin{defn}\label{def2}
Let $X$ be a Riemann surface, let ${\mathcal X}$ be a complex surface, and ${\mathcal P}$
a holomorphic proper submersion from ${\mathcal X}$ onto $X$, such that
each fiber $\mathcal{P}^{-1}(x)$ is a closed Riemann surface of genus $\sf g$.
Suppose $\mathbold{E}$ is a complex one-dimensional submanifold of $\mathcal{X}$ that intersects each fiber $\mathcal{P}^{-1}(x)$ along a
set $E_x$ of $\textsf {m}$ distinguished points. Then the
tuple ${\mathfrak F}_{\textsf{g},\textsf{m}} = ({\mathcal X} , {\mathcal P}, \mathbold{E}, X)$ is called a
holomorphic $(\textsf{g}$,$\textsf{m})$ fiber bundle over $X$.
\end{defn}
Notice that the mapping $x\to E_x$
locally defines $\sf{m}$ holomorphic sections of the $({\sf g},0)$-bundle $({\mathcal X} , {\mathcal P}, X)$ obtained by forgetting all distinguished points.

We will call two smooth oriented (holomorphic, respectively) $(\textsf{g},\,\textsf{m})$ fiber
bundles, ${\mathfrak F}^0 = ({\mathcal X}^0 ,
{\mathcal P}^0 , \mathbold{E}^0, X)$ and ${\mathfrak F}^1 = ({\mathcal X}^1 , {\mathcal
P}^1 ,\mathbold{E}^1, X)$, smoothly isomorphic (holomorphically isomorphic, respectively) if there are smooth (holomorphic,
respectively) homeomorphisms $\Phi: \mathcal{X}^0\to \mathcal{X}^1$ and $\phi:X^0\to X^1$
such that for each $x \in X^0$ the mapping $\Phi$ takes the fiber $(\mathcal{P}^0 )^{-1}(x)$
onto the fiber $(\mathcal{P}^1 )^{-1}(\phi(x))$ and the set of distinguished points in
$(\mathcal{P}^0 )^{-1}(x)$ to the set of distinguished points in $(\mathcal{P}^1 )^{-1}(\phi(x))$.
Holomorphic bundles that are holomorphically isomorphic will be considered the same holomorphic bundles.

For a smooth $({\sf g},{\sf{m}})$-bundle ${\mathfrak F} = ({\mathcal X} ,
{\mathcal P} ,\mathbold{E}, X)$ over $X$ and a homeomorphism $\omega:X\to \omega(X)$ onto a surface $\omega(X)$ we denote by ${\mathfrak F}_{\omega}$ the bundle
\begin{equation}\label{eq23}
{\mathfrak F}_{\omega} = ({\mathcal X} ,
\omega\circ{\mathcal P} ,\mathbold{E},\omega(X) ).
\end{equation}
The bundles ${\mathfrak F}$ and ${\mathfrak F}_{\omega}$ are isomorphic. Indeed, we may consider the homeomorphisms
 $\omega:X\to \omega(X)$ and ${\rm Id}:{\mathcal X}\to {\mathcal X}$.

We are interested in smooth deformations of a smooth bundle to a holomorphic one.
Two smooth (oriented) $(\textsf{g},\textsf{m})$ fiber bundles
over the same oriented
smooth base manifold $X$,  ${\mathfrak F}^0 = ({\mathcal X}^0 ,
{\mathcal P}^0 ,\mathbold{E}^0, X)$, and ${\mathfrak F}^1 = ({\mathcal X}^1 , {\mathcal
P}^1 ,\mathbold{E}^1, X)$, are called (free) isotopic if for an open interval $I$ containing $[0,1]$ there
is a smooth $(\textsf{g},\textsf{m})$ fiber bundle $({\mathcal Y} , {\mathcal P} , \mathbold{E},  X \times
I)$ over the base $X \times I$ (called an isotopy) with the following property.  For each $t \in [0,1]$ we put
${\mathcal Y}^t = {\mathcal P}^{-1} (X \times \{t\})$ and $\mathbold{E}^t= \mathbold{E}\cap {\mathcal P}^{-1} (X \times \{t\})$. The bundle ${\mathfrak F}^0\,$ is
equal to
$\,\left( {\mathcal Y}^0 , {\mathcal P} \mid {\mathcal Y}^0 , \mathbold{E}^0,   X
\times \{0\} \right)\,$, and the bundle $\,{\mathfrak F}^1\,$ is
equal to $\,\bigl( {\mathcal Y}^1 , {\mathcal P} \mid {\mathcal
Y}^1 , \mathbold{E}^1, X \times \{1\} \bigl)\,$.

\index{isotopy ! of fiber bundles}

Notice that for
each $t \in I$ the tuple $\left( {\mathcal Y}^t , {\mathcal P} \mid {\mathcal Y}^t , \mathbold{E}^t, X
\times \{t\} \right)$ is automatically a smooth  $({\sf{g}},{\sf{m}})$-fiber bundle.

Two smooth $(\sf{g,m})$-bundles over the same smooth finite open surface $X$ are isotopic if and only if they are isomorphic
(see Section \ref{sec:9.2}).

The following problem on isotopies of smooth
objects to the respective holomorphic objects concerns another version of the restricted validity of Gromov's Oka principle.

\begin{probl}\label{prob9.1}
Let $X$ be a finite Riemann surface. Can a given smooth $(\textsf{g,m)}$-fiber
bundle over $X$ be smoothly deformed to a
holomorphic fiber bundle over $X$?

More precisely,, does there exist an isotopy $({\mathcal Y} , {\mathcal P} , \mathbold{E},  X \times
I)$ over the base $X \times I$ with $[0,1]\subset I$, for which the fiber over $\{0\}\times X$ is equal to the given bundle, and the fiber over $\{1\}\times X$ can be equipped with the structure of a holomorphic bundle?
\end{probl}
The answer is in general negative. There is a similar obstruction for the existence of such
isotopies for $(\textsf{g,m)}$-bundles as for the existence of isotopies of smooth separable
quasipolynomials to holomorphic ones.
We will call this obstruction the conformal module of isotopy classes of
$(\textsf{g,m)}$-fiber bundles over the circle. It is defined as follows.
\index{conformal module ! of isotopy
classes of fiber bundles} Consider a smooth oriented $({\sf g},{\sf m})$-fiber
bundle ${\mathfrak F} = {\mathfrak F}_{({\sf g},{\sf m})} = ({\mathcal X} , {\mathcal
P} , \mathbold{E},
\partial {\mathbb D})$ over the circle. Denote by
$\widehat{\mathfrak F}$ the \index{$\widehat{\mathfrak F}$}
isotopy class of fiber bundles over $\partial {\mathbb D}$ that
contains ${\mathfrak F}$. Let $A_{r,R} = \{ z \in {\mathbb C} : r <
\vert z \vert < R \}$, $r < 1 < R$, be an annulus  containing the
unit circle. A fiber bundle on $A_{r,R}$ is said to represent
$\widehat{\mathfrak F}$ if its restriction to the unit circle
$\partial {\mathbb D}$ is an element of $\widehat{\mathfrak F}$.

\begin{defn}\label{defEl3}
{\rm (The conformal module of isotopy classes of $({\sf g},{\sf m})$-fiber
bundles.)} Let $\widehat{\mathfrak F} = \widehat{\mathfrak F}_{{\sf g},{\sf m}}$ be
the isotopy class of an oriented $({\sf g},{\sf m})$-fiber bundle over
the circle $\partial {\mathbb D}$. Its conformal module is defined
as
\begin{align}\label{eqEl8}
{\mathcal M} (\widehat{\mathfrak F}) = \sup \{ m(A_{r,R}) :\, & \mbox{there exists a
holomorphic fiber bundle} \nonumber\\ & \mbox{on}\; A_{r,R}\; \mbox{that represents}\;
\widehat{\mathfrak F}\}.
\end{align}
\end{defn}

\section{The monodromy and the mapping torus.}
\label{sec:9.2}

Consider a smooth $(\sf{g},\sf{m})$-bundle
$\mathfrak{F}_{\sf{g},\sf{m}}=(\mathcal{X},\mathcal{P},\mathbold{E},\partial{\mathbb{D}})$ over the unit
circle $\partial{\mathbb{D}}$. Denote as before the set of distinguished points $\mathbold{E}\cap \mathcal{P}^{-1}(x)$ in the fiber
over $x$ by $E_x$.
Let $v$ be the unit tangent vector field to $\partial
{\mathbb D}$. The argument used for the proof of Ehresmann's
Theorem provides a smooth vector field
$V$ on ${\mathcal X}$ which is  tangent to $\mathbold{E}$ at points of this set
and projects to $v$, i.e. $(d \, {\mathcal P})(V) = v$. Cover $\partial {\mathbb
D}$ by its universal covering ${\mathbb R}$, using the mapping $t
\to e^{2\pi i t}$, $t \in {\mathbb R}$. Lift the bundle ${\mathfrak
F}_{\sf{g},\sf{m}}$ over $\partial {\mathbb D}$ to a bundle $\widetilde{\mathfrak
F}_{\sf{g},\sf{m}} = (\widetilde{\mathcal X} , \widetilde{\mathcal P} , \widetilde{\mathbold{E}}, {\mathbb
R})$ over ${\mathbb R}$.

The fiber $\widetilde{\mathcal{P}}^{-1}(t)$ with set of distinguished points $\tilde{E}_t$ equals
$\mathcal{P}^{-1}(e^{2\pi i t})$ with set of distinguished points $E_{e^{2\pi i t}}$, hence for $k
\in \mathbb{Z}$ and each $t \in \mathbb{R}$ the sets
$\widetilde{\mathcal{P}}^{-1}(t+2\pi k)$ and   $\widetilde{\mathcal{P}}^{-1}(t)$,
are equal, and $\tilde{E}_{t+2\pi k}$ is equal to $\tilde{E}_{t}$.
Lift the vector field $V$ to a
vector field $\widetilde V$ on $\widetilde{\mathcal{X}}$ with the following property.
$\widetilde V(z_1) =\widetilde V(z_2) $, if $z_1, z_2 \in  \widetilde{\mathcal{X}}$ are
mapped to the same point $z  \in \mathcal{X}$ under the projection $\widetilde{\mathcal{X}}\overset{{p}}{-\!\!\!\longrightarrow}
{\mathcal{X}}$.

For $t \in \mathbb{R},\, \zeta \in \widetilde{\mathcal{P}}^{-1}(t_1)$, we let $\widetilde\varphi_t (\zeta)\in \widetilde{\mathcal{X}}$  be
the solution of the differential equation
\begin{equation}\label{eqEl1}
\frac\partial{\partial t} \, \widetilde\varphi_t(\zeta) = \widetilde
V (\widetilde\varphi_t (\zeta)),\; \widetilde\varphi_0(\zeta) = \zeta  .
\end{equation}
Put $S \stackrel{def}{=}{\mathcal P}^{-1} (1)\cong \widetilde{\mathcal P}^{-1} (0)$ and $E\stackrel{def}= \mathbold{E}\cap  {\mathcal P}^{-1} (1) \subset S$. Let
$\zeta \in S$. Then $\widetilde \varphi_t(\zeta) \in  \widetilde{\mathcal P}^{-1} (t)$. The
time $t$ map $\widetilde\varphi_t (\zeta)$, $\zeta \in S ,$
of the vector field $\widetilde V$ defines a homeomorphism from the
fiber $S \cong \widetilde{\mathcal P}^{-1} (0)$ onto the fiber
$\widetilde{\mathcal P}^{-1} (t)$, that maps the set of distinguished points $E\subset S$ to the set of distinguished points $\widetilde{E }_t\subset \widetilde{\mathcal P}^{-1} (t)$. In the same way $\widetilde \varphi_t$ defines a
homeomorphism from the fiber $\widetilde{\mathcal P}^{-1} (t_1)$ onto the fiber
$\widetilde{\mathcal P}^{-1} (t_1+t)$ that maps the distinguished points $\tilde{E }_{t_1}$ to the distinguished points $\tilde{E }_{t_1+t}$. The mappings $\widetilde \varphi_t$ form a group:
\begin{equation}\label{eqEl2}
\widetilde \varphi_{t_1+t_2}(\zeta) = \widetilde \varphi_{t_1}(\widetilde \varphi_{t_2}(\zeta)),
\;  \zeta \in \widetilde{\mathcal{X}}, \; t_1, t_2 \in \mathbb{R}.
\end{equation}
Let $\big(\mathbb{R} \times S ,\, {\rm pr}_1\,, \mathbb{R}\times E,\, \mathbb{R}\big)$ be the trivial bundle.
Here ${\rm pr}_1:\mathbb{R} \times S \to \mathbb{R}$ is the projection onto the first factor.
The mapping $\phi$,
\begin{equation}\label{eqEl3}
\mathbb{R}\times S \ni (t,\zeta) \to \Phi(t,\zeta)\stackrel{def}{=}
\widetilde{\varphi}_t(\zeta)\in \widetilde{\mathcal{X}}
\end{equation}
provides an isomorphism from the trivial bundle
$\big(\mathbb{R} \times S ,\, {\rm pr}_1\,, \mathbb{R}\times E,\, \mathbb{R}\big)$
to the lifted
bundle $\widetilde{\mathfrak F}_{\sf{g},\sf{m}}$. Define $\varphi_t$ by the projection
\begin{equation}\label{eqEl5}
\tilde{\mathcal{X}}\supset \widetilde{\mathcal P}^{-1} (t)
\ni \widetilde{\varphi}_t(\zeta) \overset{{p}}{-\!\!\!\longrightarrow} \varphi_t(\zeta)= p(\widetilde{\varphi}_t(\zeta))\in {\mathcal P}^{-1}(e^{2\pi i t}) \in \mathcal{X}\,.
\end{equation}
Since $\widetilde{\mathcal{P}}^{-1}(t)={\mathcal{P}}^{-1}(e^{2\pi i t}) $
and $\tilde{E}_t=E_{e^{{2\pi i t}}}$,
we obtain a smooth family
of homeomorphisms $\varphi_t : {\mathcal P}^{-1} (1) \to {\mathcal
P}^{-1} (e^{2\pi it})$ that map distinguished points to distinguished points.
We call the family $\varphi_t, \,t \in \mathbb{R}, $ a trivializing family of homeomorphisms.
Consider the time-$1$ map
\begin{equation}\label{eqEl4}
\varphi_1 : {\mathcal P}^{-1} (1) \to {\mathcal P}^{-1} (e^{2\pi i})
= {\mathcal P}^{-1} (1) \,
\end{equation}
which is a  self-homeomorphism of the fiber over $1$ that maps the set of distinguished points $E\subset S\cong {\mathcal P}^{-1} (1)$ to itself.
For $k \in \mathbb{Z}$ we have $\varphi_{t+k}=\varphi_t \circ \varphi_1^k.$
For each $n \in \mathbb{Z}$ the projection \eqref{eqEl5}
maps the points $\Phi(t,\zeta)=\widetilde{\varphi}_t(\zeta)$ and $\Phi(t+n,
\widetilde{\varphi}_1^{-n}(\zeta))=\widetilde{\varphi}_{t+n}(\varphi_1^{-n}(\zeta))$
to the same point $ \varphi_t(\zeta) \in {\mathcal P}^{-1} (e^{2\pi i t})$.
The group ${\mathbb Z}$ of integer numbers acts on $\Phi({\mathbb R}
\times S) = \widetilde{\mathcal X}$ by
\begin{equation}\label{eqEl6}
\Phi({\mathbb R} \times S) \ni \Phi (t,\zeta) \to
\Phi(t+n,\varphi_1^{-n}(\zeta)) \, , \quad n \in {\mathbb Z} \, .
\end{equation}
The quotient
\begin{equation}\label{eqEl6a}
{\mathcal X}_1 \stackrel{def}{=} {\mathbb R} \times S \diagup \big((t,\zeta) \sim (t+1,
\varphi_1^{-1}(\zeta))\big)
\end{equation}
is a smooth manifold that is diffeomorphic to $\mathcal{X}$.
Indeed, the mapping ${{p}}\circ\Phi: \mathbb{R}\times S\to \mathcal{X}$,
${{p}}\circ\Phi(t,\zeta)= \varphi_t(\zeta)\in {\mathcal P}^{-1} (e^{2\pi i t})$, satisfies the condition ${{p}}\circ\Phi\big(t+1,\varphi_1^{-1}(\zeta)\big)={{p}}\circ\Phi(t,\zeta)$. Hence, this mapping descends to a mapping from $\mathcal{X}_1$ to $\mathcal{X}$, that is a local diffeomorphism and is bijective by the construction. Hence, it is a diffeomorphism from $\mathcal{X}_1$ onto $\mathcal{X}$.
Put
\begin{equation}\label{eqEl6b}
\mathbold{E}_1\stackrel{def}= \mathbb{R}\times E\diagup \Big((t,\zeta) \sim (t+1,
\varphi_1^{-1}(\zeta))\Big)\,,
\end{equation}
and recall that ${p}(\Phi(t,E))={ p}(\widetilde{\varphi}_t(E))= \mathbold{E}\cap \mathcal{P}^{-1}(e^{2\pi i t})$, since $\widetilde{\varphi}_t$ maps the distinguished points $E\subset S= \mathcal{P}^{-1}(1)$ to the distinguished points $\mathbold{E}\cap \mathcal{P}^{-1}(e^{2\pi i t})$ in the fiber $\widetilde{\mathcal{P}}^{-1}(t)=\mathcal{P}^{-1}(e^{2\pi i t})$.
Define the projection $\mathcal{P}_1:{\mathcal X}_1 \to
\partial \mathbb{D}$, so that $\mathcal{P}_1$ takes the value $e^{2\pi it}$ on the class
containing $(t,\zeta)$. We obtain a bundle $({\mathcal X}_1,\mathcal{P}_1,\mathbold{E}_1, \partial\mathbb{D})$ that is smoothly isomorphic to
$({\mathcal X},\mathcal{P},\mathbold{E}, \partial \mathbb{D})$.
We will also say for short that
${\mathfrak F}_{\sf{g},\sf{m}}$ is smoothly isomorphic to the
mapping torus
\begin{equation}\label{eqEl7}
\left([0,1] \times S\right) \diagup \Big((0 , \zeta) \sim
(1,\varphi (\zeta)) \Big)\,
\end{equation}
where $\varphi = \varphi_1^{-1}$ is a smooth orientation preserving
self-homeomorphism of the fiber $S=\mathcal{P}^{-1}(1)$
with set of distinguished points $E\subset S$.
The mapping $\varphi$ depends
on the trivializing vector field. However, the mapping class of $\varphi$ is independent on
this vector field, it is merely determined by the bundle.
We denote it by $\mathfrak{m}_{\mathfrak{F}}$ (with $\mathfrak{F}={\mathfrak
F}_{\sf{g},\sf{m}}$).
The mapping class $\mathfrak{m}_{\mathfrak{F}} \in \mathfrak{M}(S;\emptyset, E)$
is called the monodromy mapping class
of the bundle ${\mathfrak F}$ over the circle, or the monodromy, for short.

We will consider now the relation between isotopy classes and isomorphism classes of $({\sf g,m})$-bundles, and their relation to the monodromies of the bundles.

Isotopic bundles are isomorphic. The latter fact can be proved as follows.
Let $I \supset [0,1]$ be an open interval, and let   $(\mathcal{Y},\mathcal{P},\mathbold{E},X\times I)$ be a smooth $\sf(g,m)$-bundle over $X\times I$ such that the restrictions to $\{0\}\times X$ and to $\{1\}\times X$ are equal to given smooth $\sf(g,m)$-bundles over $X$. Here $X$ is a smooth finite open surface.
Consider the vector field $v$ on $I\times X$ that equals the unit vector in positive direction of $I$ at each point of $I\times X$, and let $V$ be a smooth vector field on $\mathcal{Y}$ that projects to $v$ under $\mathcal{P}$ and is tangent to ${\mathbold{E}}$  at points of this set. For each $t\in I$ we let $\varphi_t\in \mathcal{Y}$ be the time $t$ map of the flow of $V$, more precisely,
\begin{align*}
\frac{\partial}{\partial t}  \varphi_t(y)=V(\varphi_t(y)),\;\; \varphi_0(y)=y, \;\;\; t\in I,\; y \in \mathcal{P}^{-1}(\{0\}\times X)\,.
\end{align*}
Then $ \varphi_1$ is a diffeomorphism from $\mathcal{P}^{-1}(\{0\}\times X)$ onto $\mathcal{P}^{-1}(\{1\}\times X)$, that maps the fiber over $(0,x)$ onto the fiber over $(1,x)$ and maps distinguished points to distinguished points.

Two isomorphic bundles $\mathfrak{F}_0=(\mathcal{X}_0,\mathcal{P}_0,\mathbold{E}_0,\partial \mathbb{D})$ and $\mathfrak{F}_1=(\mathcal{X}_1,\mathcal{P}_1,\mathbold{E}_1,\partial \mathbb{D})$  over the circle with the same fiber $S$ and the same set of distinguished points $E$ over the base point have conjugate monodromy mapping classes. Indeed, let $\varphi:\mathcal{X}_0\to \mathcal{X}_1$ be a diffeomorphism that maps each fiber $\mathcal{P}_0^{-1}(e^{2\pi i t})$ of the first bundle onto the fiber $\mathcal{P}_1^{-1}(e^{2\pi i t})$ of the second bundle, and let $\varphi_t:S\to \mathcal{P}_0^{-1}(e^{2\pi i t}) $ be a trivialising family of homeomorphisms for the bundle $\mathfrak{F}_0$. Then $\psi_t\stackrel{def}= \varphi|\mathcal{P}_0^{-1}(e^{2\pi i t}) \circ \varphi_t \circ (\varphi|\mathcal{P}_0^{-1}(1))^{-1}$ is a trivialising family of homeomorphisms of the second bundle, and $\psi_1\in {\rm Hom}^+(S;\emptyset,E)$ is conjugate to $\varphi_1\in {\rm Hom}^+(S;\emptyset,E)$. Let $\mathfrak{m}_0\in \mathfrak{M}(S;\emptyset,E)$ be the mapping class represented by $\varphi_1 $, $\mathfrak{m}_1\in \mathfrak{M}(S;\emptyset,E)$ the mapping class represented by $\psi_1$, and $\mathfrak{m}\in \mathfrak{M}(S;\emptyset,E)$ the mapping class of $\varphi|\mathcal{P}_0^{-1}(1)$. Then $\mathfrak{m}_1=\mathfrak{m} \mathfrak{m}_0 \mathfrak{m}^{-1}$, i.e. the monodromy mapping classes are conjugate.

For the case when $\mathfrak{F}_0$ and $\mathfrak{F}_1$ are allowed to have different fibers $S_1$ and $S_2$ over the base point and distinguished points $E_1$ and $E_2$, we notice that a diffeomorphism $\varphi:S_1\to S_2$ with $\varphi(E_1)=\varphi(E_2)$ induces an isomorphism between the mapping class groups $\mathfrak{M}(S_1;\emptyset,E_1)$ and
$\mathfrak{M}(S_2;\emptyset,E_2)$.
The isomorphism between the mapping class groups is defined up to conjugation.
Identifying $\mathfrak{M}(S_1;\emptyset,E_1)$ and
$\mathfrak{M}(S_2;\emptyset,E_2)$ by  a diffeomorphism $\varphi:S_1\to S_2$,
we see as before that the bundles have conjugate monodromies.
The general case, when the bundle isomorphism is defined by two diffeomorphims $\Phi:\mathcal{X}_0\to\mathcal{X}_1$ and $\phi:\partial {\mathbb{D}}\to X_1$ follows easily.
We see that smoothly isomorphic $({\sf{g}},{\sf m})$-bundles over the circle have conjugate monodromy mapping classes.

Vice versa, suppose two $(\sf{g,m})$-bundles $\mathfrak{F}_0$ and $\mathfrak{F}_1$ over the circle have the same fiber $S$ over the base point $1$ and the same set of distinguished points $E\subset S$. Suppose the monodromies of both bundles are equal to the same mapping class $ \mathfrak{m}$. Each of the two bundles is isomorphic to a mapping torus corresponding to the surface $S$ and
a mapping $\varphi_j$ representing $ \mathfrak{m}$ (i.e. to the bundle over the circle defined by \eqref{eqEl6a} and \eqref{eqEl6b}).
The mapping tori of isotopic mappings are isotopic, and, hence, isomorphic.
We showed that two bundles over the circle with equal fiber and equal monodromy are isomorphic.

Suppose the monodromies $ \mathfrak{m}_0$ and $ \mathfrak{m}_1$ of two bundles $\mathfrak{F}_0$ and $\mathfrak{F}_1$ over the circle  with equal fiber $S$ over the base point and equal set of distinguished points $E$ in $S$ are conjugate, i.e.  $\mathfrak{m}_1=\mathfrak{m} \mathfrak{m}_0 \mathfrak{m}^{-1}$. Let $\varphi\in \mathfrak{m}$ and $\varphi_0\in \mathfrak{m}_0$. Then $\mathfrak{F}_0$ is isomorphic to
the mapping torus of $\varphi_0$. This mapping torus has total space $\mathcal{X}_0=\big(\mathbb{R}\times S\big)\diagup \big((t,\zeta)\sim (t+1,\varphi_0(\zeta))\big)$ and set of distinguished points
$E_0=\big(\mathbb{R}\times E\big)\diagup \big((t,\zeta)\sim (t+1,\varphi_0(\zeta))\big)$. The mapping torus $\mathfrak{F}_0'$ with total space
$\mathcal{X}_0'=\big(\mathbb{R}\times \varphi(S)\big)\diagup \big((t,\zeta)\sim (t+1,\varphi\circ \varphi_0\circ \varphi^{-1} (\zeta))\big)$ and set of distinguished points $E_0'=\big(\mathbb{R}\times \varphi(E)\big)\diagup \big((t,\zeta)\sim (t+1,\varphi \circ \varphi_0 \circ \varphi^{-1} (\zeta))\big)$ is isomorphic to $\mathfrak{F}_0$ and has monodromy map $\mathfrak{m} \mathfrak{m}_0 \mathfrak{m}^{-1}$. Indeed, the mapping 
$\mathbb{R}\times S \ni(t,\zeta)\to (t,\varphi(\zeta)) \in \mathbb{R}\times S$ descends to
a diffeomorphism $\varphi:\mathcal{X}_0\to \mathcal{X}_0'$ that maps fibers to fibers.
The mapping torus $\mathfrak{F}_0'$ with total space $\mathcal{X}_0'$ is isomorphic to $\mathfrak{F}_1$, because the monodromy mappings coincide. Since the mapping torus $\mathfrak{F}_0'$ is also isomorphic to $\mathfrak{F}_0$ , the mapping tori $\mathfrak{F}_0$ and $\mathfrak{F}_1$ are isomorphic.

Suppose the bundles have different fibers $S_1$ and $S_2$ and  distinguished points $E_1$ and $E_2$ over the base point $1$.

Fix a diffeomorphism $\varphi:S_1\to S_2\,$, $\varphi(E_1)=E_2\,$. It induces an isomorphism ${\rm Is}:\mathfrak{M}(S_1;\emptyset,E_1)\to \mathfrak{M}(S_2;\emptyset,E_2)$.
Suppose the bundles have (after identifying the mapping class groups in the fiber over the base point by the isomorphism ${\rm Is}$)
conjugate monodromy mappings $\mathfrak{m}_j\in
\mathfrak{M}(S_j;\emptyset,E_j),\, j=1,2,$.
Take representing homeomorphisms $\varphi_j\in \mathfrak{m}_j$.
The mapping $\varphi \varphi_1 \varphi^{-1}$ is in $\mathfrak{M}(S_2;\emptyset,E_2)$ and is conjugate to $\varphi_2$.
By the preceding arguments the mapping torus of $\varphi_2\in \mathfrak{m}_2$ is isomorphic to the mapping torus of $\varphi \varphi_1 \varphi^{-1}$. Since by the same reasoning as above the mapping tori of $\varphi_1$ and $\varphi \varphi_1 \varphi^{-1}$ are isomorphic we obtain that in general two bundles over the circle with conjugate monodromy are isomorphic.

We will use the mapping class group $\mathfrak{M}(S;\emptyset,E)$ on a reference (Riemann) surface of type $(\sf{g},\sf{m})$. It is denoted by ${\rm Mod}(\sf{g},\sf{m})$ and called the modular group.
Since for each mapping class there is a mapping torus with monodromy equal to this mapping class,
we obtain a bijective correspondence between isomorphism classes of smooth $({\sf{g}},{\sf{m}})$-bundles over the circle and conjugacy classes of elements of the modular group
${\rm Mod}(\sf{g},\sf{m})$.

\smallskip

Let now $X$ be a smooth open surface of genus $g$ with $m$ holes with base point $x_0$. Suppose $\mathfrak{F}=(\mathcal{X},\mathcal{P}, \mathbold{E},X)$ is a smooth  $({\sf g},{\sf m})$-bundle over $X$.
Take smooth closed curves in $X$ parameterized by $\gamma_j:[0,1]\to X$ that represent the generators $e_j$ of the fundamental group $\pi_1(X,x_0)$, $j=1,\ldots, 2g+m-1$. Associate to each generator $e_j$ of $\pi_1(X,x_0)$ the monodromy mapping class of the restricted bundle $\mathfrak{F}|{{\gamma_j}} $. This is
a mapping class on the fiber $S=\mathcal{P}^{-1}(x_0)$ over $x_0$ with distinguished points $E=\mathbold{E}\cap \mathcal{P}^{-1}(x_0)$, that depends only on the $e_j$ and on
the bundle.
We obtain a well-defined mapping that associates to each generator of the fundamental group $\pi_1(X,x_0)$ a mapping class. This mapping extends to a homomorphism
from the fundamental group $\pi_1(X,x_0)$ to the mapping class group $\mathfrak{M}(S, \emptyset, E)$.
In the same way as in the case of bundles over the circle we may associate to each isomorphism class of bundles over $X$ a conjugacy class  of homomorphisms from $\pi_1(X,x_0)$ to ${\rm Mod}(\sf{g,m})$.

Vice versa, take a homomorphism $h$ from the fundamental group $\pi_1(X,x_0)$ to the modular group ${\rm Mod}(\sf{g,m})$.
A bundle over $X$
whose monodromy homomorphism coincides with $h$ can be obtained as follows. We may assume that $X$ is a domain on a closed Riemann surface $X^c$.

Let $\alpha_j$, $\beta_j$, $j=1,\ldots g,$ be a system of smooth simple closed curves on $X^c$
with base point $x_0$, that are contained in $X\subset X^c$, do not intersect any of the holes, and represent generators $e_{2j-1,0},\, e_{2j,0}$ of $\pi_1(X^c,x_0)$ (identified with $\pi_1(X,x_0)$ under the embedding $X\to X^c$). We may assume, that any pair among these curves has the only intersection point $x_0$.

If the number of holes ${m}$ is bigger than $1$, we choose smooth simple closed curves $\gamma_k,\, k=1,\ldots, {m}-1,$ with base point $x_0$ in $X$, such that each $\gamma_k$ divides $X^c$, and one connected component of its complement contains $\mathcal{C}_{k}$ and does not contain other holes. The loop $\gamma_k$ surrounds $\mathcal{C}_k$ counterclockwise.
Each $\gamma_k$ represents an element $e_{2g+k}\in \pi_1(X,x_0)$.

The choice may be done so that, when
labeling the rays of the loops emerging from the base point $x_0$ by  $\alpha_j^-,\,\beta_j^-$  $\gamma_k^-$, and the incoming rays by $\alpha_j^+,\,\beta_j^+$  $\gamma_k^+$,
then moving in counterclockwise direction along a small circle around $x_0$ we meet the rays in the order
\begin{align*}
\ldots, \alpha^-_j,\beta^-_j,\alpha^+_j,\beta^+_j,\ldots, \gamma^-_k,\gamma^+_k,\ldots\;.
\end{align*}

The union of the $\alpha_j$, $\beta_j$, and $\gamma_k$ is a bouquet of circles contained in $X$ with base point $x_0$. The collection of the elements $e_j,\, j=1,\ldots,2g+m-1,$  represented by the collection of the loops is
a set of generators $\mathcal{E}$
of $\pi_1(X,x_0)$.

\begin{figure}[h]
\begin{center}
\includegraphics[width=10cm]{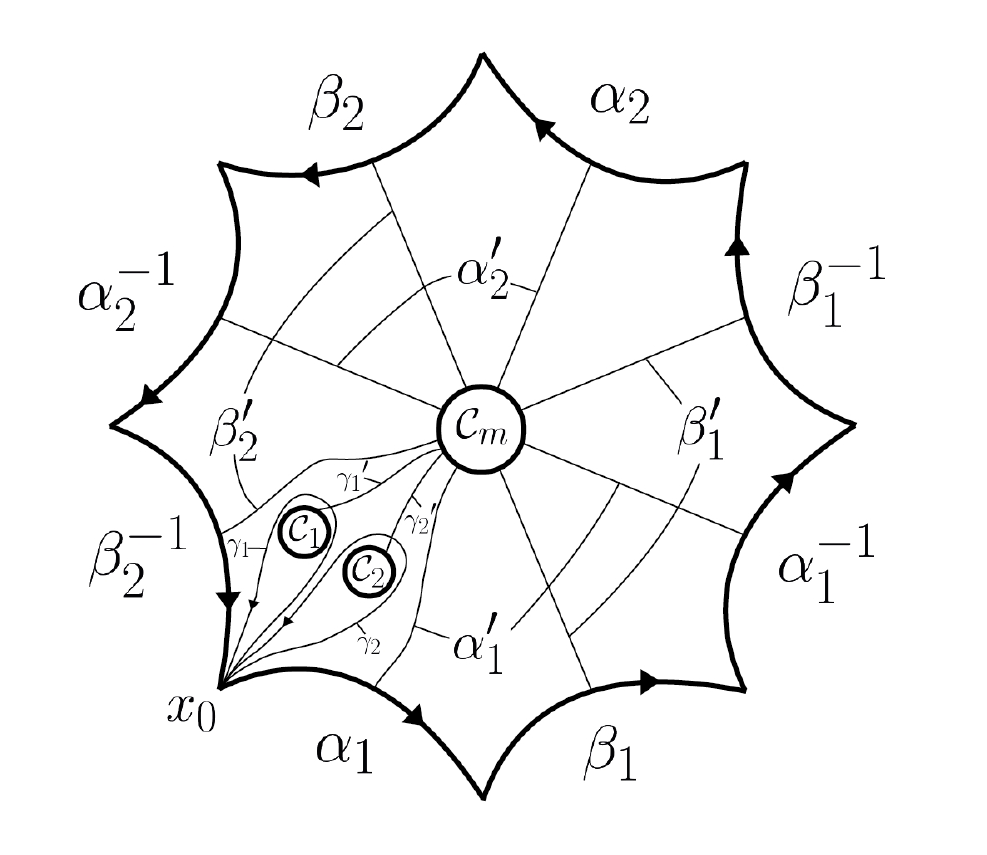}
\end{center}
\caption{Two systems of curves on a Riemann surface of genus $g$ with $m$ holes.}\label{fig5}
\end{figure}

Cut $X^c$ along the union of the $\alpha_j$ and $\beta_j$. We obtain a fundamental polygon $F$ for the universal covering $P:\tilde{X^c}\to X^c$ of $X^c$. See also Figure \ref{fig5}. The projection $P$ takes all vertices of the polygon to the point $x_0\in X\subset X^c$, it takes the sides labeled $\alpha_j$ to the loops $\alpha_j$, $P$ takes
the sides labeled $\alpha_j^{-1}$ (with inverted orientation) to the loops $\alpha_j$, and the respective facts hold for the $\beta_j$.

If the number of holes ${m}$ is bigger than $1$, we take ${ m}-1$ pairwise disjoint simple closed arcs $\gamma'_k,\, k=1,\ldots, { m}-1,$ in $X^c$,
such that $\gamma'_k$ has one endpoint on the boundary $\partial\mathcal{C}_{k}$ of a hole and the other endpoint on the boundary of the hole $\mathcal{C}_{m}$, and the endpoints are the only points of the arcs that are not contained in $X$. We may choose the $\gamma_k'$ so that it intersects $\gamma_k$ once, and does not intersect any other loop among the $\alpha_j$, $\beta_j$, and $\gamma_l$. For each $j=1,\ldots, g,$ we take simple closed arcs $\alpha'_j$ and $\beta'_j$ in $X^c$ with endpoints on $\partial \mathcal{C}_m$, so that $\alpha'_j$ intersects $\alpha_j$ exactly once and $\beta_j$ intersects $\beta'_j$ exactly once, and $\alpha_j$ and $\beta_j$ do not intersect any other of the loops of the bouquet, and all $\alpha'_j$, $\beta'_j$, $\gamma'_k$ are pairwise disjoint.

Let $V_j,\, j=1,\ldots,2{g}+{ m}-1,$ be a system of disjoint relatively closed topological discs in $X$, such that each $V_{2j-1},\,j=1,\ldots,g,$ contains $\alpha'_{j}$ in its interior and does not intersect any other of the $\alpha'$, $\beta'$, $\gamma'$, each $V_{2j},\,j=1,\ldots,g,$ contains $\beta'_{j}$ in its interior and does not intersect any other of the $\alpha'$, $\beta'$, $\gamma'$, each $V_{2j+k},\,k=1,\ldots,m-1,$ contains $\gamma'_{k}$ in its interior and does not other of the $\alpha'$, $\beta'$, $\gamma'$,
and the boundary of $V_j$ in $X$ consists of two smooth arcs, each of whose endpoints is on $\partial \mathcal{C}_j$ or on $\partial \mathcal{C}_m$.
We call each $V_j$ a band (or half-open rectangle) with horizontal sides being the mentioned arcs, and each vertical side being either a point on $\partial \mathcal{C}_j$ or $\partial \mathcal{C}_m$ , or an arc on either $\partial \mathcal{C}_j$ or $\partial \mathcal{C}_m$.
(The horizontal sides are contained in $V_j$, the vertical sides are not.) We may assume that each $V_{2j-1}$ intersects $\alpha_j$ along a connected set and does not intersect any other curve among the $\alpha$'s, $\beta$'s, and $\gamma$'s, each $V_{2j}$ intersects $\beta_j$ along a connected set and does not intersect any other $\alpha$-curve, $\beta$-curve or $\gamma$-curve, and each $V_{2g+k}$ intersects $\gamma_k$ along a connected set and does not intersect any other $\alpha$-curve, $\beta$-curve or $\gamma$-curve.
The set $D\stackrel{def}=X\setminus \cup \bigcup _{j=1}^{2{ g}+{ m}-1} V_j$ is a simply connected domain contained in $X$. Indeed, after removing the $V_{2g+k}$, we obtain a surface of genus $g$ with one hole, denoted by $\mathcal{C}'$. The $\alpha_j'$ and $\beta_j'$ can be extended to loops in $X^c$ with some base point $\mathring{x}\in \mathcal{C}_m$. They form a system of curves in $X^c$ which is dual to the $\alpha_j$ and $\beta_j$. Cutting $X^c$ along this system of loops, we obtain another fundamental domain for the covering $\tilde{X^c}\to X^c$. It is now clear that removing the hole $\mathcal{C}'$
and the $V_j$ from $X^c$ we obtain a disc.

Moreover, $X$ is obtained from the disc $D$ by attaching the ${2g+ m-1}$ bands $V_j$ (each horizontal side of a band is glued to a smooth arc in the boundary of $D$). The bands may be chosen so that the disc $D$ contains the base point $x_0$ of $X$.

Consider the universal covering $\tilde{X} \overset{{ p}}{-\!\!\!\longrightarrow} X$.
Take a point $\tilde{x}_0\in \tilde{X}$ with $p(\tilde{x}_0)=x_0$.
For each $j$ we take the lifts $\widetilde{\gamma}_j$ of $\gamma_j$ to $\tilde X$ with initial point $\tilde{x}_0$, and the lift $\tilde{V}_j$ of $V_j$
that intersects $\widetilde{\gamma}_j$.
Let $\widetilde{D}_0$ be the lift of $D$, that contains $\tilde{x}_0$, and  let $\widetilde{D}_j,\, j=1,\ldots,2{g}+{m}-1, $ be the lift of $D$ to $\tilde U$, that contains the endpoint $\tilde{x}_j\stackrel{def}=\widetilde{\gamma}_j(1)$. Let $\tilde U$ be the domain in $\tilde X$ that is the union of $\tilde{D}_j,\, j=0,1,\ldots, 2g+m-1,$ and $\tilde{V}_j,\, j=1,\ldots,2g+m-1$.

Choose for each $j\geq 1$ a mapping $\varphi_j$ in the mapping class $h(e_j)$, and define a bundle over $U$ by taking the trivial bundle over $\tilde U$ and
making the following identifications. For each $x'\in D$ and each $j$ we glue the  fiber over the point $\tilde{x}'_j\in \widetilde{D}_j$ for which $p(\tilde{x}'_j)=x'$ to the fiber over the point $\tilde{x}'_0\in \widetilde{D}_0$ for which $p(\tilde{x}'_0)=x'$  using the mapping $\varphi_j$. We obtained a bundle with the given monodromy homomorphism.

As in the case of bundles over the circle one can see that each bundle is isomorphic to a bundle that is constructed in the just described way, and bundles constructed in this way with equal conjugacy class of monodromy homomorphisms are isomorphic.
We obtain the
following theorem
(see, e.g. \cite{FaMa} for the case of closed fibers of genus $g\geq 2$ and paracompact Hausdorff spaces $X$).
\begin{thm}\label{thmEl1}
Let $X$ be a smooth finite open surface. The set of isomorphism classes \index{isotopy class ! of fiber
bundles}of smooth oriented $({\sf g},{\sf m})$-bundles on $X$ is in
one-to-one correspondence to the set of conjugacy classes of
homomorphisms from the fundamental group $\pi_1 (X,x_0)$ into the
modular group ${\rm Mod}(\sf{g},\sf{m})$.
\end{thm}

\smallskip

We saw that isotopic smooth $({\sf g,m})$-bundles over the circle are smoothly isomorphic.
Vice versa, consider two smooth $({\sf g,m})$-bundles $\mathfrak{F}_0$ and $\mathfrak{F}_1$ over the circle
with fiber and distinguished points over the base point $1$ being equal. If the bundles have the same monodromy, then they are isotopic by an isotopy that fixes the fiber over the base point and the set of distinguished points in this fiber (by a based isotopy, for short). Indeed, denote the fiber over $1$ by $S$ and the set of distinguished points in $S$ by $E$. Lift the bundle $\mathfrak{F}_j,\,j=0,1,$ to a bundle over the horizontal line $\mathbb{R}\times \{j\}\subset \mathbb{R}^2$ by the covering $\mathbb{R}\times\{j\} \ni (t,j)\to e^{2\pi i t}$.
Let $\varepsilon$ be a small positive number.
Restrict each lift to an open interval $ (-\varepsilon,1+\varepsilon)\times\{j\}\supset [0,1]\times\{j\}$, $\varepsilon>0$, and extend the bundles to bundles over the rectangles $(-\varepsilon,1+\varepsilon)\times (j-\varepsilon,j+\varepsilon)$, so that the restriction of the extended bundle $\mathfrak{F}_j'$ to each vertical segment is a product bundle.
Consider the trivial bundle with fiber $S$ and set of distinguished points $E$ over the vertical line segments $\{k\}\times (-\varepsilon,1+\varepsilon)$, $k=0,1$. Take an extension of the trivial bundle over each of the vertical line segments $\{k\}\times (-\varepsilon,1+\varepsilon)$ to an open rectangle $(k-\varepsilon,k+ \varepsilon)\times (-\varepsilon,1+\varepsilon)$  around the segment, so that for $x \in (-\varepsilon,\varepsilon)$ and $y\in (-\varepsilon,1+\varepsilon)$ the fibers and the distinguished points over $(x,y)$ and $(x+1,y)$ are the same, and
the pieces over the four rectangles match to a smooth bundle $\tilde{\mathfrak{F}}_0$ over the neighbourhood $V\stackrel{def}=\Big((-\varepsilon, 1+\varepsilon)\times (-\varepsilon, 1+\varepsilon)\Big) \, \setminus \,\Big((\varepsilon, 1-\varepsilon)\times (\varepsilon, 1-\varepsilon)\Big) $  of the boundary of the square $[0,1]  \times[0,1]$. Let $(0,0)$ be the base point of $V$. The monodromy of the obtained bundle over $V$
is the identity. Hence, the bundle is smoothly isomorphic to the trivial bundle. Let $\varphi$ be a diffeomorphism of the total space of the bundle to the total space of the trivial bundle over $V$, that maps the fiber over each point $x$ to the fiber over $x$ and maps distinguished points to distinguished points.
The trivial bundle over $V$
extends to the trivial bundle over the neighbourhood $(-\varepsilon, 1+\varepsilon)\times (-\varepsilon, 1+\varepsilon)$ of the square. Glue the
trivial bundle over the square $(\frac{\varepsilon}{2}, 1-\frac{\varepsilon}{2})\times (\frac{\varepsilon}{2}, 1-\frac{\varepsilon}{2})$, that is relatively compact in the open unit square, to the bundle $\tilde{\mathfrak{F}}_0$ over $V$ along $V\cap  \big((\frac{\varepsilon}{2}, 1-\frac{\varepsilon}{2})\times (\frac{\varepsilon}{2}, 1-\frac{\varepsilon}{2})\big)$
using the mapping $\varphi^{-1}$. We obtain a smooth bundle over
the neighbourhood $(-\varepsilon, 1+\varepsilon)\times (-\varepsilon, 1+\varepsilon)$ of the unit square. For $(x,y)\in(-\varepsilon, \varepsilon)\times (-\varepsilon,1+\varepsilon)$ we glue the fiber over $(x,y)$ to the fiber over $(x+1,y)$
using the identity mapping. We get a bundle over the product of the circle with the interval $(-\varepsilon, 1+\varepsilon)$, which provides
an isotopy with fixed fiber over the base point and fixed set of distinguished points (a based isotopy for short).

In the same way it can be shown that two bundles over the circle with the same fiber and set of distinguished points over the base point and conjugate monodromy are free isotopic.

If two bundles over the circle have different fibers $S_1$ and $S_2$ and sets of distinguished points $E_1$ and $E_2$, respectively, and (after an isomorphism) their monodromy maps are conjugate, then the bundles are free isotopic. To reduce this fact to the previous statement, we use a smooth bundle over an interval whose fibers over two different points are equal to $S_1$ and $S_2$ with sets of distinguished points $E_1$ and $E_2$, respectively.

In a similar way one can see that smooth $({\sf g,m})$-bundles over a smooth finite open surface $X$ with conjugate monodromy homomorphism (after an isomorphism between mapping classes) are isotopic. We obtain a bijective correspondence between smooth isomorphism classes and free isotopy classes of $({\sf g,m})$-bundles over finite open surfaces.

A smooth $(\sf{g},\sf{m})$-bundle $(\mathcal{X},\mathcal{P}, \mathbold{E},X)$ over a smooth finite open surface $X$ admits a smooth section, i.e. a smooth mapping $X\ni x \to s(x)\in \mathcal{X}\setminus\mathbold{E} $ for which $\mathcal{P}\circ s(x)=x$ for all $x\in X$.
Indeed, let $D$ and $V_j$ be the same sets as before. There is a smooth section over $D\cup \beta_{j^{\pm}}$, where $\beta_{j^{\pm}}$ are the arcs of $\partial D$ to which the horizontal sides of the $V_j$ are glued. Indeed, there is a simply connected neighbourhood of this set, and the bundle is smoothly trivial over any simply connected domain.
Extend the smooth section to a smooth section over a neighbourhood of $D\cup \beta_{j^{\pm}}\cup (\gamma_j\cap V_j) $. Since a neighbourhood of each $V_j$ is simply connected and there is a smooth deformation retraction of a neighbourhood of $V_j$ to a neighbourhood of $\beta_{j^{\pm}}\cup (\gamma_j\cap V_j)$, we obtain a smooth section on the whole surface $X$.

Recall that
an admissible system of
curves on a closed surface $S$ of genus $\sf g$ with set of $\sf m$ distinguished points $E$ is said to reduce a family of mapping classes $\mathfrak{m}_j \in \mathfrak{M}(S; \emptyset, E) $ if it reduces each $\mathfrak{m}_j$.

Similarly, a $({\sf{g}},{\sf{m}})$-bundle over a finite open surface with fiber $S$ over the base point $x_0$ and set of distinguished points $E\subset S$ is called reducible if there is an admissible system of curves in the fiber over the base point
that reduces all monodromy mapping classes simultaneously. Otherwise the bundle is called
irreducible.

Consider any isotopy class of reducible $(\sf{g},\sf{m})$-bundles with fiber $S$ over the
base point $x_0$. Take an admissible system $\mathcal{C}$ of simple closed curves on $S$
that reduces each monodromy mapping class and is maximal in the sense that there is no
strictly larger admissible system with this property. For each monodromy mapping class
$\mathfrak{m}_j$ we choose a representing homeomorphism $\varphi_j$ that fixes $\mathcal{C}$
setwise. Then all homeomorphisms $\varphi_j$ fix $S\setminus \mathcal{C}$. Decompose
$S\setminus \mathcal{C}$ into disjoint open sets $S_k$, each of which being a minimal union
of connected components of $S\setminus \mathcal{C}$ that is invariant under each
$\varphi_j$.
Each $S_k$ is a union of surfaces with holes and possibly with distinguished points and is homeomorphic to a union
$\overset{\circ}{S}_k$ of closed surfaces with a set $E_k$ of finitely many punctures and, possibly, a set $E_k'$ of finitely many distinguished points. The restrictions of the $\varphi_j$ to
each $S_k$ (possibly, with distinguished points) are conjugate to self-homeomorphisms of $\overset{\circ}{S}_k$ (possibly with set of distinguished points $E_k'$).

For each $k$ we denote by $S_k^c$ the (possibly non-connected) closed Riemann surface which is the disjoint union of the closure of the connected components of $\overset{\circ}{S}_k$. The Riemann surface $S_k^c$ is equipped with distinguished points $E_k\cup E'_k$. For each $k$ we obtain
a conjugacy class of homomorphisms from the fundamental group $\pi_1(X,x_0)$
to the modular group of the possibly non-connected Riemann surface $S_k^c$ with set of distinguished points $E_k\cup E'_k$. By Theorem \ref{thmEl1} (more precisely, by its analog for possibly non-connected fibers) we obtain for each $k$ an isotopy class of bundles over
$X$ with fiber a union of closed surfaces with distinguished points. The obtained isotopy classes of bundles are
irreducible and are called irreducible bundle components of the  $(\sf{g},\sf{m})$-bundle
with fiber $S$. (There may be different decompositions into irreducible bundle components.)

\section{$(0,4)$-bundles over genus $1$ surfaces with a hole}
\label{sec:9.3}
The $(0,n)$-bundles over a manifold $X$ are closely related to the separable quasipolynomials on $X$. In this section we will
state a theorem that concerns $(0,4)$-bundles and is related to Theorem \ref{thm8.0}. We
will prove it together with Theorem \ref{thm8.0}.

Let $X$ be a finite open Riemann surface. (Notice that for the smooth version of the following settings it is enough to suppose that $X$ is a finite open smooth surface.)
By a holomorphic (smooth, respectively) $(0,n)$-bundle with a section over $X$ we mean a holomorphic (smooth, respectively) $(0,n+1)$-bundle
$(\mathcal{X},\mathcal{P}, \mathbold{E}, X)$,
such that the complex manifold (smooth oriented manifold, respectively) $\mathbold{E}\subset \mathcal{X}$ is the disjoint union of two complex manifolds (smooth oriented manifolds, respectively) $\mathring{\mathbold{E}}$ and $\mathbold{s}$, where $\mathring{\mathbold{E}}\subset \mathcal{X}$
intersects each fiber $\mathcal{P}^{-1}(x)$ along a set $\mathring{E}_x$ of $n$ points, and $\mathbold{s}\subset \mathcal{X}$ intersects each fiber $\mathcal{P}^{-1}(x)$ along a
single point $s_x$. In other words, the mapping $x\to s_x,\, x\in X$, is a (holomorphic, smooth, respectively) section of the $(0,n)$-bundle with set of distinguished points $\mathring{E}_x$ in the fiber over $x$. Two smooth $(0,n)$-bundles with a section are called
isomorphic if they are isomorphic as smooth $(0,n+1)$-bundles and the diffeomorphism between the total spaces takes the section of one bundle to the section of the other bundle.
Two smooth $(0,n)$-bundles with a section are called isotopic if they are isotopic as $(0,n+1)$-bundles with an isotopy that joins the sections of the bundles.
Two smooth $(0,n)$-bundles with a section are isotopic iff they are isomorphic.

A special smooth (holomorphic, respectively) $(0,n+1)$-bundle over $X$ is a smooth (holomorphic, respectively) bundle over $X$ of the form $(X\times \mathbb{P}^1, {\rm pr}_1, \mathbold{E}, X)$,
where ${\rm pr}_1 : X\times \mathbb{P}^1 \to X$ is the projection onto the first factor, and the smooth submanifold $\mathbold{E}$ of $X\times \mathbb{P}^1$  is equal to the disjoint union
$\mathring{\mathbold{E}}\cup \mathbold{s}^{\infty}$,
where $\mathring{\mathbold{E}}$ is a smooth submanifold of $X\times \mathbb{P}^1$ that intersects each fiber along $n$ finite points and $\mathbold{s}^{\infty}$ is the smooth submanifold of $X\times \mathbb{P}^1$ that intersects each fiber $\{x\}\times \mathbb{P}^1$ along the point $\{x\}\times \{\infty\}$. A special $(0,n+1)$-bundle is, in particular, a $(0,n)$-bundle with a section.

Each smooth mapping $f:X \to C_n ({\mathbb C}) \diagup
{\mathcal S}_n$ from $X$ to the symmetrized configuration space defines a smooth special $(0,n+1)$-bundle over $X$. The smooth submanifold $\mathbold{E}$ of $X\times \mathbb{P}^1$  is equal to the union $\bigcup_{x \in X}\,\big(x,f(x)\cup\{\infty\}\big)$. Vice versa, for each special $(0,n+1)$-bundle the mapping $X\ni x\to \mathring{\mathbold{E}}\cap ({\rm pr}_1)^{-1}(x)$ defines a smooth separable quasi-polynomial $f$ of degree $n$ on $X$.
The special $(0,n+1)$-bundle is holomorphic if and only if the mapping $f$ is
holomorphic.

Choose a base point $x_0\in X$. Denote by $f_*$ the mapping from  $\pi_1(X,x_0)$ to the fundamental group $\pi_1(C_n ({\mathbb C}) \diagup {\mathcal S}_n, f(x_0))\cong \mathcal{B}_n$ induced by $f$.
Two smooth mappings $f_1$ and $f_2$ from $X$ to $C_n ({\mathbb C}) \diagup {\mathcal S}_n$
define special $(0,n+1)$-bundles that are isotopic through special $(0,n+1)$-bundles
if and only if the quasipolynomials $f_1$ and $f_2$ are free isotopic.
They define bundles that are isomorphic (equivalently, they are isotopic through $(0,n)$-bundles with a section) if and only if for a set of generators $e_j$ of the fundamental group $\pi_1(X,x_0)$, for a braid $w\in \mathcal{B}_n$ and integer numbers $k_j$ the equalities
$(f_1)_*(e_j)=w^{-1}(f_2)_*(e_j)w\Delta_n^{2k_j}$ hold.
Indeed, the bundles are isomorphic iff their monodromy homomorphisms are conjugate.
The monodromy mapping classes of the bundles are elements of $\mathfrak{M}(\mathbb{P}^1; \{\infty\},f(x_0))$ and
the braid group on $n$ strands modulo its center $\mathcal{B}_n\diagup \mathcal{Z}_n$
is isomorphic to the mapping class group $\mathfrak{M}(\mathbb{P}^1;\{\infty\},f(x_0))$.

Each smooth $(0,n)$-bundle $\mathfrak{F}= (\mathcal{X},\mathcal{P},\mathring{\mathbold{E}}\cup \mathbold{s},X)$ with a section
over a smooth finite open surface $X$ is isomorphic (equivalently, free isotopic) to a smooth special $(0,n+1)$-bundle. Indeed,
put $\mathring{E}=\mathring{\mathbold{E}}\cap\mathcal{P}^{-1}(x_0)$, $s=\mathbold{s}\cap\mathcal{P}^{-1}(x_0)$ for the base point $x_0\in X$.
A homeomorphism $\varphi: \mathcal{P}^{-1}(x_0)\to \{x_0\}\times\mathbb{P}^1$ with $\varphi(s)=(x_0,\infty)$ and $\varphi(\mathring{E})= \{x_0\}\times \mathring {E}'$ for a set
$\mathring{E}'\subset C_n(\mathbb{C})\diagup \mathcal{C}_n$ induces
an isomorphism ${\rm Is}:\mathfrak{M}(\mathcal{P}^{-1}(x_0); s, \mathring{E})  \to\mathfrak{M}(\mathbb{P}^1;\{\infty\},\mathring {E}')$. The bundle $\mathfrak{F}$ corresponds to a conjugacy class of homomorphisms $\pi_1(X,x_0)\to \mathfrak{M}(\mathcal{P}^{-1}(x_0); s, \mathring{E}) $. The isomorphism between mapping class groups gives us a conjugacy class of homomorphisms $\pi_1(X,x_0)\to \mathfrak{M}(\mathbb{P}^1;\{\infty\},\mathring {E}')$. There is a special $(0,n+1)$-bundle that corresponds to the latter conjugacy class of homomorphisms.
This bundle is isomorphic to $\mathfrak{F}$.

\begin{lemm}\label{lem*}
Each holomorphic $(0,n)$-bundle $\mathfrak{F}\stackrel{def}=(\mathcal{X},\mathcal{P},\mathring{\mathbold{E}}
\cup\mathbold{s},X)$ with a section over a finite open Riemann surface $X$ is holomorphically isomorphic
to a special holomorphic $(0,n+1)$-bundle.
\end{lemm}

\noindent {\bf Proof}.
Let $D$, $\tilde{D}_j,\, j=0,\ldots,$ and $\tilde U$ be as in Section \ref{sec:9.2}. Lift the bundle $\mathfrak{F}$ to a bundle on the universal covering $\tilde X$ and restrict the lift to a bundle $\tilde{\mathfrak{F}}= (\tilde{\mathcal{P}}^{-1}(\tilde{U}),\tilde{\mathcal{P}},\widetilde{\mathring{\mathbold{E}}}
\cup\widetilde{\mathbold{s}},\tilde{U})$  on $\tilde U$. The lift $\widetilde{\mathring{\mathbold{E}} }$ of ${\mathring{\mathbold{E}} }$
to $\tilde U$ is the union of $n$ connected components $\widetilde{\mathring{\mathbold{E}} }_j,\,j=1,\ldots,n,$ each intersecting each fiber along a single point.

Consider the holomorphic $(0,0)$-bundle $(\tilde{\mathcal{P}}^{-1}(\tilde{U}),\tilde{\mathcal{P}},\tilde{U})$
that is obtained from $\tilde{\mathfrak{F}}$
by forgetting the
distinguished points and the section. All fibers of the $(0,0)$-bundle $(\tilde{\mathcal{P}}^{-1}(\tilde{U}),\tilde{\mathcal{P}},\tilde{U})$ are conformally equivalent compact Riemann surfaces.
Hence, by a theorem of Fischer and Grauert \cite{FiGr} the bundle is locally holomorphically
trivial. This means that for each point $\tilde{x}\in \tilde{U}$ there exists a simply connected open subset $U_{\tilde{x}} \subset \tilde{U}$
containing $\tilde{x}$, and a biholomorphic map $\varphi_{U_{\tilde{x}}}:(\tilde{\mathcal{P}})^{-1}(U_{\tilde{x}})\to U_{\tilde{x}}\times
\mathbb{P}^1$ such that the diagram
\[
\xymatrix{
{\tilde{\mathcal{P}}^{-1}(U_{\tilde{x}})}      \ar[r]^{\varphi_{U_{\tilde{x}}}     } \ar[d]_{\tilde{\mathcal{P}}         } &
{U_{\tilde{x}} \times \mathbb{P}^1}   \ar[dl]^{{\rm  pr}_1} \\
  U_{\tilde{x}}
}
\]
commutes. For each $\tilde{x}'\in U_{\tilde{x}}$ and $j=1,\ldots,n,$ the mapping
$\varphi_{U_{\tilde{x}}}$
takes the point
$\widetilde{\mathring{\mathbold{E}}}_j \cap (\tilde{\mathcal{P}})^{-1}(\tilde{x}')$,
to a point denoted by $(\tilde{x}', \zeta_j(\tilde{x}'))$, and takes the point of the section $\mathbold{s} \cap (\tilde{\mathcal{P}})^{-1}(\tilde{x}')$ to a point denoted by $(\tilde{x}',\zeta_{n+1}(\tilde{x}'))$. The $\zeta_j,\, j=1,\ldots,n+1,$ are holomorphic.
Composing if necessary $\varphi_{U_{\tilde{x}}}$ with a holomorphic mapping from $U_{\tilde{x}} \times \mathbb{P}^1$ onto itself which preserves each fiber, we may assume that
$\varphi_{U_{\tilde{x}}}$ maps no point in $\tilde{\mathcal{P}}^{-1}(\tilde{x})\cap (\widetilde{\mathring{\mathbold{E}}} \cap \mathbold{s})$ to $(\tilde{x},\infty)$.

Maybe, after shrinking the $U_{\tilde x}$, the mapping
$$
U_{\tilde x}\times \mathbb{P}^1\ni (\tilde{x}',\zeta)\to w_{\tilde x}(\tilde{x}',\zeta)=\Big(\tilde{x}', -1+2\frac{(\zeta_2(\tilde{x}')-\zeta)( \zeta_{n+1}(\tilde{x}')-\zeta_1(\tilde{x}'))          }{(\zeta_2(\tilde{x}')- \zeta_{1}(\tilde{x}')) (\zeta_{n+1}(\tilde{x}')-\zeta)}\Big)
$$
provides a biholomorphic map $ U_{\tilde x}\times \mathbb{P}^1\toitself$, that maps for each $\tilde{x}'\in U_{\tilde{x}}$ the point $\zeta_1(\tilde{x}')$ to $1$, it maps $\zeta_2(\tilde{x}')$ to $-1$, and $\zeta_{n+1}(\tilde{x}')$ to $\infty$.

Denote by $\psi_{\tilde x}$ the mapping $w_{\tilde x}\circ \varphi_{U_{\tilde{x}}}$
from  $\tilde{\mathcal{ P}}^{-1}(U_{\tilde{x}})$ onto $U_{\tilde{x}}\times \mathbb{P}^1$.
If for two points $\tilde{x},\tilde{y}\in \tilde{U}$ the intersection $U_{\tilde{x}}\cap U_{\tilde{y}}$ is not empty and connected, then the mappings $\psi_{\tilde{ x}}$ and $\psi_{\tilde{ y}}$ coincide on $\tilde{\mathcal{ P}}_1^{-1}(U_{\tilde{x}}\cap U_{\tilde{y}})$.
We get a holomorphic isomorphism $\tilde{\psi}$ from the total space $\tilde{\mathcal{U}}$ of the bundle $\tilde{\mathfrak{F}}$
onto $\tilde{U}\times \mathbb{P}^1$
that maps the first component of $\widetilde{\mathring{\mathbold{E}}}$ to $\tilde{U}\times \{-1\}$, the second to $\tilde{U}\times \{1\}$, and the section to  $\tilde{U}\times \{\infty\}$. The isomorphism takes the set $\widetilde{\mathring{\mathbold{E}}}$ to a set $\widetilde{\mathring{\mathbold{E'}}}\subset \tilde{U}\times \mathbb{P}^1$, and the section $\widetilde{\mathbold{s}}$ to the set
$\widetilde{\mathbold{s}}^{\infty}=\tilde{U}\times \{\infty\}$.
The identity mapping from $\tilde U$ onto itself together with the mapping $\tilde{\psi}$ define
a holomorphic isomorphism from the bundle $\tilde{\mathfrak{F}}$ onto the bundle $(\tilde{U}\times \mathbb{P}^1,\,{\rm pr}_1,\,\widetilde{\mathring{\mathbold{E'}}}\cup \widetilde{\mathbold{s}}^{\infty}\,, \tilde{U})$.
Recall that the original bundle $\mathfrak{F}$ is obtained from the bundle $\tilde{\mathfrak{F}}$ by gluing for each $j$ the restriction $\tilde{\mathfrak{F}}|\tilde{D}_j$ to the restriction $\tilde{\mathfrak{F}}|\tilde{D}_0$ using a holomorphic isomorphism
$\varphi_j$ from $\tilde{\mathcal{P}}^{-1}(\tilde{D}_0)$ onto $\tilde{\mathcal{P}}^{-1}(\tilde{D}_j)$. A bundle that is holomorphically isomorphic to $\mathfrak{F}$ is obtained as follows. Take the bundle $(\tilde{U}\times \mathbb{P}^1,\,{\rm pr}_1,\,\widetilde{\mathring{\mathbold{E'}}}\cup \widetilde{\mathbold{s}}^{\infty},\, \tilde{U})$ and glue 
for each $j=1,\ldots,n,$ its restrictions to $\tilde{D}_j$ and  $\tilde{D}_0$
together
using the
biholomorphic mapping $\tilde{\psi}|\tilde{\mathcal{P}}^{-1}(\tilde{D}_j))\circ \varphi_j\circ (\tilde{\psi}|\tilde{\mathcal{P}}^{-1}(\tilde{D}_0))^{-1}$ between subsets of the total space of the bundle  $(\tilde{U}\times \mathbb{P}^1,\,{\rm pr}_1,\,\widetilde{\mathring{\mathbold{E'}}}\cup \widetilde{\mathbold{s}}^{\infty},\, \tilde{U})$.
The last mapping restricts to each fiber as a conformal mapping that maps the point of the section in one fiber to the point of the section in the other fiber. Hence, on each fiber $\{x\}\times \mathbb{P}^1$
we get a mapping $(x,\zeta)\to (x,a_j(x) +\alpha_j(x)\zeta)$ that restricts to $\{x\}\times\mathbb{C}$ as a complex affine mapping in the second coordinate. Here $x=p(\tilde{x}_j)=p(\tilde{x}_0)$ for the projection $p$ from $\tilde X$ to $X$. The mappings $x\to a_j(x)$ and $x\to \alpha_j(x)$ depend holomorphically on $x\in D$, and $\alpha_j$ is nowhere vanishing.

It remains to find a holomorphic function $\tilde a$ and a nowhere vanishing holomorphic function $\tilde{\alpha}$ on $\tilde U$, such that for $\tilde{x}_j\in \tilde{D}_j, \,\tilde{x}_0\in \tilde{D}_0$ such that $p(\tilde{x}_j)=p(\tilde{x}_0)$, the equations
\begin{align}\label{eqEl50'}
\;\;\;\;\;\;\;\;\tilde{\alpha}(\tilde{x}_j) & =\alpha_j(x)\,\cdot \, \tilde{\alpha}(\tilde{x}_0),\,\;\; j =1,\ldots\,
\end{align}
\begin{align}\label{eqEl50}
\tilde{a}(\tilde{x}_j) & =a_j(x)+\alpha_j(x)\, \tilde{a}(\tilde{x}_0),
\end{align}
\noindent hold. Indeed, let $w$ be the biholomorphic mapping from $\tilde{U}\times \mathbb{P}^1$ that acts on the fiber over $\tilde{x}$ by the mapping $(\tilde{x},\zeta)\to (\tilde{x},\tilde {a}(\tilde{x}) + \tilde{\alpha}(\tilde{x}) \cdot \zeta)$.
Then the bundle $(\tilde{U}\times \mathbb{P}^1,\,{\rm pr}_1,\, w(\widetilde{\mathring{\mathbold{E'}}})\cup (\mathbold{s}^ {\infty}),\, \tilde{U})$ descends to a special $(0,4)$-bundle that is holomorphically isomorphic to $\mathfrak{F}$.

The equation \eqref{eqEl50'} leads to
a second Cousin problem for $\tilde{\alpha_j}$.
Equation \eqref{eqEl50} can be rewritten as
\begin{equation}\label{eqEl50''}
\frac{\tilde{a}(\tilde{x}_j)}{\tilde{\alpha}(\tilde{x}_j)}=\frac{a_j(x)}{\tilde{\alpha}(\tilde{x}_j)}
+\frac{\tilde{a}(\tilde{x}_0)}{\tilde{\alpha}(\tilde{x}_0)}\,.
\end{equation}
For convenience of the reader we provide the reduction of the equations \eqref{eqEl50'}
to the second Cousin problem.
Let $D$, $V_j$ and $\gamma_j$ be the same objects as in Section \ref{sec:9.2}.
Put $U_0=D$. We consider an open cover of $X$ as follows.
$U_0$ is an open topological disc in $X$ that contains the base point $x_0$. Cover each $V_j$ by two simply connected open sets $U_{j^+}$ and $U_{j^-}$
with connected and simply connected intersection, so that the  $U_{j^{\pm}}$ are disjoint from the $U_{k^{\pm}}$ for $j\neq k$, and
for each $U_{j^{\pm}}$ its intersection with $U_0$ is connected and simply connected. We may
also assume that the intersection of at least three of the sets is empty, the intersections of each $\gamma_j$ with $U_{j^+}$ and $U_{j^-}$ are connected, and each $\gamma_j$ is disjoint from $U_{k^+} \cup U_{k^-}$  for $k\neq j$.
Label the  $U_{j^{\pm}}$ so that walking on $\gamma_j\setminus U_0$ (which is contained in $V_j$) in the direction of the orientation of $\gamma_j$
we first meet $U_{j^-}$.

\begin{figure}[h]
\begin{center}
\includegraphics[width=65mm]{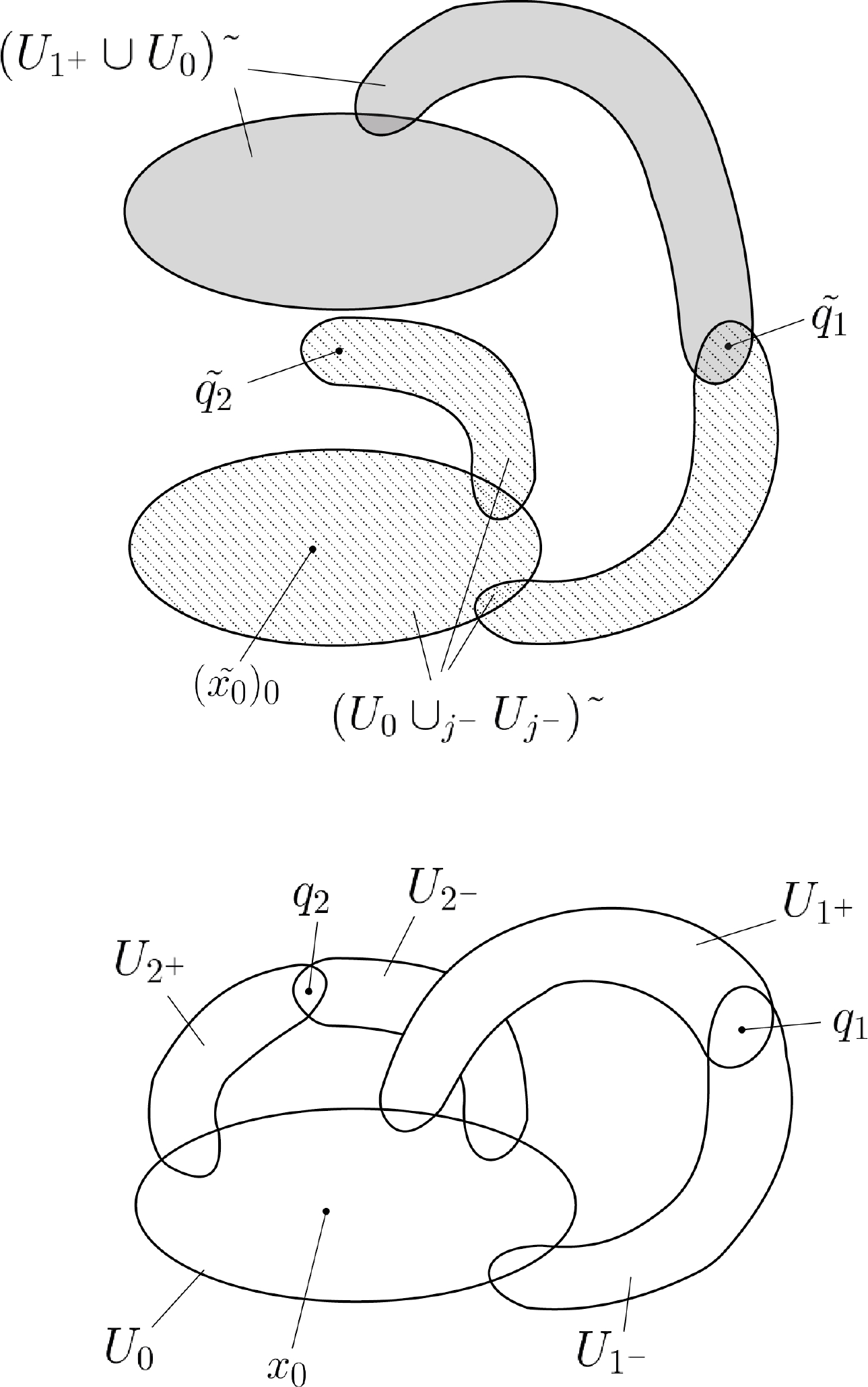}
\end{center}
\caption{The open cover of $X$ and its lift to $\tilde X$}\label{fig92}
\end{figure}

We consider the following subsets of $\tilde X$. Let
$(U_0\cup_{j^-}\, U_{j^-}){\;\widetilde{}}\;$  be a lift of $U_0\cup_{j^-}\, U_{j^-}$ to
$\tilde X$. Let $(\tilde{x_0})_0$ be the point in the lifted set that projects to $x_0$. For
each $j\geq 1$ we consider a point ${q}_j\in U_{j^-} \cap U_{j^+}$ and the point
$\tilde{q}_j$ in 
$(U_0\bigcup \cup_{j^-} U_{j^-}){\;\widetilde{}}\;$ that projects to $q_j$.
For each $j$ we denote by
$(U_{j^+}\cup U_0){\;\widetilde{}}\;$  the lift of  $U_{j^+}\cup U_0$ to $\tilde X$ that
contains $\tilde{q}_j$. Let $\widetilde{U_0^{j^+}}$ be the subset of
$(U_{j^+}\cup U_0){\;\widetilde{}}\;$ that projects to $U_0$.(See Figure \ref{fig92}.)

The sets $U_0$, $U_{j^{\pm}}$ cover $X$.
The intersection of three distinct sets of the cover is empty. On the intersection of
pairs of covering sets we define a non-vanishing holomorphic function as follows.
Put
\begin{align}\label{eqEl26}
\alpha_{0,j^-}\;\;=(\,\alpha_{j^-,0}\;)^{-1}= 1 \quad \quad \quad \;&\mbox{on}\; U_0 \;\;
\cap U_{j^-} \quad  \,j =1,\ldots, 2g+m\,,\nonumber\\
\alpha_{j^-,j^+}=(\alpha_{j^+,j^-})^{-1}= 1 \quad \quad \quad \;&\mbox{on}\; U_{j^-} \cap
U_{j^+} \quad  \,j =1,\ldots, 2g+m\,,\nonumber\\
\alpha_{j^+,0}\;\; = (\,\alpha_{0,j^+}\;)^{-1}= \alpha_j \quad \;\;
\quad \,&\mbox{on}\;U_{j^+}\cap  U_0
\;\; \quad \,j =1,\ldots, 2g+m\,.
\end{align}
Here $\alpha_j(x), \; z \in U_{j^+} \cap U_0,$ was defined by equation \eqref{eqEl50}.

Since there are no triple intersections, the equations \eqref{eqEl26} define a Cousin II distribution.
A second Cousin distribution on a complex manifold defines a holomorphic line bundle. The second Cousin problem has a solution if and only if the line bundle is holomorphically trivial. Since on an open Riemann surface each holomorphic line bundle is holomorphically trivial the Cousin problem has a solution (see e.g. \cite{Fo}).

The solution of the Cousin problem for the Cousin distribution
\eqref{eqEl26}
provides non-vanishing holomorphic functions $\alpha_0$ on $U_0$ and
$\alpha_{j^\pm}$ on $U_{j^\pm}$ such that
\begin{align}\label{eqEl27}
\alpha_0 \;\;\;\;=\alpha_{j^-} \quad\quad \;\;\;\; &\mbox{on}\; U_0\;\,\cap
U_{j^-},\nonumber\\
\alpha_{j^-}\;\;= \alpha_{j^+} \quad\quad \;\;\;\; &\mbox{on}\; U_{j^-}\cap
U_{j^+},\,\nonumber\\
\alpha_{j^+,0}= \alpha_{j^+} \alpha_{0}^{-1}\;\;\;\;\;\, &\mbox{on}\; U_{j^+}\cap U_0\;\;.
\end{align}

Put
\begin{equation}
\tilde{\alpha}(\tilde{x})=\begin{cases}
\alpha_0({x})& \tilde{x}\in\tilde{U}_0\\
\alpha_{j^-}({x})& \tilde{x}\in\tilde{U}_{j^-}\\
\alpha_{j^+}({x})& \tilde{x}\in\tilde{U}_{j^+}\\
\alpha_j(x)\cdot \alpha_0({x})& \tilde{x}\in \widetilde{U_0^{j^+}}\,.
\end{cases}
\end{equation}
Here $x=p(\tilde{x})$ for the projection $p:\tilde{X}\to X$.
By \eqref{eqEl50} the function $\tilde{\alpha}$ is a well-defined holomorphic function on $\tilde U$ with the required property.

The reduction of \eqref{eqEl50} to a first Cousin problem
is similar. The first Cousin problem is solvable since open Riemann surfaces are Stein manifolds (see e.g. \cite{Fo} and \cite{H1}).
\hfill $\Box$

\medskip

We will consider now deformations of a smooth $(0,3)$-bundle with a section over a finite Riemann surface $X$  to a holomorphic bundle.
Statement $(1)$ of Theorem \ref{thmEl.0} below  is the analog of Theorem \ref{thm8.0} for $(0,3)$-bundles with a section.
\begin{defn}\label{defnEl.3'}
We will say that a smooth $({\sf g,m})$-bundle $\mathfrak{F}$ over a smooth finite open surface $X$ is isotopic to a holomorphic bundle for the conformal structure $\omega:X\to \omega(X)$, if the pushed forward bundle $\mathfrak{F}_{\omega}$ (see equation \eqref{eq23}) is isotopic to a holomorphic bundle on $\omega(X)$. If the bundle is isotopic to a holomorphic bundle for each conformal structure with only thick ends on $X$, then $\mathfrak{F}$ is said to have the Gromov-Oka property.
\end{defn}
A $(0,n)$-bundle with a section, in particular, a special $(0,n+1)$-bundle, has the Gromov-Oka property, if and only if for each conformal structure on $X$ with only thick ends the bundle is isotopic through
$(0,n)$-bundles with a section to a holomorphic bundle.

In the following theorem we start with a smooth special $(0,4)$-bundle, since each smooth $(0,3)$-bundle with a section is isotopic to a special $(0,4)$-bundle.
Recall that each finite open Riemann surface $\mathcal{X}$ is conformally equivalent to a domain in a closed Riemann surface $\mathcal{X}^c$.

\begin{thm}\label{thmEl.0}
$(1)$ Let $X$ be a smooth surface of genus one with a hole with base point $x_0$, and with a chosen set
$\mathcal{E}=\{e_1,e_2\}$ of generators of $\pi_1(X,x_0)$. Define the set
$\mathcal{E}_0=
\{e_1,e_2, e_1 e_2^{-1},  e_1 e_2^{-2}, e_1e_2 e_1^{-1} e_2^{-1}\}$ as in Theorem {\rm \ref{thm8.0}}. Consider a smooth special $(0,4)$-bundle $\mathfrak{F}$ over $X$. Suppose for each $e\in \mathcal{E}_0$ the restriction of the bundle $\mathfrak{F}$ to an annulus
in $X$ representing $\hat e$ has the Gromov-Oka property,
Then the bundle over $X$ has the Gromov-Oka property.\\
$(2)$ If a bundle $\mathfrak{F}$ as in $(1)$ is irreducible, then it is isotopic to an isotrivial bundle, and, hence,
for any conformal structure $\omega$ on $X$
the bundle $\mathfrak{F}_{\omega}$ is isotopic to a bundle that extends to a holomorphic bundle on $\omega(X)^c$. In particular, the bundle is isotopic to a holomorphic bundle for any conformal structure on $X$ (including conformal structures of first kind).\\
$(3)$ If $\mathfrak{F}$ is any smooth reducible $(0,4)$-bundle over $X$, then each irreducible bundle component is isotopic to an
isotrivial bundle. There exists a Dehn twist in the fiber over the base point 
such that the bundle $\mathfrak{F}$ can be recovered from the irreducible bundle components up to composing each monodromy with a power of this Dehn twist. $\mathfrak{F}$ is isotopic to a holomorphic bundle for each conformal structure of second kind on $X$.\\
$(4)$ Any reducible holomorphic $(0,3)$-bundle with a holomorphic section over a punctured Riemann surface is holomorphically trivial.
\end{thm}

A bundle is called isotrivial if its lift to a finite covering of $X$ is the trivial bundle.
We will prove now Theorems \ref{thm8.0} and \ref{thmEl.0} using Theorem \ref{thm8.1}.
\medskip

\noindent {\bf Proof of Theorem \ref{thmEl.0}.}
The bundle $\mathfrak{F}$ defines a smooth quasipolynomial $f$ on $X$. For each $e\in \mathcal{E}_0$ the restriction of the bundle to an annulus $A_{\hat{e}}$ in $X$ representing $\hat e$ has the Gromov-Oka property. By Lemma \ref{lem*}
we may assume that for a conformal structure of conformal module bigger than $\frac{\pi}{2}\log(\frac{3+\sqrt{5}}{2})^{-1}$ on $A_{\hat{e}}$ the restriction $\mathfrak{F}\mid A_{\hat{e}}$ is isotopic (through $(0,3)$-bundles with a section) to a special holomorphic $(0,4)$-bundle $\mathfrak{F}_{\hat{e}}$ (not merely to a holomorphic $(0,3)$-bundle with a section). The bundle $\mathfrak{F}_{\hat{e}}$ defines a holomorphic quasipolynomial $f_{\hat{e}}$ on $A_{\hat{e}}$. The monodromies of the isotopic (hence, isomorphic) bundles  $\mathfrak{F}\mid A_{\hat{e}}$ and $\mathfrak{F}_{\hat{e}}$ differ by conjugation (after identification of the mapping class groups on the fiber over the base point by an isomorphism).
Hence, for each $e\in\mathcal{E}$ there exists an integer number $k_e$,
such that the quasipolynomial $f \mid A_{\hat{e}}$ is (free) isotopic to the holomorphic quasipolynomial
$z \to e^{2\pi ik_e z} \,f_{\hat{e}}(z),\, z\in  A_{\hat{e}}$. We proved that the quasipolynomial $f$ satisfies the conditions of Theorem \ref{thm8.1}.

The monodromy mapping classes of the quasipolynomial $f$ are elements of the braid group $\mathcal{B}_3$, the monodromy mapping classes of the bundle
$\mathfrak{F}$ can be identified with elements of $\mathcal{B}_3\diagup \mathcal{Z}_3$.
In terms of the $(0,4)$-bundle $\mathfrak{F}$ the Theorem \ref{thm8.1}
implies the following. The isotopy class of the special $(0,4)$-bundle
$\mathfrak{F}$ on $X$ corresponds to the conjugacy class of a homomorphism $\mathbold{\Phi}:\pi_1(X,x_0)\to
\mathbold{\Gamma} \subset \mathcal{B}_3\diagup \mathcal{Z}_3$, where
$\mathbold{\Gamma}$
is generated either by $\sigma_1 \,\sigma_2\diagup \mathcal{Z}_3$, or by
$\Delta_3\diagup\mathcal{Z}_3$, or by $\sigma_1\diagup\mathcal{Z}_3$.

The group $\mathbold{\Gamma}$ is Abelian. Let $\omega$ be any conformal structure on $X$ (maybe, of first kind). The
fundamental group $\pi_1 (\, \omega(X)^{c}, \omega(x_0))$ of the closed torus
$\,\omega(X)^{c}$ containing (a conformal copy of) $\omega(X)$ is the Abelianization of $\pi_1 (\omega(X),\omega(x_0))\cong \pi_1 ( \, X , x_0)$. Hence
$\mathbold{\Phi}$ defines a homomorphism from $\pi_1 (\, \omega(X)^{c}, \omega(x_0))$ to
$\mathbold{\Gamma}$ which we also denote by $\mathbold{\Phi}$.

Consider first the case when the generator of $\mathbold{\Gamma}$ is either
$\sigma_1\sigma_2 \diagup \mathcal{Z}_3$, or
$\Delta_3\diagup \mathcal{Z}_3$. In these cases
the bundle $\mathfrak{F}$ is irreducible
by Lemma \ref{lemEl.0}.
The generator of $\mathbold{\Gamma}$ corresponds to an element of a
mapping class group  $\mathfrak{M}(\mathbb{P}^1; \{\infty\},\mathring{ E}_3)$,
that is represented by a periodic mapping $\zeta\to \theta \cdot \zeta,\,
\zeta \in \mathbb{P}^1$ (see \cite{Be1}). In the first case $\theta = e^{\frac{2\pi i}3}$ and $\mathring{E}_3$ consists of
$3$ equidistributed points on a circle with center zero,
in the second case  $\theta= e^{\frac{2\pi i}{2}}$ and $\mathring{E}_3$ consists of the origin and two
equidistributed points on a circle with center zero.

Consider the homomorphism from $\mathbold{\Gamma}$ to the group of complex
linear transformations of ${\mathbb P}^1$, which assigns to the
generator $\mathbold{b} \in \mathbold{\Gamma}$ the complex linear transformation $\zeta \to \theta
\cdot \zeta$, $\zeta \in {\mathbb P}^1$.
Composing $\mathbold{\Phi}$ with
this homomorphism we obtain a homomorphism $\Xi$ from $\pi_1 (\, \omega(X)^{c}, \omega(x_0))$
to the group of complex linear transformations of ${\mathbb
P}^1$.

Represent the closed torus $\, \omega(X)^{c}$ as quotient $ \,\omega(X)^{c}  = {\mathbb C} \diagup \Lambda$ for a lattice $\Lambda
=\{\lambda_1 \, k_1 + \lambda_2 \, k_2 : k_1 , k_2 \in {\mathbb
Z}\}$ with real linearly independent numbers $\lambda_1 , \lambda_2
\in {\mathbb C}$. Denote by $p_1$ the covering map $p_1 : {\mathbb
C} \to {\mathbb C} \diagup \Lambda$. Identify the fundamental group
$\pi_1 ({\mathbb C} \diagup \Lambda , z_0) = \pi_1 ( \, \omega(X)^{c},
\omega(z_0))$ with the group of covering transformations which we also denote
by $\Lambda$.

The fundamental group $\Lambda$ of $ \, \omega(X)^{c} \cong {\mathbb C}
\diagup \Lambda$ acts on ${\mathbb C} \times {\mathbb P}^1$ as
follows. Associate to $\lambda \in \Lambda$ the mapping
\begin{equation}
\label{eq8.1} (z,\zeta) \overset{\lambda}{\longrightarrow} (\lambda
(z) , \Xi (\lambda) (\zeta)) \, , \quad (z,\zeta) \in {\mathbb C}
\times {\mathbb P}^1 \, .
\end{equation}
The thus defined mapping $\lambda$ takes the set $\mathbb{C}\times\mathring{E}_3$ to itself and fixes the set $\mathbb{C}\times\{\infty\}$.
The action is free and properly discontinuous. Hence the quotient
$({\mathbb C} \times {\mathbb P}^1) \diagup \Lambda$ is a complex
manifold. For each element of the quotient its projection to
${\mathbb C} \diagup \Lambda$ is well-defined:
$$
\xymatrix{
(z,\zeta) \ar[d] &\sim&(\lambda (z) , \Xi (\lambda) (\zeta)) \ar[d] \\
z &\sim&\lambda(z) }
$$
The equivalence relation in the upper line gives the quotient
$({\mathbb C} \times {\mathbb P}^1) \diagup \Lambda$, the relation
in the lower line gives ${\mathbb C} \diagup \Lambda$.

Denote by $ \mathcal{P}$
the projection $\mathcal{P}:({\mathbb C} \times {\mathbb P}^1)\diagup \Lambda
\to {\mathbb C} \diagup \Lambda$.
We obtain a holomorphic $(0,3)$-bundle $\mathfrak{F}_{\omega}\stackrel{def}=\Big(({\mathbb C} \times {\mathbb P}^1)\diagup \Lambda,\,\mathcal{P},\, (\mathbb{C}\times\mathring{E}_3)\diagup\Lambda \cup(\mathbb{C}\times\{\infty\})\diagup \Lambda,{\mathbb C} \diagup \Lambda\Big)$ with a holomorphic section over the closed torus $\omega(X)^c$.
By construction the restriction of this bundle to a representative of each generator
$e_j$ of $\pi_1(\omega(X)^{c},x_0)$ gives the mapping torus corresponding to $\Xi(e_j)$. Hence, the
monodromy mapping classes of the bundle are equal to $\mathbold{\Phi}(e_j)$. We proved that
in the irreducible case of Theorem \ref{thmEl.0} for a conformal structure $\omega$ on $X$ the bundle $\mathfrak{F}_{\omega}$ is smoothly isomorphic (equivalently, isotopic) to a holomorphic
bundle that extends to a holomorphic bundle on the closed torus $\omega(X)^c$.

Since the monodromy mapping classes are periodic, the lift of the bundle $(\mathfrak{F}^c)_{ \omega}$ to a finite covering of $X$ is isotopic to the trivial bundle.
Statement (2) is proved and also statement (1) in the irreducible case.

Suppose $\mathfrak{F}$ is an arbitrary reducible smooth special $(0,4)$-bundle.
There is an admissible curve $\gamma\subset \mathbb{P}^1\setminus (\mathring{E}_3\cup \{\infty\})$ that reduces each monodromy mapping class of the bundle. Hence, by (the proof of) Lemma \ref{lemEl.0} each monodromy mapping class of the bundle is a non-trivial
power of (a single conjugate of) the mapping class $\mathfrak{m}({\sigma_1\diagup \mathcal{Z}_3})$ in
$\mathfrak{M}(\mathbb{P}^1;\infty,\mathring{E}_3)$ that corresponds to  ${\sigma_1\diagup \mathcal{Z}_3}$. We are in the case when $\mathbold{\Gamma}$ is generated by $\sigma_1\diagup \mathcal{Z}_3$.
Label the points of $\mathring{E}_3$ as in Lemma \ref{lemEl.0}.
The simple closed curve $\gamma$ separates $\{\zeta_1,\zeta_2\}$ from $\{\zeta_3,\infty\}$.
Let $\mathcal{C}_1$ be the connected component of the complement of $\gamma$ that contains $\{\zeta_1,\zeta_2\}$, and let $\mathcal{C}_2$ be the connected component of the complement of $\gamma$ that contains $\{\zeta_3,\infty\}$.

There are two irreducible components of the
mapping class $\mathfrak{m}({\sigma_1\diagup \mathcal{Z}_3})$. The irreducible component  that corresponds to $\mathcal{C}_1$ is described as follows.
The connected component $\mathcal{C}_1$ is
homeomorphic to $\mathbb{P}^1\setminus \{\infty\}$ with set of distinguished points $\{\zeta_1,\zeta_2\}$. Each self-homeomorphism of $\mathbb{P}^1\setminus \{\infty\}$ with these two distinguished points extends to a self-homeomorphism of $\mathbb{P}^1$ that maps the set $\{\zeta_1,\zeta_2, \infty\}$ of three distinguished points to itself.
The irreducible component of the mapping class $\mathfrak{m}({\sigma_1\diagup
\mathcal{Z}_3})$ corresponding to $\mathcal{C}_1$ is an element of
the mapping class group $\mathfrak{M}(\mathbb{P}^1; \{\infty\}, \{\zeta_1,\zeta_2\})$.
This mapping class group is generated by the element $\mathbold{\sigma}$ which
corresponds to the braid $\sigma$ in the braid group $\mathcal{B}_2$ on two strands. The
square of $\mathbold{\sigma}$ is the identity. We put  $\mathbold{\Gamma}_1=
\mathfrak{M}(\mathbb{P}^1; \{\infty\}, \{\zeta_1,\zeta_2\})$.

To describe the irreducible component corresponding to $\mathcal{C}_2$ of the mapping class $\mathfrak{m}(\sigma_1\diagup \mathcal{Z}_3)$, we notice that
$\mathcal{C}_2$ is homeomorphic to $\mathbb{P}^1\setminus \{\zeta_4\}$ with distinguished points $\{\zeta_3,\infty\}$ for a point $\zeta_4\in \mathbb{P}^1, \, \zeta_4\neq \zeta_3,\infty$.
Each self-homeomorphism of $\mathcal{P}^1\setminus \{\zeta_4\}$ with distinguished points
$\{\zeta_3,\infty\}$ extends to a self-homeomorphism of
$\mathbb{P}^1$ with set of distinguished points $\{\zeta_3,\zeta_4,\infty\}$. The respective irreducible component of $\mathfrak{m}({\sigma_1\diagup \mathcal{Z}_3})$ is the identity in $\mathfrak{M}(\mathbb{P}^1; \{\zeta_4, \infty\},z_3)$.

All monodromy mapping classes of a bundle that is isomorphic to $\mathfrak{F}$ are contained in $\mathbold{\Gamma}$, hence are powers of $\mathfrak{m}({\sigma_1\diagup \mathcal{Z}_3})$.
This implies that there are two irreducible bundle components, corresponding to $\mathcal{C}_1$ and to $\mathcal{C}_2$, respectively.
They are described as follows.

The irreducible bundle component of the bundle $\mathfrak{F}$, corresponding to $\mathcal{C}_2$, is associated to the trivial homomorphism. Hence, this irreducible bundle component is smoothly isomorphic to the trivial $(0,3)$-bundle over $X$.

The irreducible bundle component of $\,\mathfrak{F}\,$ related to $\mathcal{C}_1$ corresponds to the conjugacy
class of the homomorphism $\;\mathbold{\Phi}_1: \pi_1(X,x_0)\to \mathbold{\Gamma}_1$
induced by the homomorphism $\mathbold{\Phi}$. The same proof as in the first option for
$\mathbold{\Gamma}$ shows that for any conformal structure $\omega$ on $X$ there is a holomorphic
$(0,2)$-bundle $(\mathfrak{F}'_1)^{\omega}$
with a holomorphic section
on the closed torus $\omega(X)^{c}$ whose isomorphism class corresponds to
the conjugacy class of the homomorphism $\pi_1(\omega(X)^c,\omega(x_0)) \to
\mathbold{\Gamma}_1$ that is induced by $\mathbold{\Phi}_1$.
The restriction of $(\mathfrak{F}'_1)^{\omega}$ to a punctured torus $\omega({X})'$, $\omega(X)\Subset \omega({X})'$,
is isomorphic to a special holomorphic $(0,3)$-bundle $(\mathfrak{F}_1)^{\omega}$ with set of distinguished points
$\{\zeta_1(x),\zeta_2(x),\infty\}$ in the fiber $\{x\}\times \mathbb{P}^1$ over $x$. For
$x\in \omega(X) \Subset \omega(X)'$ the points $\zeta_1(x),\zeta_2(x)$ are contained in a large closed disc in
$\mathbb{C}$ centered at the origin.

Take a point $\zeta_3\in \mathbb{C}$ outside this closed disc.
Consider the special holomorphic $(0,4)$-bundle $(\mathfrak{F})^{\omega}
\stackrel{def}=\Big(X \times \mathbb{P}^1, {\rm pr}_1, \mathring{\mathbold{E}}_3\cup\mathbold{s}^{\infty}, X\Big)$,
where the set of finite distinguished points $\mathring{\mathbold{E}}_3\cap(\{x\}\times \mathbb{P}^1)$ in the fiber over $x$ equals
$\{x\}\times\{\zeta_1(x) ,\zeta_2(x), \zeta_3\}$ and $\mathbold{s}^{\infty}\cap(\{x\}\times \mathbb{P}^1)=(x,\infty)$.
The monodromy mapping class of the bundle $(\mathfrak{F})^{ \omega}$ along $e_j$ differs from the monodromy mapping class
of the push-forward $\mathfrak{F}_{ \omega}$ to $\omega(X)$ of the bundle $\mathfrak{F}$ by the $k_j$-th power of a Dehn twist about the curve $\gamma$.
Consider again the punctured torus $\omega(X)'\supset \omega(X)$, and assume
that $\omega(X)'=(\mathbb{C}\setminus \Lambda) \diagup \Lambda$
with
$\Lambda = \{ \lambda_1 \, n_1 + \lambda_2 \, n_2 \, , \ n_1 ,
n_2 \in {\mathbb Z}\}$, and $x_0=\frac{\lambda_1+\lambda_2}{2}\diagup \Lambda$, and that the segments $\frac{\lambda_1+\lambda_2}{2}+[0,\lambda_1]$ and $\frac{\lambda_1+\lambda_2}{2}+[0,\lambda_2]$ are lifts to $\mathbb{C}$ of curves representing a pair of generators of $\pi_1 ({\mathbb C} \diagup \Lambda , x_0)$.

For all integers $\ell_1$ and $\ell_2$ there is a holomorphic function $F$ on ${\mathbb C}
\backslash
\Lambda$ such that
\begin{equation}
\label{eq8.18} F(x+\lambda_1) = \ell_1 \quad {\rm and} \quad F(x +
\lambda_2) = \ell_2, \; x \in \mathbb{C}\setminus \Lambda\,.
\end{equation}
This follows from the proof of Lemma \ref{lem*}.
An elementary way to see this for $\ell_2 = 0$ and any $\ell_1 \in
{\mathbb Z}$ is the following. For each $x\in \mathbb{C}$ there are uniquely determined
real numbers $t_1(x)$ and $t_2(x)$ depending smoothly on $x$ such that $x=\lambda_1 t_1(x)+\lambda_2 t_2(x)$.
Consider a function $\chi_0$ on the real axis which vanishes near zero such that
$\chi_0(t+1)=\chi_0(t)+\ell_1$, $t \in \mathbb{R}$. Consider the function $\chi(x)\stackrel{def}=\chi_0(t_1(x)),\, x\in \mathbb{C}$.

The
1-form $\overline\partial \, \chi$ descends to a smooth 1-form
$\delta$ on the torus ${\mathbb C} \diagup \Lambda$. Let $\mathring{g}$ be a
solution of the equation $\overline\partial \, \mathring{g} = \delta$ on the
punctured torus ${\mathbb C} \backslash \Lambda \diagup \Lambda$.
Let $g$ be the lift of $\mathring{g}$ to ${\mathbb C}
\backslash \Lambda$. Put $F = \chi - {g}$.

Let $F$ be the holomorphic function on
${\mathbb C} \backslash \Lambda$ that satisfies equation
\eqref{eq8.18}
for the integers $\ell_1=-k_1$ and $\ell_2=-k_2$.
Then the monodromy mapping classes along the $e_j$ of the special holomorphic $(0,4)$-bundle $(\mathfrak{F}')^{\omega}$ over $\omega(X)$
with set of distinguished points
$\{\zeta_1(x) ,\zeta_2(x), \zeta_3 \,e^{2\pi i F(x)},\infty\}$ in each fiber $\{x\}\times \mathbb{P}^1$, $x\in \omega(X)$,
are the same as the monodromy mapping classes of the bundle $\mathfrak{F}_{ \omega}$. Statement (3) and, hence, statement (1) in the reducible case are proved.

It remains to prove that any reducible holomorphic $(0,3)$-bundle $\mathfrak{F}$
with a holomorphic section over a punctured Riemann surface $X$ is holomorphically isomorphic to the trivial bundle. The respective statement is known in more general situations,
but the proof of the present statement is a simple reduction to Picard's Theorem and we include it.

By Lemma \ref{lem*} we may assume that the bundle is equal to a special $(0,4)$-bundle
$\mathfrak{F}=(X\times \mathbb{P}^1,{\rm pr}_1, \mathring{\mathbold{E}}\cup \mathbold{s}^{\infty},X)$.
There is a finite unramified covering $\hat X$ of the punctured Riemann surface $X$
so that for the lift $\hat {\mathfrak{F}}= \big(\hat{X}\times \mathbb{P}^1,{\rm pr}_1,\widehat{\mathring{\mathbold{E}}}\cup \hat{\mathbold{s}}^{\infty},\hat{X}\big)$ of the bundle $\mathfrak{F}$ to $\hat X$ the set
$\widehat{\mathring{\mathbold{E}}}\cup \hat{\mathbold{s}}^{\infty}$ is the union of four disjoint complex curves each intersecting each fiber along a single point. After a holomorphic isomorphism (see the proof of Lemma \ref{lem*}) we may assume that the connected components of $\widehat{\mathring{\mathbold{E}}}$ are $\hat{X}\times\{-1\}$, $\hat{X}\times\{1\}$, and a component that intersects each fiber $\hat{x}\times \mathbb{P}^1$ along a point $(\hat{x},g(\hat{x}))$ for a holomorphic function $g$ on $\hat X$ that omits the values $-1,\,1$, and $\infty$.

Since the bundle ${\mathfrak{F}}$ is reducible, the monodromy mapping classes
of the bundle $\hat{\mathfrak{F}}$
are powers of a conjugate of $\sigma_1^2\diagup \mathcal{Z}_3$.
Hence, $g:\hat{X}\to \mathbb{C}\setminus \{-1,1\}$ is free homotopic to a mapping with image in the punctured disc $\{z\in\mathbb{C}:|z+1|<1\}$.

By Picard's Theorem the mapping $g$ extends to a meromorphic mapping $g^c$ from the closure $\hat{X}^c$ of $\hat X$ to $\mathbb{P}^1$.
Then the meromorphic extension $g^c$ omits the value $1$. Indeed, if $g^c$ was equal to $1$ at some puncture of $\hat X$, then $g$ would map a simple closed curve on $\hat X$ that surrounds the puncture positively to a loop in $\mathbb{C}\setminus \{-1,1\}$ that surrounds $1$ with positive winding number. This
contradicts the fact that $g$ is homotopic to a mapping onto a disc punctured at $-1$ and contained in $\mathbb{C}\setminus \{-1,1\}$. Hence, $g^c$ is a meromorphic function on a compact Riemann surface that omits a value, and, hence $g$ is constant.

This means that the bundle $\hat{\mathfrak{F}}$ over $\hat X$ is holomorphically isomorphic to the trivial bundle. Hence, the monodromy mappings of the original bundle $\mathfrak{F}$ are periodic. By Lemma \ref{lemEl.0}
this is possible for a reducible bundle only if the monodromies are equal to the identity. Repeat the same reasoning as above for the bundle $\mathfrak{F}$ instead of $\hat{\mathfrak{F}}$,
we see that the bundle $\mathfrak{F}$ is holomorphically trivial.
Theorem \ref{thmEl.0} is proved. \hfill $\Box$

\medskip

\noindent {\bf Proof of Theorem \ref{thm8.0}.}
By Theorem \ref{thm8.1} the isotopy class of the separable quasipolynomial $f$ corresponds to the
conjugacy class of a homomorphism $\Phi: \pi_1(X,x_0) \to \Gamma \subset \mathcal{B}_3$,
where $\Gamma$ is generated either by $\sigma_1\sigma_2$, or by
$\Delta_3$, or by $\sigma_1$ and $\Delta_3^2$. Consider the special $(0,4)$-bundle
associated to the quasipolynomial.
The isotopy class of the bundle over $X$
(through $(0,3)$-bundles with a section)
corresponds to the conjugacy class of the homomorphism $\mathbold{\Phi}: \pi_1(X,x_0) \to
\mathbold{\Gamma}=\Gamma\diagup\mathcal{Z}_3 $ that is obtained from $\Phi$ by taking the
quotient with respect to $\mathcal{Z}_3$.

By Theorem \ref{thmEl.0} for any a priori given conformal structure $\omega$ of second kind on $X$ the isotopy class of the bundle contains a special holomorphic $(0,4)$-bundle.
The special holomorphic
$(0,4)$-bundle determines a holomorphic quasipolynomial $f^{\omega}$ for this conformal structure by
assigning to each $x\in \omega(X)$ the triple of finite distinguished points in the fiber
$\{x\}\times \mathbb{P}^1$ and associating to it the monic polynomial with this set of zeros.
The monodromies of the quasipolynomial $f^{\omega}$ differ from the monodromies of the push forward $f_{\omega}$ of the original quasipolynomial by powers of $\Delta_3^2$.
Using the function $F$ that satisfies
\eqref{eq8.18} for suitable integers $\ell_1$ and $\ell_2$ we obtain a holomorphic quasipolynomial $f^{\omega}(x) e^{2\pi i F(x)},\, x \in \omega(X),$    that is isotopic to $f_{\omega}$.
Theorem \ref{thm8.0} is proved. \hfill $\Box$

\section{Finite open Riemann surfaces and mappings to the twice punctured complex plane}
\label{sec:9.4a}
In this section we prove Theorem \ref{prop0}.
Let $X$ be an open Riemann surface of genus $g\geq 0$ with $m>0$ holes, equipped with a base point $x_0$.
The fundamental group $\pi_1(X,x_0)$ of $X$ is a free group in $x\stackrel{def}=2g+m-1$ generators.
We take the system $\mathcal{E}$ of generators of $\pi_1(X,x_0)$ as in Section \ref{sec:9.4a}.
Recall that the Riemann surface $X$ is conformally equivalent to an open subset of a closed Riemann surface $X^c$ of genus $g$ (\cite{Sto}).
in Section \ref{sec:9.4a} we chosed $2g$ generators $e_{j}, \, j=1,\ldots ,2g,$
of the fundamental group $\pi_1(X,x_0)$  whose images
under the homomorphism induced by inclusion $X \to X^c$
are generators of the fundamental group of $X^c$, such that
each pair $e_{2j-1}, e_{2j},\, j=1,\ldots g,$ corresponds to a handle.
Each
pair of generators $e_{2j-1},\,e_{2j}$ corresponding to a handle was represented by simple closed loops $\alpha_j,\,\beta_j$ with base point $x_0$.
The connected components of $X^c \setminus X$ are again denoted by $\mathcal{C}_1, \mathcal{C}_2,\ldots, \mathcal{C}_{m}$. Each component is either a point or a closed disc.
For each component $\mathcal{C}_\ell,\, \ell=1,\ldots,m-1,$ the generator $e_{2g+\ell}$ of the fundamental group of $X$ is represented by a simple loop $\gamma_{\ell}$ with base point $x_0$
that is contractible in $X\cup \mathcal{C}_\ell$ and divides $X$ into two connected
components, one of them containing $\mathcal{C}_\ell$. We will say that $\gamma_{\ell}$ surrounds $\mathcal{C}_\ell$. Each curve $\gamma_{\ell}$ is oriented so that $\mathcal{C}_{\ell}$ is on the left when walking along the curve $\gamma_{\ell}$ equipped with this orientation.

\begin{figure}[h]
\begin{center}
\includegraphics[width=80mm]{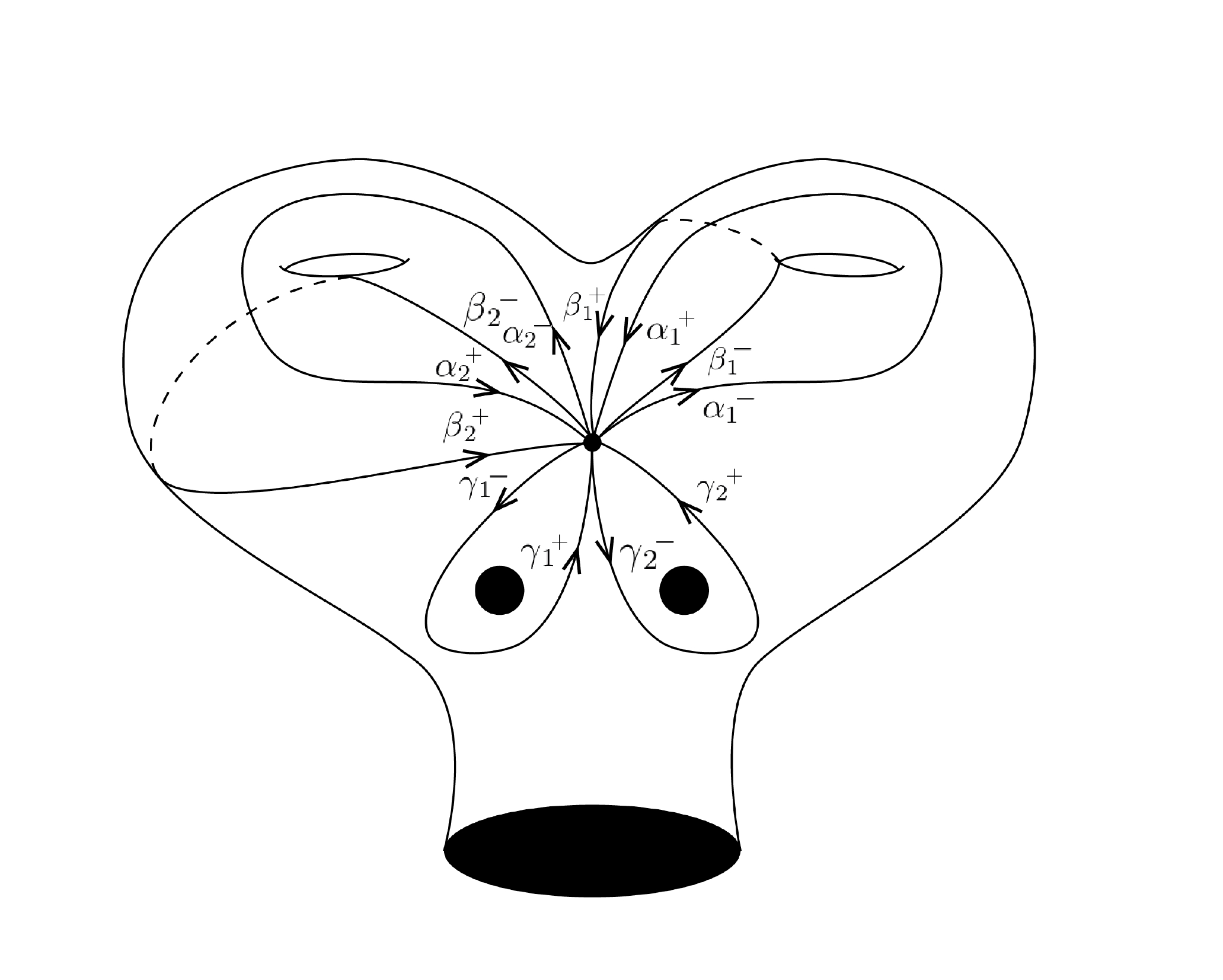}
\end{center}
\caption{A system of generators of the fundamental group of a surface of genus $g$ with $m$ holes}\label{fig8.3}
\end{figure}

Further, the generators and the representatives are chosen so that the only intersection point of any pair of loops is $x_0$, and when
labeling the rays of the loops emerging from the base point $x_0$ by  $\alpha_j^-,\,\beta_j^-$  $\gamma_j^-$, and the incoming rays by $\alpha_j^+,\,\beta_j^+$  $\gamma_j^+$,
then moving in counterclockwise direction along a small circle around $x_0$ we meet the rays in the order
\begin{align*}
\ldots, \alpha^-_j,\beta^-_j,\alpha^+_j,\beta^+_j,\ldots, \gamma^-_k,\gamma^+_k,\ldots\;.
\end{align*}
See also Figure \ref{fig8.3}.

We will now describe a set $\mathcal{E}'$ of elements of $\pi_1(X,x_0)$ such that the conjugacy class of
each element of  $\mathcal{E}'$ can be represented by a simple closed curve.
Later we will prove that
a mapping $X\to \mathbb{C}\setminus \{-1,1\}$ has the Gromov-Oka property if and only if for each $e\in \mathcal{E}'$ the restriction of $f$ to an annulus in $X$ that represents $\hat e$ has the Gromov-Oka property.

Let first $g>0$. For the $j$-th handle we choose three elements
$e_{2j-1}$, $e_{2j}$, and $[e_{2j-1},e_{2j}]$. The conjugacy class of each of them can be represented by a simple closed curve. For the $\ell$-th hole, $\ell=1,\ldots,m-1,$ we take the element $e_{2j+\ell}$. There is a simple closed curve that represents its conjugacy class. The obtained $3g+m-1$ elements of the fundamental group will be contained in $\mathcal{E}'$.

We convert each unordered pair of different elements of $\mathcal{E}$, that is different from a pair $\{e_{2j-1}$, $e_{2j}\}$ corresponding to a handle, into an ordered pair $(e',e'')$.
There is a simple closed curve representing $\widehat{e'\,e''}=\widehat{e''\,e'}$.
The products $e' e''$
constitute a set of
$\frac{1}{2}(2g+m-1)(2g+m-2) - g$
elements of the fundamental group. They will be contained in $\mathcal{E}'$.

For each pair $e_{2j-1},e_{2j}$ of elements of $\mathcal{E}$ that corresponds to a handle, and each other element $e'$ among the chosen generators (if there is any) we consider the elements $e_{2j-1}^2\,e_{2j}\, e'$, $e_{2j-1}^3\,e_{2j}\, e'$,
$e_{2j-1}\,e_{2j}^2\, e'$, and $e_{2j-1}\,e_{2j}^3\, e'$. For the conjugacy class of each such element there is a simple closed curve in $X$ representing it. The
obtained $4g\, (2g+m-3)$ elements
of the fundamental group will be contained in $\mathcal{E}'$.

\begin{figure}[h]
\begin{center}
\includegraphics[width=100mm]{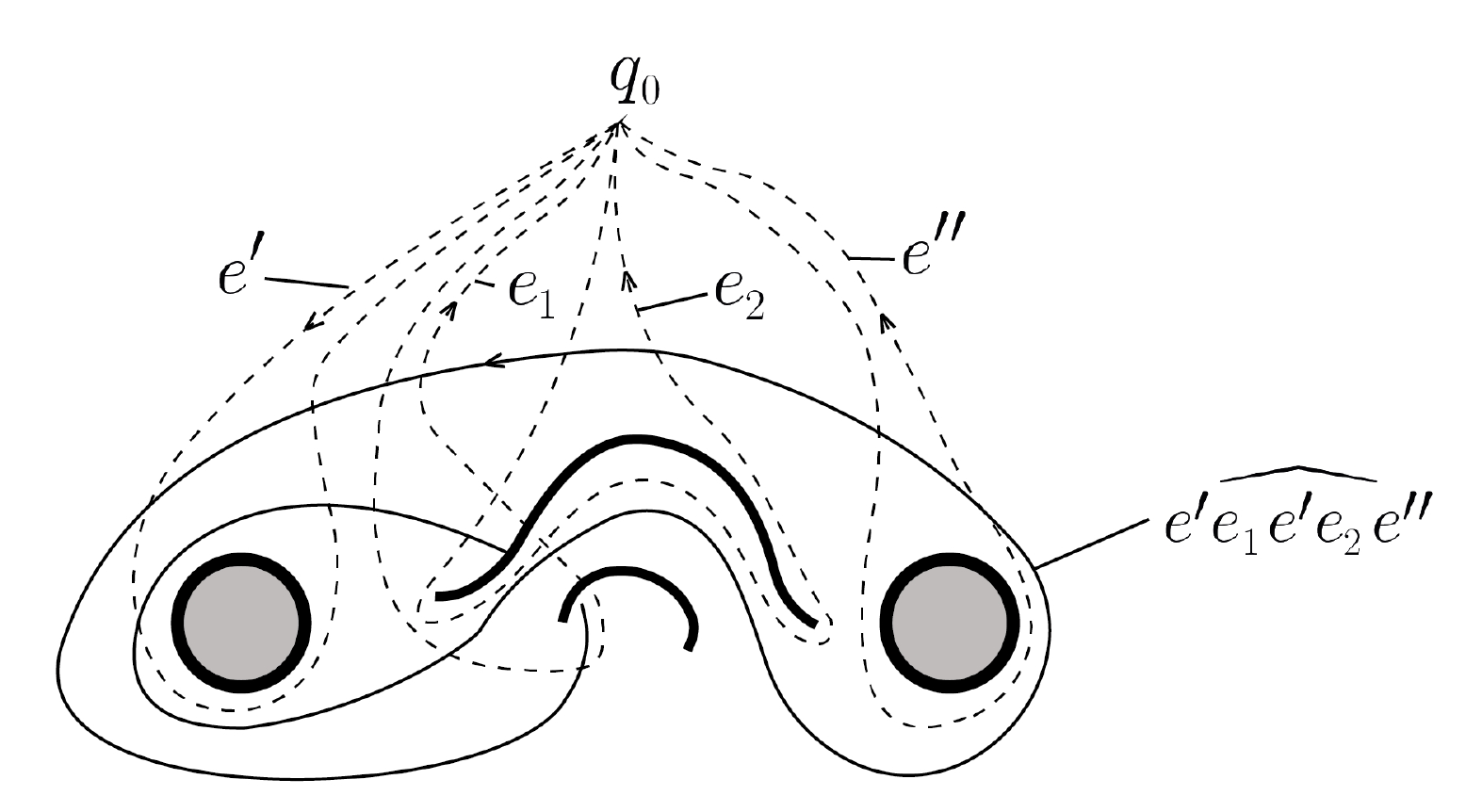}
\end{center}
\caption{A simple closed curve representing the conjugacy class of ${e'e_1 e' e_2 e''}$}\label{FigGrom1}
\end{figure}

Suppose $m>2$. Let $(e_1,e_2)$ be the pair of elements of $\mathcal{E}$ that corresponds to the first handle. We convert each unordered pair of different elements of $\mathcal{E}$
corresponding to holes to an ordered pair $(e',e'')$, so
that the conjugacy class of the element $e'e_1 e' e_2 e''$ can be represented by a simple closed curve in $X$ (see Figure \ref{FigGrom1}). We obtain no more than 
$\frac{1}{2}(m-1)(m-2)$ elements of $\pi_1(X,x_0)$.
They will be contained in $\mathcal{E}'$.

Suppose $g(X)=0$. If $m=1$ the fundamental group is trivial. If $m=2$ it has one generator. The set $\mathcal{E}'$ will be equal to the set $\mathcal{E}$ consisting of the generator. Let $m\geq 3$. The element $e_m\stackrel{def}=(\prod _{j=1}^{m-1}e_j)^{-1}$ is represented by a simple closed curve that surrounds $\mathcal{C}_m$ positively (i.e. walking along the curve the set $\mathcal{C}_m$ is on the left).
All elements $e\in \mathcal{E}\cup \{e_m\}$ will be contained in $\mathcal{E}'$. Further,
each unordered pair of different elements of $\mathcal{E}\cup \{e_m\}$ is converted into an ordered pair $(e',e'')$ and
the product $e'e''$ will be contained in $\mathcal{E}'$. Finally,
we convert each unordered triple of different elements of $\mathcal{E}\cup \{e_m\}$ into an ordered triple $(e',e'',e''')$ so that $\widehat{e'\,e''\,e'''}$ can be represented by a simple closed curve.
The product $e'e''e'''$
will be contained in $\mathcal{E}'$.

We described all elements of $\mathcal{E}'$. The conjugacy class of each of them can be represented by a simple closed curve. We estimate now the total number $N$ of elements of $\mathcal{E}'$.

Let first $g>0$. With $x=2g+m-1\geq 2$ the number $N$ does not exceed
\begin{equation*}
\begin{cases}
x+g+ \frac{1}{2} x(x-1)-g + 4g (x-2), & 0< m\leq 2,\\
x+g+ \frac{1}{2} x(x-1)-g+  4g (x-2)+ \frac{1}{2}(m-1)(m-2),& \quad \;\;\; m>2.\\
\end{cases}
\end{equation*}

Since $2g\leq x$, for $m\leq 2$ the sum
does not exceed $x+x(x-1) +2x(x-2)$.
Hence, $N\leq 3 x^2 - 4x < x^3$.
For $m>2$ we use the inequality $m-1= x-2g \leq x-2$, hence $\big(4g+\frac{1}{2}(m-2)\big)(x-2)\leq 2 x(x-2)$.
Again $N\leq x^3$.

Let $g=0$. If $m=2$ the set $\mathcal{E}'$ contains a single element, and $N=1=(m-1)^3=x^3$.
If $x=m-1\geq 2$ the number $N$ is not bigger than $x+1 + \frac{1}{2}(x+1)x + \frac{1}{6}(x+1)x(x-1)= \frac{1}{2}(x+1)(x+2)+ \frac{1}{6} x(x^2-1)=\frac{1}{6}x^3 +\frac{1}{2}x^2 +(\frac{3}{2}-\frac{1}{6})x+1$ and since $x\geq 2$, also in this case $N\leq x^3$.
We estimated
the number of elements of $\mathcal{E}'$
from above by $x^3=(2g+m-1)^3$.

\medskip

For the proof of Theorem \ref{prop0} we need two lemmas.
We denote by $a_1$ ($a_2$, respectively) the generator
of $\pi_1(\mathbb{C}\setminus\{-1,1\},0)$ that
is represented by a curve with base point $0$ that surrounds $-1$ ($1$, respectively) positively. We call $a_1$ and $a_2$ the standard generators of $\pi_1(\mathbb{C}\setminus\{-1,1\},0)$. The fundamental group $\pi_1(\mathbb{C}\setminus\{-1,1\},f(x_0))$ with base point $x_0$ is isomorphic to the fundamental group $\pi_1(\mathbb{C}\setminus\{-1,1\},0)$ with base point $0$ by an isomorphism that is determined up to conjugation.
The restriction of $f$ to an annulus in $X$ that represents $\hat e$ has the Gromov-Oka property if and only if
the conjugacy class $\widehat{f_*(e)}$ of the image $f_*(e)$ has infinite conformal module.
The conjugacy classes of elements of $\pi_1(\mathbb{C}\setminus\{-1,1\},0)$ with infinite conformal module are represented by integer powers of the elements
\begin{equation}\label{eq8.19''}
a_1,\; a_2,\;{\mbox{ or}}\;\, (a_1a_2)^{-1}\,.
\end{equation}

\begin{lemm}\label{lem20}
Let $X$ be an oriented smooth finite open surface and $\mathcal{E}'$ the set of elements of the fundamental group $\pi_1(X,x_0)$ chosen above. Assume the restriction of a continuous mapping $f:X\to \mathbb{C}\setminus \{-1,1\}$ to each annulus in $X$ representing the conjugacy class $\hat e$ of an element $e\in \mathcal{E'}$ has the Gromov-Oka property. Then either the image of the monodromy homomorphism $f_*:\pi_1(X,x_0)\to \pi_1(\mathbb{C}\setminus \{-1,1\},f(x_0))\cong
\pi_1(\mathbb{C}\setminus \{-1,1\},0)$ is contained in the group generated by a conjugate of a power of one of the elements of the form \eqref{eq8.19''},
or $X$ is equal to the oriented $2$-sphere denoted by $S^2$
with at least three holes removed and the mapping $f$ is homotopic to
a mapping that extends to a
homeomorphism $F:S^2\to \mathbb{P}^1$
that maps three points of $S^2$ contained in pairwise different holes to the three points $-1,1,$ and $\infty$, respectively.
\end{lemm}

\begin{lemm}\label{lem21} Let $X$ be a finite open Riemann surface with only thick ends.
For each homomorphism $h:\pi_1(X,x_0)\to \pi_1(\mathbb{C}\setminus \{-1,1\},0)$ with image generated by a power of one of the elements of the form \eqref{eq8.19''}
there exists a holomorphic mapping $f:X\to  \mathbb{C}\setminus \{-1,1\}$, whose monodromy homomorphism $f_*:\pi_1(X,x_0)\to \pi_1(\mathbb{C}\setminus \{-1,1\},0)$ is conjugate to $h$.

Let $X$ be the Riemann sphere $\mathbb{P}^1$ with at least three holes (maybe, of first kind). Suppose $F:\mathbb{P}^1\to \mathbb{P}^1$ is an
orientation preserving homeomorphism,
that maps three points of $\mathbb{P}^1$ contained in pairwise different holes to the three points $-1,1,$ and $\infty$, respectively. Then
the restriction $F\mid X$ is homotopic to a mapping that extends to a
conformal mapping $\mathbb{P}^1\to \mathbb{P}^1$ that maps $X$ to $\mathbb{C}\setminus\{-1,1\}$.
\end{lemm}
Recall that a continuous mapping $f:X\to \mathbb{C}\setminus\{-1,1\}$ is irreducible if it is not homotopic to a mapping with image in a punctured disc contained in $\mathbb{C}\setminus\{-1,1\}$ (the puncture may be equal to $\infty$). Each such mapping defines a special $(0,4)$-bundle $(X\times \mathbb{P}^1, {\rm pr}_1, \mathbold{E},X)$, where $\mathbold{E}$ intersects the fiber over $x$ along the set $\{x\}\times \{-1,1,{f}(x),\infty\}$. The mapping $f$ is irreducible iff the bundle is irreducible.

The following facts are easy to see. For a mapping $f$ from a smooth finite open surface $X$ to $\mathbb{C}\setminus\{-1,1\}$ the image $f_*(X,x_0)$ is contained in a
group generated by a conjugate of a power of one of the elements of the form \eqref{eq8.19''} if and only if $f$ is reducible.

The restriction $F\mid X$ of a homeomorphism $F:\mathbb{P}^1\to \mathbb{P}^1$,
with $F(X)\subset \mathbb{C}\setminus\{-1,1\}$ is irreducible.

\medskip

\noindent {\bf Proof of Theorem \ref{prop0}}.
If a mapping $f:X\to \mathbb{C}\setminus \{-1,1\}$ has the Gromov-Oka property then by Lemma \ref{lem0''} the restriction of $f$ to each annulus in $X$ representing the conjugacy class of an element of $\mathcal{E}'$ has the Gromov-Oka property.

By Theorem \ref{thmEl-1} the Lemma \ref{lem20} and the Lemma \ref{lem21} imply, that vice versa, if $X$ has positive genus and the restriction of a continuous mapping $f:X\to \mathbb{C}\setminus \{-1,1\}$ to each annulus in $X$ representing the conjugacy class $\hat e$ of an element $e\in \mathcal{E'}$ has the Gromov-Oka property, then $f$ has the Gromov-Oka property on $X$.

If $g=0$ a continuous orientation preserving mapping $f:X\to \mathbb{C}\setminus \{-1,1\}$ has the Gromov-Oka property on $X$, if its restriction to each annulus in $X$ representing the conjugacy class $\hat e$ of an element $e\in \mathcal{E'}$ has the Gromov-Oka property.
\hfill $\Box$

\medskip

\noindent {\bf Remark}. Notice that the two lemmas describe the homotopy classes of all mappings with the Gromov-Oka property from finite open smooth surfaces to the twice punctured complex plane.

The mappings with the Gromov-Oka property are the reducible mappings, and in case $X$ is the oriented two-sphere with at least three holes there are also irreducible homotopy classes of mappings with the Gromov-Oka property. Each consists of mappings with the following property. For each orientation preserving homeomorphism $\omega:X\to \omega(X)$ onto a Riemann surface $\omega(X)$ (maybe, of first kind) $f\circ \omega^{-1}$ is homotopic to a holomorphic mapping that extends to a conformal mapping $\mathbb{P}^1\to \mathbb{P}^1$ that maps three points in different holes to $-1,1$ and $\infty$, respectively (in some order depending on the class). These homotopy classes are the only irreducible homotopy classes with the Gromov-Oka property, and they are the only homotopy classes of mappings $X\to \mathbb{C}\setminus \{-1,1\}$ that contain a holomorphic mapping for any conformal structure on $X$ including conformal structures of first kind. Indeed,
if the image of a non-trivial monodromy homomorphism $f_*$ is generated by a conjugate of a power of one of the elements of the form \eqref{eq8.19''} the mapping $f$ cannot be homotopic to a holomorphic mapping for a conformal structure of first kind (see the proof of statement $(4)$ of Theorem \ref{thmEl.0}).

\medskip

\noindent {\bf Proof of Lemma \ref{lem20}.}
Identifying the fundamental groups $\pi_1(\mathbb{C}\setminus\{-1,1\},0)$ and $\pi_1(\mathbb{C}\setminus\{-1,1\},f(x_0))$ by a fixed isomorphism, we obtain for each $e\in\mathcal{E}'$ the equation $f_*(e)=w_e^{-1}\,b_e\,w_e$ for an element
$w_e\in \pi_1(\mathbb{C}\setminus\{-1,1\},0)$ and an element $b_e$ that is a power of an element of the form \eqref{eq8.19''}.
The images $f_*(e_{2j-1})$ and $f_*(e_{2j})$ of two elements corresponding to a handle commute. This is a particular case of the situation treated in the proof of Theorem \ref{thm8.1}, but can be seen easily directly.
Indeed, the sum of exponents of the terms of a word representing a commutator is equal to zero.
Since the commutator $[f_*(e_{2j-1}),f_*(e_{2j})]=f([e_{2j-1},e_{2j-1}])$ must be a power of a conjugate of an element of the form \eqref{eq8.19''} it must be equal to the identity.

A word in the generators of a free group is called reduced, if neighbouring terms are powers of different generators. We will identify elements of a free group (in particular of  $\pi_1(\mathbb{C}\setminus\{-1,1\},0)$) with reduced words in generators of the group. A word is called cyclically reduced, if either the word consists of a single term, or it has at least two terms and the first and the last term of the word are powers of different generators.

Suppose for two generators $e^{(1)}, e^{(2)} \in \mathcal{E} \subset \mathcal{E}'$
we have $b_{e^{(\ell)}}=a_{j_{\ell}}^{k_{\ell}}$, $\ell=1,2,$ where each $a_{j_{\ell}}$ is one of the generators $a_1$, $a_2$ of $\pi_1(\mathbb{C}\setminus\{-1,1\},0)$. Put $w={w_{e^{(2)}}}(w_{e^{(1)}})^{-1}$. The monodromy  $f_*(e^{(1)}e^{(2)})=f_*(e^{(1)})f_*(e^{(2)})$ is conjugate to $a_{j_1}^{k_1} w^{-1} a_{j_2}^{k_2} w$. The element $a_{j_1}^{k_1} w^{-1} a_{j_2}^{k_2} w$ is conjugate to $a_{j_1}^{k_1} {w'}^{-1} a_{j_2}^{k_2} w'$ for an element $w'= a_{j_2}^{-k'_2}w a_{j_1}^{-k'_1}
\in \pi_1(\mathbb{C}\setminus\{-1,1\},0)$ that is either equal to the identity or
it can be written as a reduced word that starts with a power of $a_{j_1}$ and ends with a power of $a_{j_2}$. (In the case $a_{j_1}= a_{j_2}$ we allow that $w'$ is equal to a power of the other generator of $\pi_1(\mathbb{C}\setminus\{-1,1\},0)$.)
If $w'$ is not the identity then $a_{j_1}^{k_1} {w'}^{-1} a_{j_2}^{k_2} w'$ can be written as cyclically reduced word that contains powers of different sign of generators of $\pi_1(\mathbb{C}\setminus\{-1,1\},0)$. Any cyclically reduced word that represents a conjugate of a power of an element of the form \eqref{eq8.19''} contains only powers of equal sign of the generators. Hence $w'={\rm Id}$, and $w=a_{j_2}^{k'_2} a_{j_1}^{k'_1}$. Put $w_{{e^{(1)}{e^{(2)}}}}= a_{j_2}^{-k'_2} w_{e^{(2)}} = a_{j_1}^{k'_1} w_{e^{(1)}} $. Then $w_{e^{(1)}e^{(2)}} f_*(e^{(1)}) w_{e^{(1)}e^{(2)}}^{-1}= a_{j_1}^{k_1}$ and $w_{e^{(1)}e^{(2)}} f_*(e^{(2)}) w_{e^{(1)}e^{(2)}}^{-1} = a_{j_2}^{k_2}$.

As a corollary we see that $f_*(e^{(1)}e^{(2)})$ is conjugate to $a_{j_1}^{k_1} a_{j_2}^{k_2}$. Hence, either at least one of the monodromies is the identity, or $j_1=j_2$, or $k_1=k_2=\pm 1$. Further, if the monodromies along a collection of elements of $\mathcal{E}$ are conjugate to a power of a single generator $a_j$,
there is a single element of $\pi_1(\mathbb{C}\setminus\{-1,1\},0)$  that conjugates each monodromy along an element of this collection to a power of $a_j$. Indeed, in this case we obtain
the following for each pair $e^{(1)},\,e^{(2)},$ of elements of this collection. If $f_*(e^{(\ell)})=w_{e^{(\ell)}}^{-1} a_j^{k_{\ell}} w_{e^{(\ell)}}$ for $\ell=1,2,$ then
$w_{e^{(2)}}=a_j^{k'} w_{e^{(1)}}$ for the same $a_j$, and $w_{e^{(1)}}$ conjugates both, $f_*(e^{(1)})$ and $f_*(e^{(2)})$, to a power of $a_j$.

The same arguments apply if for two elements $e^{(1)},e^{(2)}$ the monodromy $f_*(e^{(1)})$ is conjugate to a power of $a_j$ and the mondromy $f_*(e^{(2)})$ is conjugate to a power of $a_1a_2$. We replace the generators $a_1$, $a_2$ of the free group $\pi_1(\mathbb{C}\setminus\{-1,1\},0)$ by the generators
$A_1=a_j$ and $A_2=(a_1a_2)^{-1}$. Notice that $A_2$ can be considered as element of the fundamental group of the thrice punctured plane that is represented by loops surrounding $\infty$ positively.

We obtained the following.
If the monodromy $f_*(\tilde{e})$ along an element $\tilde{e}\in \mathcal{E}$ is conjugate to the $k$-th power of an element among  $a_1$, $a_2$, $a_1a_2$,
with $|k|>1$, then for any $e'\in \mathcal{E}$, $e'\neq \tilde e$, the monodromy $f_*(e')$ is a (trivial or non-trivial) power of the same element.

Moreover, we claim that if $g>0$ and the monodromy along
an element $e_j \in \mathcal{E},\, j\leq 2g,$ (i.e. $e_j$ belongs to a pair corresponding to a handle) is not trivial, then
all monodromies are powers of the same conjugate of an element of the form \eqref{eq8.19''}.
Indeed, by the preceding argument we may assume that all non-trivial monodromies along elements $e$ of $\mathcal{E}$ are conjugate to $v_e^{\pm 1}$ for an element $v_e$ of the form \eqref{eq8.19''}.
Consider two elements $e_{2j-1}, e_{2j}\in \mathcal{E}$ corresponding to a handle. Conjugating all monodromies by a single element, we suppose for instance that $f_*(e_{2j})=a_1$. Then $f_*(e_{2j-1})$ is a power of $a_1$, hence, by our assumption $f_*(e_{2j-1})$ is equal to $a_1^{\pm 1}$, or it is the identity.
If $f_*(e_{2j-1})$ equals $a_1$ or the identity, then
$f_*(e_{2j-1}\, e_{2j}^2\,e')$ can only be conjugate to a power of an element of the form \eqref{eq8.19''}, if $f_*(e')$ is a power of $a_1$.
If $f_*(e_{2j-1})$ is equal to $a_1^{-1}$, then
$f_*(e_{2j-1}\, e_{2j}^3\,e')$ can only be conjugate to a power of an element of the form \eqref{eq8.19''}, if $f_*(e')$ is a power of $a_1$.

The case when $f_*(e_{2j})$ equals $a_1^{-1}$, $a_2^{\pm 1}$, or $(a_1 a_2)^{\pm 1}$,
or when the role of $e_{2j-1}$ and $e_{2j}$ is interchanged, can be treated similarly.
The claim is proved.

Suppose $g(X)>0$ but all monodromies along elements of $\mathcal{E}$ corresponding to handles are trivial.
We claim that still all monodromies are powers of the same conjugate of an element of form \eqref{eq8.19''}. If the monodromies along all but possibly one element of $\mathcal{E}$ are trivial, there is nothing to prove. Hence, we may assume that $m>2$.
By the preceding argument we may assume that all non-trivial monodromies along elements of $\mathcal{E}$ are conjugate to an element of the form \eqref{eq8.19''} or are inverse to a conjugate of an element of the form \eqref{eq8.19''}. Let $e_1$ and $e_2$ be the elements of $\mathcal{E}$ corresponding to the first labeled handle.
Suppose there is an unordered pair of different elements of $\mathcal{E}$  so that the monodromy along each element of the pair is non-trivial. Note that  none of the elements equals
$e_1$ or $e_2$. Convert the unordered pair into an ordered pair $(e',e'')$, so that the conjugacy class $\reallywidehat{e'e_1e'e_2e''}$ can be represented by a simple closed curve.
Assume that $f_*(e')=a_1$.
If $f_*(e'')$ is not a power of $a_1$ then
$f_*(e'e_1e'e_2e'')= f_*(e')^2 f_*(e'')$ cannot be a power of an element of the form \eqref{eq8.19''}. The remaining cases when $f_*(e')$ is of the form \eqref{eq8.19''} or is inverse to an element of the form \eqref{eq8.19''} are treated in the same way. The claim is proved.

Let $g(X)=0$, i.e. $X$ equals the oriented $2$-sphere $S^2$ with holes.
Suppose $f:X\to \mathbb{C}\setminus\{-1,1\}$ is a continuous mapping whose restrictions to each annulus representing an element $e\in\mathcal{E}'$ have the Gromov-Oka property.
Assume that the monodromies are not all conjugate to powers of a single element of the form \eqref{eq8.19''}. Then
no monodromy is conjugate to the $k$-th power of an element of the form  \eqref{eq8.19''} with $|k|>1$, and there are at least three holes, and
at least two elements $e_{j'}$ and $e_{j''}$ in $\mathcal{E}$ with monodromies $f_*(e_{j'})$ and $f_*(e_{j''})$ being non-trivial powers of conjugates of different elements of the form \eqref{eq8.19''}. By the previous arguments we may assume after conjugating all monodromies by a single element of $\pi_1(\mathbb{C}\setminus\{-1,1\},0)$ that the momodromies along two generators in $\mathcal{E}$ are equal to different elements of the form \eqref{eq8.19''}, or to the inverses of different elements of the form \eqref{eq8.19''}.

Recall that the element $e_m= (\prod_{j=1}^{m-1} e_j)^{-1} \in \pi_1(X,x_0)$ is represented by a loop with base point $x_0$ that surrounds the last hole $\mathcal{C}_m$
counterclockwise. The product $\prod _{j=1}^m f_*(e_j)$ is equal to the identity. Hence, since
the product of two different elements of the form \eqref{eq8.19''}
cannot be equal to the identity,
there must be an integer number $j''' \in [1,m]$ different from $j'$ and $j''$,
for which $f_*(e_{j'''})\neq {\rm Id}$.

The monodromies along two different elements from $\mathcal{E}\cup\{e_m\}$ cannot be non-trivial powers of the same element of the form \eqref{eq8.19''}.
Otherwise there would be an ordered triple of different elements of $\mathcal{E}\cup\{e_m\}$ whose product is in $\mathcal{\mathcal{E}'}$  but the monodromy along the product cannot be a conjugate of a power of an element of the form \eqref{eq8.19''}.
Hence the monodromy along at most three elements of $\mathcal{E}\cup\{e_m\}$ is nontrivial.
Hence, there are exactly three elements $e_{j'}$, $e_{j''}$, and $e_{j'''}$ among the $e_k,\, k=1,\ldots,m$, with non-trivial monodromy.
The monodromies along two of them are equal to either $a_1$ and $a_2$, respectively, or to $a_1^{-1}$ and $a_2^{-1}$, respectively. Order the three elements by $(e',e'',e''')$ so that the product is in $\mathcal{E}'$.
After a cyclic permutation which does not change the conjugacy class of the product, we may assume that the monodromy along $e'''$ is not equal to a power of an $a_j$. Then the ordered triple of mondromies along $(e',e'',e''')$ is either $(a_1,a_2,(a_1a_2)^{-1})$, or $(a_2,a_1,(a_2a_1)^{-1})$, or $(a_1^{-1},a_2^{-1},a_2a_1)$, or $(a_2^{-1},a_1^{-1},a_1 a_2)$. The monodromies $(a_1,a_2,(a_1a_2)^{-1})$ are represented by curves that surround $(-1,1,\infty)$ positively, the monodromies
$(a_1^{-1},a_2^{-1},a_2a_1)$ are represented by curves that surround $(-1,1,\infty)$ negatively, and similarly for the remaining cases.

The case $(a_1,a_2,(a_1a_2)^{-1})$
for the monodromies along $(e',e'',e''')$  corresponds to the homotopy class containing the following mapping. Let $\mathcal{C}',\mathcal{C}'',\mathcal{C}'''$ be the holes of $X$ corresponding to $e',e'', e'''$, respectively.
Take points $p'\in \mathcal{C}',\, p''\in \mathcal{C}'',$ and $p'''\in \mathcal{C}'''$, respectively. Denote by $F$ an orientation preserving
homeomorphism from $S^2$ onto $\mathbb{P}^1$,
that maps  $p'$ to $-1$, $p''$ to $1$, and $p'''$ to $\infty$. It is straightforward to check that $F$ has the required monodromies along all generators. Hence, $f$ is homotopic to $F\mid X$. The case 
$(a_2,a_1,(a_2a_1)^{-1})$
for the monodromies along $(e',e'',e''')$ is similar. In these two cases
for any conformal structure $\omega:X\to \omega(X)$, including conformal structures of first kind, the mapping $(F|X)_{\omega}$ is homotopic to a
mapping hat extends to $\mathbb{P}^1=\omega(X)^c$ as a conformal self-homeomorphism of $\mathbb{P}^1$ that maps $X$ into $\mathbb{C}\setminus\{-1,1\}$.

If $f$ has monodromies $(a_1^{-1},a_2^{-1},a_2a_1)$ along $(e',e'',e''')$,
then $f$ is homotopic to $F\mid X$ for an orientation reversing
homeomorphism $F$ from $S^2$ onto $\mathbb{P}^1$,
that maps  $p'$ to $-1$, $p''$ to $1$, and $p'''$ to $\infty$. The case $(a_2^{-1},a_1^{-1}, a_1a_2)$ is similar.
In these cases for any conformal structure $\omega:X\to \omega(X)$, including conformal structures of first kind, the mapping $(F|X)_{\omega}$ is homotopic to a
mapping hat extends to $\mathbb{P}^1=\omega(X)^c$ as an anti-conformal self-homeomorphism of $\mathbb{P}^1$ that maps $X$ into $\mathbb{C}\setminus\{-1,1\}$. But the mapping $f$ does not have the Gromov-Oka property.

The latter fact can be seen as follows.
Assume the contrary. Let $\omega_n:X\to \omega_n(X)$ be a sequence of conformal structures of second kind on $X$ such that $X_n\stackrel{def}= \omega_n(X)$ can be identified with an increasing sequence of domains in $\mathbb{C}\setminus\{-1,1\}$ whose union equals $\mathbb{C}\setminus\{-1,1\}$. We may choose the $\omega_n$ uniformly converging on compact subsets of $X$.
We identify the fundamental groups of $X_n$ and of $\mathbb{C}\setminus\{-1,1\}$ by the isomorphism induced by inclusion.

Consider the case when $f$ has monodromies $(a_1^{-1},a_2^{-1},a_2a_1)$ along $(e',e'',e''')$.
By our assumption for each $n$ the mapping $f\circ \omega_n^{-1}$ is homotopic to a holomorphic mapping $f_n:X_n\to \mathbb{C}\setminus\{-1,1\}$ with monodromies $a_1^{-1}$ along $a_1$ and $a_2^{-1}$ along $a_2$.
By Montel's Theorem there is a subsequence $f_{n_k}$ that converges locally uniformly on $\mathbb{C}\setminus \{-1,1\}$. The limit function $F$ cannot be a constant (including $-1$, $1$, or $\infty$), since for curves $\gamma_j$ representing $a_j$, $j=1,2$, the set $f_{n_j}(\gamma_1)\cup f_{n_j}(\gamma_2)$ separates $-1$, $1$ and $\infty$.
The limit function $F$ is a holomorphic mapping from
$\mathbb{C}\setminus \{-1,1\}$ to itself.
Hence it
extends to a meromorphic function on $\mathbb{P}^1$, and defines  therefore a branched covering of $\mathbb{P}^1$. This is impossible, since the curves $\gamma_j$ are mapped to curves that surround $-1$ (and $1$, respectively) negatively.
\hfill $\Box$

\medskip

\noindent {\bf Proof of Lemma \ref{lem21}.}
Let $X$ be a finite open Riemann surface with only thick ends, and let $h:\pi_1(X,x_0)\to \Gamma$ be a homomorphism into the group $\Gamma$ generated by an element of $\pi_1(\mathbb{C}\setminus\{-1,1\},0)$ of the form \eqref{eq8.19''}. 
The Riemann surface
$X$ is conformally equivalent to a domain on a closed Riemann surface $X^c$ such that all connected components of its complement are closed topological discs. Hence, there exists an open Riemann surface $X_1\subset X^c$ such that $X$ is relatively compact in $X_1$. Let $\omega:X\to X_1$ be a homeomorphism from $X$ onto $X_1$ that is homotopic on $X$ to the inclusion $X\hookrightarrow X_1$.
Identify the fundamental groups $\pi_1(X,x_0)$ and $\pi_1(X_1,\omega(x_0))$ by the mapping induced by $\omega$. For $e\in \pi_1(X_1, \omega(x_0))$ we denote by $\sigma_e$ the covering transformation for the universal covering $\tilde{X}_1\to X_1$ that corresponds to $e$ under an isomorphism from the fundamental group onto the group of covering transformations (see \cite{Fo}).
For each collection of integer numbers $k_j,\, j=1,\ldots , m+2g-1,$ there exists a holomorphic function $\tilde{ F}:\tilde{X}_1\to \mathbb{C} $ on the universal covering $\tilde{X}_1\cong \mathbb{C}_+$ such that for each generator $e_j\in \mathcal{E}$ of $\pi_1(X,x_0)$ (identified with the respective generator of $\pi_1(X_1, {\omega}(x_0))$  the equality  $\tilde{F}(\sigma_{e_j}(\tilde{z}))=
k_j + \tilde{F}(\tilde{z}),\, z\in \tilde{X},$ holds. Such a function can be obtained as in the proof of Lemma \ref{lem*}. The function $e^{2\pi i {\tilde F}}$ descends to a holomorphic function $F_1$ on $X_1$.
For each $j$ and a representative $\gamma_j:[0,1]\to X_1$ of $e_j$ the closed curve $F_1\circ \gamma_j$ in $\mathbb{C}\setminus \{0\}$ represents ${\sf e}^{k_j}$ for the generator ${\sf e}$ of the fundamental group of $\mathbb{C}\setminus \{0\}$ with base point $F_1(x_0)$.
The restriction of $F_1$ to $X\Subset X_1$ is bounded. Hence, for a positive number $C$ the mapping $-1+ \frac{1}{C} F_1$ maps $X$ into $\mathbb{C}\setminus \{-1,1\}$ so that the monodromy along each $e_j$ equals $a_1^{k_j}$ (after identifying $\pi_1(\mathbb{C}\setminus \{-1,1\},-1+ \frac{1}{C} F_1(x_0))$
with  $\pi_1(\mathbb{C}\setminus \{-1,1\},0)$ by a suitable isomorphism).
The case when the monodromies are powers of other elements of form \eqref{eq8.19''} is treated by composing $F_1$ with a suitable conformal self-mapping of the Riemann sphere that permutes the points $-1$, $1$, and $\infty$.

The second statement of Lemma \ref{lem21} is clear.
\hfill $\Box$

\section{Smooth elliptic fiber bundles and differentiable families of complex manifolds.}
\label{sec:9.4}
\index{fiber bundle !  elliptic}

Recall that Kodaira (see \cite{Kodai}) defines a complex analytic family (also called holomorphic family) $\mathcal{M}_t, \,t
\in \mathcal{B},$ of compact complex manifolds over the base $\mathcal{B}$
as a triple $(\mathcal{M}, \mathcal{P}, \mathcal{B})$, where $\mathcal{M}$ and $\mathcal{B}$
are complex  manifolds, and $\mathcal{P}$ is a proper holomorphic submersion with
$\mathcal{P}^{-1}(t) = \mathcal{M}_t$. Thus, each holomorphic genus $\sf{g}$ fiber bundle
over a Riemann surface is a complex analytic family of compact Riemann surfaces of genus
$\sf{g}$.

We will consider holomorphic $(\sf{g},\sf{m})$-bundles as complex analytic families of
Riemann surfaces of genus $\sf{g}$ equipped with $\sf{m}$ distinguished points, for short, as complex analytic (or holomorphic) families of Riemann surfaces of type $({\sf g,m})$,
or holomorphic $(\sf{g},\sf{m})$-families.

Kodaira (s. \cite{Kodai}) defines a differentiable family of compact complex manifolds as
a smooth fiber bundle with compact fibers and the following additional structure. The fibers
are equipped with complex structures that depend smoothly on the parameter and induce on
each fiber the smooth fiber structure. The formal definition of a differentiable family of
Riemann surfaces of type $(\sf{g},\sf{m})$ is the following.

\begin{defn}\label{def4}
A differentiable
family of Riemann surfaces $\mathcal{M}_t,\, t \in \mathcal{B}$, of type
$(\sf{g},\sf{m})$ is a tuple $(\mathcal{M},\, \mathcal{P},\,\mathbold{E},\, \mathcal{B})$, where $\mathcal{M}$ and $\mathcal{B}$
are oriented $C^{\infty}$ manifolds, and $\mathcal{P}$ is a smooth proper orientation preserving submersion with the
following property. For each $t \in \mathcal{B}$ the fiber $\mathcal{P}^{-1}(t)$, denoted by
$\mathcal{M}_t$, is a compact Riemann surface of genus $\sf{g}$, 
such that the complex structure
induces the differentiable structure defined on $\mathcal{P}^{-1}(t)$ as a submanifold of
$\mathcal{M}$. $\mathbold{E}\subset \mathcal{M}$ is a smooth submanifold of $\mathcal{M}$ that intersects each fiber $\mathcal{M}_t$ along a set $E_t$ of $\sf{m}$ distinguished points. Moreover, the complex structures of the $\mathcal{M}_t$ depend smoothly on the parameter
$t$, more precisely, there is an open locally finite cover $U_j,\, j=1,2,\ldots$, of
$\mathcal{M}$  and complex valued $C^{\infty}$ functions $(z_1^j,\ldots, z_n^j)$ on each
$U_j$ (with $n$ the complex dimension of $\mathcal{M}_t$) such that for each $j$ and $t$
these functions define holomorphic coordinates on $U_j \cap \mathcal{M}_t$.
 \end{defn}
We will also call such families smooth families of Riemann surfaces of type $(\sf{g},\sf{m})$, or smooth $(\sf{g},\sf{m})$-families.

Two smooth families of Riemann surfaces of type $(\sf{g},\sf{m})$ will be called isomorphic (as families of Riemann surfaces) if they are isomorphic as smooth bundles (smoothly isomorphic, for short) and there is a bundle isomorphism that is holomorphic on each fiber.

We will show in this section that each smooth elliptic fiber bundle with a section is smoothly isomorphic (equivalently, isotopic) to a bundle that carries the structure of a differentiable
family of Riemann surfaces of type $(1,1)$, also called a differentiable family of
closed Riemann surfaces of genus $1$ with a smooth section.

Each compact Riemann surface of genus $1$ is conformally equivalent to the quotient
$\mathbb{C} \diagup \Lambda$ for a lattice $\Lambda=a\mathbb{Z}+b\mathbb{Z}$ where $a$ and
$b$ are two complex numbers that are independent over the real numbers. A torus given in the
form  $\mathbb{C} \diagup \Lambda$ will be called here a canonical torus.
The torus $\mathbb{C} \diagup (\mathbb{Z}+i \mathbb{Z})$ will be called standard.

A family of lattices $\Omega \ni x \to \Lambda(x)$ on a smooth manifold $\Omega$ is called
smooth if for each $x_0 \in X$
there is a neighbourhood $U(x_0) \subset \Omega$ such that the lattices
$\Lambda (x)$ can be written as
\begin{equation}\label{eqE110a}
\Lambda(x) = a(x){\mathbb Z}  + b(x){\mathbb Z}\,.
\end{equation}
for smooth functions $a$ and $b$ on $U(x_0)$.
The family is called holomorphic if it can be locally represented by
\eqref{eqE110a} with holomorphic functions $a$ and $b$.

If on $U(x_0)$ we have $\widetilde\Lambda (x) = \Lambda (x),\,$ with
$\widetilde\Lambda (x) = \widetilde a(x){\mathbb Z} + \widetilde
b(x){\mathbb Z}$ for other smooth functions $\widetilde a$ and
$\widetilde b$, then there is a matrix $A= \begin{pmatrix} \alpha & \beta \\
\gamma & \delta \end{pmatrix} \in {\rm SL}_2 ({\mathbb Z})\,$ not depending on $x$,
such that
\begin{equation}\label{eq2a}
\widetilde a(x) = \alpha \, a(x) \,+\gamma \,b(x),\;\;
\widetilde b(x) = \beta \,a(x) \,+ \delta\, b(x)\,.
\end{equation}
Indeed, $\Lambda(x)$ and $\widetilde{\Lambda}(x)$ are equal lattices
and $(a(x),b(x))$ and $(\widetilde a(x), \widetilde b(x))$ are pairs
of generators of the lattice. Put  $B(x)= \begin{pmatrix} {\rm{Re}}\,a(x) & {\rm{Re}}\,b(x) \\
{\rm{Im}}\,a(x) & {\rm{Im}}b(x) \end{pmatrix}$ and
$\widetilde{ B}(x)= \begin{pmatrix} {\rm{Re}}\,\widetilde a(x) & {\rm{Re}}\,\widetilde b(x) \\
{\rm{Im}}\,\widetilde a(x) & {\rm{Im}} \, \widetilde b(x)
\end{pmatrix}\,.$
Then the real linear self-map of $\mathbb{C}$, defined by $B(x)$ (by
$\widetilde {B}(x)$, respectively) takes $1$ to $a(x)$ and $i$ to $b(x)$ ($1$ to $\widetilde{a}(x)$ and $i$ to $\widetilde{b}(x)$, respectively). Hence, both maps take
the standard lattice $\mathbb{Z} + i \mathbb{Z}$
onto $\Lambda (x) = \widetilde \Lambda (x)$. Then $ B^{-1}(x)\circ\widetilde B (x)$ maps the standard lattice onto itself, hence, it is
in ${\rm SL}_2 ({\mathbb Z})$. Since $B(x)$ and $\widetilde B (x)$
depend continuously on $x$, the matrices $ B ^{-1}(x) \circ \widetilde B(x) $ depend
continuously on $x$ and, hence, are equal to a matrix $A= \begin{pmatrix} \alpha &
\beta \\
\gamma & \delta \end{pmatrix} \in {\rm SL}_2 ({\mathbb Z})\,$ that does not depend on $x$.
Equation \eqref{eq2a} is satisfied for the entries of this matrix.

Two smooth families of lattices $\Lambda_0(x)$ and $\Lambda_1(x)$ depending on a parameter $x$ in a
smooth manifold $\Omega$ are called isotopic, if for an interval $I\supset [0,1]$ there is a smooth family of lattices
${\Lambda}(x,t)$,  $(x,t)\in \Omega \times I$, such that
${\Lambda}(x,0)=\Lambda_0(x)$ and  ${\Lambda}(x,1)=\Lambda_1(x)$ for all $x\in \Omega$.

Let $X$ be a finite open Riemann surface and
$U$ an open subset of $X$ where there are smooth function
$a(x)$ and $b(x)$ such that
$\Lambda(x)=a(x)\mathbb{Z}+b(x)\mathbb{Z},\, x\in U$. Then the family $\Lambda= \Lambda(x),\, x \in U,$ defines a free and properly discontinuous group action
\begin{equation}\label{eqE111a}
U\times \mathbb{C} \ni (x,\zeta) \to \big(x, \zeta+ a(x)n+ b(x) m\big),\,n,m\in \mathbb{Z}\,,
\end{equation}
on $U\times \mathbb{C}$. The quotient of  $(U\times \mathbb{C})$ by this action depends only on $\Lambda$, not on the choice of the generators $a(x)$ and $b(x)$ of the lattice.
Denote by $(U\times \mathbb{C})\diagup \Lambda$ the quotient of $U\times
\mathbb{C}$ by this action \eqref{eqE111a}. Let
\begin{equation}\label{eqE111b}
\mathcal{P}_{U,\Lambda}:(U\times \mathbb{C})\diagup \Lambda \to U
\end{equation}
be the mapping that takes the equivalence class of $(x,\zeta)$ to $x$.
Consider in each fiber $\mathcal{P}_{\Lambda}^{-1}(x)$ the distinguished point $\;s_{\Lambda,x}=(x,0)\diagup \Lambda(x)\;$. The mapping $\;x\to s_{\Lambda,x},\, x \in U,\; $ is smooth, hence defines a smooth section. Put $\;\;\;\,\mathbold{s}_{\Lambda,U}\stackrel{def}=\cup_{x\in U}\{s_{\Lambda,x}\}$.
Then
the tuple $\mathfrak{F}_{\Lambda,U}\stackrel{def}= \big((U\times \mathbb{C}) \diagup \Lambda,\mathcal{P}_{U,\Lambda}, \,\mathbold{s}_{\Lambda,U}\,,U\big)$
defines a smooth $(1,1)$-bundle over $U$. Moreover, $\mathfrak{F}_{\Lambda,U}$
is equipped with the structure of a differentiable family of Riemann surfaces of
type $(1,1)$.
If the family $\Lambda$ of lattices on $X$ is holomorphic, then
$\mathfrak{F}_{\Lambda,U}$ is
a holomorphic $(1,1)$-bundle, equivalently, a
holomorphic family of Riemann surfaces of type $(1,1)$.

Take an
open covering of $X$ by sets $U_j$ on which there are smooth functions
$a_j(x)$ and $b_j(x)$ such that
$\Lambda(x)=a_j(x)\mathbb{Z}+b_j(x)\mathbb{Z},\, x\in U_j$.
If $U_j\cup U_k \neq \emptyset$, we glue the restrictions of $\mathfrak{F}_{\Lambda,U_j}$ and
$\mathfrak{F}_{\Lambda, U_k}$ to $U_j\cup U_k$ together by using the identity map in each fiber.
We obtain a smooth  family of Riemann surfaces of type $(1,1)$ on $X$, which we denote by
$\mathfrak{F}_{\Lambda,X}=(\mathcal{X}_{\Lambda},\mathcal{P}_{\Lambda}, \,\mathbold{s}_{\Lambda}\,,X)$.

If the family $\Lambda$ of lattices on $X$ is holomorphic, then
$\mathfrak{F}_{\Lambda}$ is
a holomorphic $(1,1)$-bundle, equivalently, a
holomorphic family of Riemann surfaces of type $(1,1)$.

The following lemma holds.
\begin{lemm}\label{lem1}
Each smooth $(1,1)$-bundle $\mathfrak{F}$ over a smooth finite open surface $X$ is
smoothly isomorphic (equivalently, isotopic) to a bundle of the form
\begin{equation}\label{eqEl10b}
\mathfrak{F}_{\Lambda,X}=(\mathcal{X}_{\Lambda},\mathcal{P}_{\Lambda}, \,\mathbold{s}_{\Lambda}\,,X)\,.
\end{equation}
\end{lemm}
Recall that the bundle \eqref{eqEl10b} carries the structure of a smooth family of Riemann surfaces of type $(1,1)$.

\medskip

\noindent {\bf Proof.} We need to find a smooth bundle of the form \eqref{eqEl10b} whose monodromy homomorphism is conjugate to that of the original bundle $\mathfrak{F}$. This can be done as in Section \ref{sec:9.2}. We represent $X$ as the union of an open disc $D$
containing the base point of $X$, and a collection of $\ell$  attached 
bands $V_j$ that are relatively closed in $X$ and correspond to the generators of the fundamental group of $X$. Let as in Section \ref{sec:9.2} the set $\tilde{U}$ be a simply connected domain on the universal covering of $X$ such that each point of $D$ is covered $\ell+1$ times and each other point of $X$ is covered once. Label the preimages of $D$ under the projection $\tilde{U}\to X$ by $\tilde{D}_j,\; j=0,1,\ldots \ell,$ so that the preimage $\tilde{V}_j$ of the band $V_j,\,j=1,\ldots,\ell,$ is attached to $\tilde{D}_0$ and $\tilde{D}_j$. Denote by $\mathfrak{m}_j$ the monodromy of $\mathfrak{F}$ along the generator $e_j$ of the fundamental group of $X$ that corresponds to the $j$-th band $V_j$.

Map
the fiber of $\mathfrak{F}$ over the base point diffeomorphically onto the standard torus
$\mathbb{C}\diagup(\mathbb{Z}+i \mathbb{Z})$ so that the distinguished point is mapped to $0\diagup (\mathbb{Z}+i \mathbb{Z})$. By an isomorphism between mapping class groups
we identify each mapping class $\mathfrak{m}_j, \, j=1,\ldots,\ell,$ with a mapping class denoted by the same letter on the standard torus
$\mathbb{C}\diagup(\mathbb{Z}+i \mathbb{Z})$ with distinguished point $0\diagup (\mathbb{Z}+i \mathbb{Z})$.
Represent the new class $\mathfrak{m}_j$
by a mapping
$ \varphi_j^{-1}:\mathbb{C}\diagup(\mathbb{Z}+i \mathbb{Z})\toitself$,                    such
that $\varphi_j$ lifts  to  a real linear self-mapping $\widetilde{\varphi}_j$ of $\mathbb{C}$
that maps the lattice $\mathbb{Z}+i \mathbb{Z}$ onto itself. In other words,
$\widetilde{\varphi}_j$ corresponds to a $2 \times 2$ matrix $A_j$ with
integral entries and determinant $1$, $\widetilde {\varphi}_j (x + iy) = A_j
\begin{pmatrix} x \\ y \end{pmatrix}$, $A_j \in {\rm SL}_2 ({\mathbb
Z})$.

Consider on each set $\widetilde{\Omega}_j\stackrel{def}=\tilde{D}_0\cup \tilde{V}_j\cup \tilde{D}_j$ (which is a simply connected domain)
a smooth family of real linear self-maps $(\widetilde\varphi_j)_z$, $z \in \widetilde{\Omega}_j$, of the
complex plane ${\mathbb C}$ such that
\begin{align}\label{eqEl10a}
(\widetilde\varphi_j)_z=\mbox{Id} \;\mbox{for $z\in \tilde{D}_0$ and}
(\widetilde\varphi_j)_z=\widetilde \varphi_j\; \mbox{for $z\in \tilde{D}_j$}\,,
\end{align}
and put
\begin{equation}\label{eqEl10a'}
\Lambda_j(z)=(\widetilde\varphi_j)_z(\mathbb{Z}+i\mathbb{Z}),\;z \in \widetilde{\Omega}_j\,.
\end{equation}
Glue the bundles $\big(\,(\widetilde{\Omega}_j\times \mathbb{C})\diagup \Lambda_j\,,\, \mathcal{P}_{\Lambda_j}\,,\, \mathbold{s}_{\Lambda_j}\,,\, X\,\big)$ together
by identifying points $\tilde{z}_j\in \tilde{D}_j,\, j=0,\ldots \ell,$ that project to the same point in $D$,
and use the mapping $\varphi_j$ to glue the fiber
over each $\tilde{z}_j\in \tilde{D}_j$ to the fiber over $\tilde{z}_0\in \tilde{D}_0$, where $\tilde{z}_j$ and $\tilde{z}_0$ project to the same point in $D$.
We obtain a bundle
$\mathfrak{F}_{\Lambda,X}=\big(\,(X\times \mathbb{C})\diagup \Lambda\,,\, \mathcal{P}_{\Lambda}\,,\, \mathbold{s}_{\Lambda}\,,\, X\,\big)$ over $X$, which is isomorphic (as a smooth bundle) to
the original bundle over $X$, since the monodromy homomorphisms of the two bundles are conjugate to each other.

The lemma is proved. \hfill $\Box$
\medskip

By Lemma  \ref{lem1} the Problem \ref{prob9.1}  can be reformulated as follows.

\smallskip

\noindent  {\it Let $X$ be a finite open Riemann surface. Is a given
smooth family of Riemann surfaces of type $(1,1)$ on $X$
isotopic to a
complex analytic family of Riemann surfaces of type $(1,1)$?}

\section{Complex analytic families of canonical tori.      }
\label{sec:9.5}
Let now ${\mathfrak F} = ({\mathcal X} , {\mathcal P} , X)$ be a
{\it holomorphic} elliptic fiber bundle over a finite open Riemann
surface $X$.
The fiber bundle is, in particular, a smooth elliptic fiber bundle. For each disc $\Delta
\subset X$ there is a smooth family of diffeomorphisms $\varphi_t:S \to
\mathcal{P}^{-1}(t),\; t \in \Delta,$ from the reference Riemann surface $S=\mathbb{C}\diagup (\mathbb{Z}+i \mathbb{Z})$ of genus one onto the fiber over $t$.
Consider the Teichm\"uller class $[\varphi_t],\; t \in \Delta$.
The following Lemma \ref{lem3} states that these Teichm\"uller classes depend
holomorphically on the parameter.

\begin{lemm}\label{lem3}
Let $\mathfrak{F}$ be a holomorphic elliptic fiber bundle over a
Riemann surface $X$. For each small enough disc $\Delta \subset X$
there is a holomorphic map $z \to \tau (z)$, $z \in \Delta$, into the Teichm\"uller space
$\mathcal{T}(1,0)$ of the standard torus, such that each fiber ${\mathcal P}^{-1}(z) $ is conformally equivalent to $\mathbb{C}\diagup
(\mathbb{Z}+\tau(z)\mathbb{Z})$.
\end{lemm}

For convenience of the reader we will provide a proof. The key ingredient is a lemma of
Kodaira which we formulate now.

Let ${\mathcal X}$ and $X$ be complex manifolds and let ${\mathcal
P} : {\mathcal X} \to X$ be a proper holomorphic submersion such
that the fibers ${\mathcal X}_z = {\mathcal P}^{-1} (z)$, $z\in X$,
are compact complex manifolds of complex dimension $n$. For each $z
\in X$  we denote by $\Theta_z$ the sheaf of germs of holomorphic
tangent vector fields of the complex manifold ${\mathcal X}_z$.
Denote by $H^0 ({\mathcal X}_z , \Theta_z)$ the space of global
sections of the sheaf.

\begin{lemm}\label{lem4}
{\rm (Kodaira, \cite{Kodai}, Lemma 4.1, p. 204.)} If the dimension
$d \stackrel{def}= {\rm dim} (H^0 ({\mathcal X}_z , \Theta_z))$ is independent of $z$ then
for any small enough (topological) ball $\Delta$ in $X$ there is for
each $z \in \Delta$ a basis $(v_1 (z) , \ldots , v_d (z))$ of $H^0
({\mathcal X}_z , \Theta_z)$ such that $(v_1 (z) , \ldots , v_d
(z))$  depends holomorphically on $z \in \Delta$.
\end{lemm}
For a proof we refer to \cite{Kodai}.
Let $m$ be the dimension of the complex manifold $X$ in the statement of Kodaira's Lemma. The condition that the $v_j (z)$, $j =
1,\ldots , d$, depend holomorphically on $z$, means the following.
Let $U_{\alpha} \subset {\mathcal X}$ be a small open subset of
${\mathcal X}$ on which there are local holomorphic coordinates
$(\zeta_1^{\alpha} , \ldots , \zeta_n^{\alpha} , z_1 , \ldots ,
z_m)$, where $z = (z_1 , \ldots , z_m) \in \Delta$ and for fixed $z
= (z_1 , \ldots , z_m)$ the $\zeta^{\alpha} = (\zeta_1^{\alpha} ,
\ldots , \zeta_n^{\alpha})$ are local holomorphic coordinates on the
fiber over $z$. In these local coordinates the vector field $v_j$,
$j = 1,\ldots , n$, can be written as
\begin{equation}\label{eqEl17}
v_j (z) = \sum_{k=1}^n v_{j \, k}^{\alpha}
(\zeta^{\alpha} , z) \frac{\partial}{\partial \zeta_k^{\alpha}} \, ,
\end{equation}
where the $v_{j \, k}^{\alpha}$ are holomorphic in $(\zeta^{\alpha}
, z)$. The condition does not depend on the choice of local
coordinates with the described properties.

We prepare now the proof of Lemma
\ref{lem3}.
Kodaira's Lemma applies in the situation of Lemma \ref{lem3}. Indeed, in the situation of Lemma \ref{lem3}
the base $X$ has complex dimension one and the fibers are compact Riemann surfaces of genus $1$.
The space of holomorphic sections
$H^0 ({\mathcal X}_z , \Theta_z)$ of holomorphic tangent vector
fields of each fiber ${\mathcal X}_z$ has complex dimension one.
Indeed, let $\mathbb{T}$ be a compact Riemann surface of genus $1$. The
complex structure of $\mathbb{T}$ may be given by a system of local holomorphic coordinates
with transition functions $\zeta(z)=z+c$ for complex constants $c$ (so called flat
coordinates). A holomorphic vector field (equivalently, a holomorphic section in the
holomorphic tangent bundle) on $\mathbb{T}$ assigns to each chart  $\mathbb{C} \supset U_j \to V_j \subset \mathbb{T}$ a function $v_j(z_j),\, z_j\in U_j$, such that if $V_j\cap V_k\neq
\emptyset$ then $v_j(z_j(z_k)) \frac{dz_k}{dz_j}=v_k(z_k)$. Since the complex structure on $\mathbb{T}$ is given by a system of flat local
coordinates, the equality $\frac{dz_k}{dz_j}\equiv 1$ holds for each $j,k$ with $V_j\cap V_k\neq \emptyset$.
Hence, the $v_j(z_j)$ define a holomorphic function $v$ on $\mathbb{T}$. Each holomorphic
function on a compact complex manifold is constant. Hence,
the space of holomorphic tangent vector fields to any compact Riemann surface of genus one has complex dimension $1$.
In particular, for each $z \in X$ the space $H^0
({\mathcal X}_z , \Theta_z)$ is generated by a single holomorphic
tangent vector field $v(z)$ on ${\mathcal X}_z$.
By Kodaira's Lemma
$v(z)$ may be chosen to depend holomorphically on $z \in \Delta$ for
each small enough disc $\Delta$ in $X$. We obtain a holomorphic vector field
$v$ on ${\mathcal X}_{\Delta}\stackrel{def}= {\mathcal P}^{-1} (\Delta) $ whose
restriction to each fiber ${\mathcal X}_z$ equals $v(z)$.

Consider the (holomorphic) universal covering space $\widetilde{\mathcal
X}_{\Delta}$ of ${\mathcal X}_{\Delta} $. Denote by $p$ the covering map.
Consider the triple $(\widetilde{\mathcal
X}_{\Delta}, \mathcal{P}\circ p, \Delta)$.

\begin{lemm}\label{lem5}
$(\widetilde{\mathcal X}_{\Delta},\mathcal{P}\circ p, \Delta) $ is holomorphically isomorphic to the trivial
holomorphic fiber
bundle over $\Delta$ with fiber ${\mathbb C}$.
\end{lemm}

\noindent {\bf Proof of Lemma \ref{lem5}.} For each $z\in \Delta$ the fiber $\mathcal{P}^{-1}(z)$ of the bundle
$({\mathcal X}_{\Delta},\mathcal{P}, \Delta) $ is
a compact Riemann surface of genus $1$.
Its preimage under $p$ is the set $(\mathcal{P}\circ p)^{-1}(z)$ which is the fiber of the bundle $(\widetilde{\mathcal X}_{\Delta},\mathcal{P}\circ p, \Delta) $ over $z$.
The set $\mathcal{P}^{-1}(z)$
is a complex one-dimensional submanifold of ${\mathcal X}_{\Delta}$ which is the zero set of the holomorphic function $\mathcal{P}-z$.
The set $(\mathcal{P}\circ p)^{-1}(z)$ is a complex one-dimensional submanifold of $\widetilde{\mathcal X}_{\Delta}$ which is the zero set of the holomorphic function
$\mathcal{P}\circ p-z$.
Hence, the restriction $p\mid (\mathcal{P}\circ p)^{-1}(z): (\mathcal{P}\circ p)^{-1}(z)\to \mathcal{P}^{-1}(z)$ defines a holomorphic covering. Indeed, since $p:\widetilde{\mathcal X}_{\Delta}\to {\mathcal X}_{\Delta}$ is a covering, each point $z'$  of $\mathcal{P}^{-1}(z)$ has a neighbourhood $U_{z'}\subset  {\mathcal X}_{\Delta}$ such that $p^{-1}(U_{z'})$ is the disjoint union of open sets $\tilde U _{z'}^k\subset \widetilde{\mathcal{X}}$ and $p$ maps each $\tilde U _{z'}^k$ biholomorphically onto $U_{z'}$. Then $p$  maps each $\tilde U _{z'}^k \cap (\mathcal{P}\circ p)^{-1}(z) $ conformally onto $ U_{z'}\cap \mathcal{P}^{-1}(z)$.

Each fiber $(\mathcal{P}\circ p)^{-1}(z)$ is simply connected, since any loop in the fiber $(\mathcal{P}\circ p)^{-1}(z)$ is contractible in $\widetilde{\mathcal X}_{\Delta}$, and the bundle  $(\widetilde{\mathcal X}_{\Delta},\mathcal{P}\circ p, \Delta) $ is smoothly isomorphic to the trivial bundle $(\Delta\times (\mathcal{P}\circ p)^{-1}(z), {\rm pr}_1, \Delta)$.

Since for each $z\in \Delta$ the covered manifold $\mathcal{P}^{-1}(z) $
is a Riemann surface of genus $1$, the covering manifold $(\mathcal{P}\circ p)^{-1}(z) $ is conformally equivalent to the complex plane $\mathbb{C}$.
Hence, the
triple $(\widetilde{\mathcal X}_{\Delta} , {\mathcal P} \circ p ,
\Delta)$ is a holomorphic fiber bundle with fiber ${\mathbb C}$.

Take a holomorphic section $s(z),\, z\in \Delta,$ of the mapping $\widetilde{\mathcal
X}_{\Delta} \overset{{\mathcal P} \, \circ \,
p}{-\!\!\!-\!\!\!-\!\!\!\longrightarrow} \Delta$. It exists after,
perhaps, shrinking $\Delta$. Let $v$ be the holomorphic
vector field on ${\mathcal X}_{\Delta}$ from Kodaira's Lemma, and
$\widetilde v$ its lift to the universal
covering $\widetilde{\mathcal X}_{\Delta}$. Define a mapping
${\mathcal G} : \Delta \times {\mathbb C} \to \widetilde{\mathcal
X}_{\Delta}$ by
\begin{equation}\label{eqEl18}
{\mathcal G} (z,\zeta) = \gamma_{s(z)} (\zeta) \, ,
\quad (z,\zeta) \in \Delta \times {\mathbb C} \, , 
\end{equation}
where for each $z$ the mapping $\gamma_{s(z)}$ is the solution of
the holomorphic differential equation
\begin{equation}\label{eqEl19}
\gamma'_{s(z)} (\zeta) = \widetilde v (\gamma_{s(z)}
(\zeta)) \, , \quad \gamma_{s(z)} (0) = s(z) \in ({\mathcal P} \circ
p)^{-1} (z) \, , \quad \zeta \in {\mathbb C} \, .
\end{equation}
Since $\widetilde v$ is tangential to the fibers we have the
inclusion
\begin{equation}\label{eqEl20}
\gamma_{s(z)} ({\mathbb C}) \subset ({\mathcal P}
\circ p)^{-1} (z) \cong {\mathbb C} \quad \mbox{for each} \quad z
\in \Delta \, .
\end{equation}
For each $z$ the solution exists for all $\zeta \in {\mathbb C}$
since the restriction of $\widetilde v$ to the fiber over $z$ is the
lift of a vector field on a closed torus. The mapping
\begin{equation}\label{eqEl21}
(z,\zeta) \to \gamma_{s(z)} (\zeta) \, , \quad z \in
\Delta \, , \quad \zeta \in {\mathbb C} \, ,
\end{equation}
is a local holomorphic diffeomorphism.
By the Poincar\'{e}-Bendixson Theorem for each $z$ it maps
${\mathbb C}$ one-to-one onto an open subset of the fiber $({\mathcal P} \circ p)^{-1}
(z) \cong {\mathbb C}$, hence, it maps ${\mathbb C}$ {\it onto}
$({\mathcal P} \circ p) ^{-1} (z)$. Hence, ${\mathcal G}$ defines a
holomorphic isomorphism of the trivial bundle $(\Delta \times\mathbb{C},{\rm pr}_1,
\Delta)$ onto the bundle $(\widetilde{\mathcal
X}_{\Delta},{\mathcal P} \, \circ \,p,\, \Delta)$.
Lemma \ref{lem5} is
proved.  \hfill $\Box$

\medskip

\noindent {\bf Proof of Lemma \ref{lem3}.}
Consider the covering transformations of the covering
$\widetilde{\mathcal X}_{\Delta} \to {\mathcal X}_{\Delta}$.
In terms of the isomorphic bundle $(\Delta \times \mathbb{C},{\rm pr}_1,
{\mathbb C})$
the restrictions of the covering transformations to each
fiber $\{z\} \times {\mathbb C}$, are translations, hence, the
covering transformations have the form
\begin{equation}\label{eqEl22}
\psi_{n,m} (z,\zeta) =\big(z, \zeta + n \, a(z) + m \,
b(z) \big)\, , \quad (z,\zeta) \in \Delta \times {\mathbb C} \, ,
\end{equation}
for integral numbers $n$ and $m$. Here $a(z)$ and $b(z)$ are complex
numbers which are linearly independent over ${\mathbb R}$ and depend
on $z$. Since the covering transformations $\psi_{1,0}$ and $\psi_{0,1}$ are
holomorphic, the numbers $a(z)$ and $b(z)$ depend holomorphically on
$z \in \Delta$. After perhaps interchanging $a(z)$ and $b(z)$ we may assume that
$\mbox{Re}\frac{a(z)}{b(z)}  >0\,.$
Hence, $\tau (z) = \frac{a(z)}{b(z)}$ depends
holomorphically on $z \in \Delta$ and can be interpreted as a holomorphic mapping from
$\Delta$ into the Teichm\"uller space (which is identified with the upper half-plane
$\mathbb{C}_+$). Lemma \ref{lem3} is proved.  \hfill $\Box$


\begin{cor}\label{cor1}
Let $\mathfrak{F}$ be a holomorphic elliptic fiber bundle over a
simply connected Riemann surface $X$, or let $\mathfrak{F}$ be a holomorphic elliptic fiber bundle over a finite open Riemann surface $X$, such that all monodromies are equal to the identity. Then
there is a holomorphic map $z \to \tau (z)$, $z \in X$, into the Teichm\"uller space $\mathcal{T}(1,0)$, such
that each fiber ${\mathcal P}^{-1}(z) $ is conformally equivalent to $\mathbb{C}\diagup
(\mathbb{Z}+\tau(z)\mathbb{Z})$.
\end{cor}

\noindent {\bf Proof.} Let first $X$ be simply connected. Consider a smooth family of diffeomorphisms $\varphi_z: S \to
\mathcal{P}^{-1}(z),\ z \in X$, where $S=\mathbb{C}\diagup (\mathbb{Z}+ i \mathbb{Z})$ is the standard torus. Then the Teichm\"uller
classes $\tau(z) \stackrel{def}{=}[\varphi_z], z \in X,$ depend smoothly on $z$ and
$\mathcal{P}^{-1}(z)$ is conformally equivalent to $\mathbb{C}\diagup
(\mathbb{Z}+[\varphi_z]\mathbb{Z})$. By Lemma \ref{lem3}
$[\varphi_z]$ depends holomorphically on $z$. Indeed, on any small enough disc $\Delta
\subset X$ there is a holomorphic map $\tau_{\Delta}$ to the Teichm\"uller space such that
$\mathcal{P}^{-1}(z)$ is conformally equivalent to $\mathbb{C}\diagup (\mathbb{Z}+ \tau_{\Delta}(z)
\mathbb{Z})$. Hence, $[\varphi_z]=\varphi_{\Delta}^*(\tau_{\Delta}(z)),\; z \in \Delta,\;$ for a modular
transformation $\varphi^*_{\Delta}$ of the Teichm\"uller space. The corollary follows from the fact that modular
transformations are biholomorphic maps on the Teichm\"uller space.

Let $X$ be a finite open Riemann surface, and let all monodromies of the bundle $\mathfrak{F}$ be equal to the identity. Consider a lift $\tilde{\mathfrak{F}}$ of the bundle to the universal covering $\tilde X \overset{p}{-\!\!\!\longrightarrow}    X$. There exists a holomorphic mapping $\tilde{\tau}:\tilde{X}\to \mathcal{T}(1,0)$ such that for each $\tilde{x}\in\tilde{X}$
the value $\tilde{\tau}(\tilde{x})$ represents the conformal class of the fiber of $\tilde{\mathfrak{F}}$ over $\tilde{x}$ (which is equal to the conformal class of the fiber of $\mathfrak{F}$ over $p(\tilde{x})$). Let $e$ be an element of the fundamental group of $X$ with base point $x_0$. For the associated covering transformation of the covering $\tilde X  \overset{p}{-\!\!\!\longrightarrow}    X$ (denoted also by $e$) we let $\varphi_e$ be a self-homeomorphism of the fiber of $\mathfrak{F}$ over $x_0$, that represents the monodromy mapping class of the bundle along $e$. Denote by $\psi_e^*$
the modular transformation on $\mathcal{T}(1,0)$ corresponding to $\varphi_e$. Then
$\tilde{\tau}(e(\tilde{x}))=\psi_e^*(\tilde{\tau}(\tilde{x}))$ for each $e$.
Since all monodromies of the bundle $\mathfrak{F}$ are trivial, each $\psi_e^*$ is equal to the identity, and therefore $\tilde{\tau}$ descends to a well-defined mapping $\tau$ on $X$ with the required property.
\hfill $\Box$

\medskip

The following proposition holds.

\begin{prop}\label{prop1} Each holomorphic $(1,1)$-bundle
$(\mathcal{X},\mathcal{P},\mathbold{s}, X)$ over a finite open Riemann surface $X$ is holomorphically isomorphic to a
holomorphic bundle of the form
\begin{equation}\label{eqE111e}
(\mathcal{X}_{\Lambda},\mathcal{P}_{\Lambda}, \,\mathbold{s}_{\Lambda}\,,X)\,.
\end{equation}
for a holomorphic family of lattices $\Lambda=\{\Lambda(x)\}_{x\in X}$.
\end{prop}
Here as before, $\mathcal{X}_{\Lambda}=(X\times \mathbb{C})\diagup \Lambda$, and
$\mathcal{P}_{\Lambda}$ assigns  the point $x\in X$ to the class in the quotient containing the point $(x,\zeta)$.
The value at the point $x\in X$ of the holomorphic section $\mathbold{s}_{\Lambda}$ of the bundle \eqref{eqE111e}
is the class in the quotient containing $(x,0)$.
\medskip

\noindent {\bf Proof.} For each $x\in X$ we let $\Delta_x$ be a topological disc, $x \in\Delta_x \subset X$, for which Kodaira's Lemma holds. As in Lemma \ref{lem5} we let
$p_{\Delta_x}:\widetilde{\mathcal{X}}_{\Delta_x} \to {\mathcal{X}}_{\Delta_x} $ be the universal covering of
$\mathcal{X}_{\Delta_x}=\mathcal{P}^{-1}(\Delta_x)$. Let $\tilde{\mathbold{s}}_{\Delta_x}  $ be a lift of ${\mathbold{s}}_{\Delta_x} \stackrel{def}=  \mathbold{s}\cap \mathcal{X}_{\Delta_x}$ to $\widetilde{\mathcal{X}}_{\Delta_x}$.
By Lemma \ref{lem5} there exists (after perhaps shrinking $\Delta_x$)
a holomorphic bundle isomorphism
\begin{equation}\label{eq6}
\Big(\,\widetilde{\mathcal{X}}_{\Delta_x},\, (\mathcal{P}\circ p_{\Delta_x})| \widetilde{\mathcal{X}}_{\Delta_x}\,, \tilde{\mathbold{s}}_{\Delta_x},\,\Delta_x\,\Big) \to \Big(\,\Delta_x\times \mathbb{C},\, {\rm pr}_1,\,\Delta_x\times \{0\},\, \Delta_x\,\Big)\,.
\end{equation}
The bundle isomorphism is given by the identity mapping on $\Delta_x$ and a holomorphic mapping $\tilde{\Phi}_{\Delta_x}$ from $\widetilde{\mathcal{X}}_{\Delta_x}$ to $\Delta_x\times \mathbb{C}$, which maps the fiber of the first bundle over each $x'\in \Delta_x$ to the fiber of the second bundle over $x'$ and maps the section of the first bundle to the section of the second bundle.

Cover $X$ by a locally finite set of such discs $\Delta_j$.
For each $\Delta_j$ we consider the covering transformations of the projection $p_{\Delta_j}:\widetilde{\mathcal{X}}_{\Delta_j}\to \mathcal{X}_{\Delta_j}$. The bundle isomorphism
conjugates the group of these covering transformations to a group of fiber preserving transformations of $\Delta_j\times \mathbb{C}$ with free and properly discontinuous action.
The conjugated group acts on each fiber $\{x\}\times \mathbb{C}$ as a lattice
$\Lambda_j(x)$, and the lattices depend holomorphically on $x\in \Delta_j$.
Taking the quotient, we obtain for each $j$ a holomorphic bundle isomorphism
\begin{equation}\label{eq7}
\big(\mathcal{X}_{\Delta_j},\, \mathcal{P}_{\Delta_j},\; \mathbold{s}_{\Delta_j},\, \Delta_j\,\big)\to \big((\Delta_j\times \mathbb{C})\diagup {\Lambda}_{j},\; {\rm pr}_1,\, (\Delta_j\times\{0\})\diagup {\Lambda}_{j},\; \Delta_j\,\big).
\end{equation}
The bundle isomorphism  \eqref{eq7} is given by the identity mapping from $\Delta_j$ to itself and by a biholomorphic mapping $\Phi_j:\mathcal{X}_{\Delta_j}\to
(\Delta_j\times \mathbb{C})\diagup {\Lambda}_{j}$ that maps the fiber of the first bundle over each point $x\in \Delta_j$ to the fiber of the second bundle over $x$ and maps the section of the first bundle to the section of the second bundle.

Take $j$ and $k$ so that $\Delta_j\cap \Delta_k\neq \emptyset$.
Then $\Phi_j\circ (\Phi_k)^{-1}$ is a fiber preserving biholomorphic mapping from  $((\Delta_j\cap \Delta_k)\times \mathbb{C})\diagup {\Lambda}_{k}$ onto
$((\Delta_j\cap \Delta_k)\times \mathbb{C})\diagup {\Lambda}_{j}$.
Moreover, $\Phi_j\circ (\Phi_k)^{-1}$ maps $((\Delta_j\cap \Delta_k)\times\{0\})\diagup {\Lambda}_{j}$ onto $((\Delta_j\cap \Delta_k)\times\{0\})\diagup {\Lambda}_{k}$.
The mapping $\Phi_j\circ (\Phi_k)^{-1}$ lifts to a fiber preserving biholomorphic mapping  $\tilde{\Phi}_{j,k}\stackrel{def}=\tilde{\Phi}_j\circ (\tilde{\Phi}_k)^{-1}$
of $ (\Delta_j\cap \Delta_k)\times \mathbb{C}$ onto itself, that takes $(\Delta_j\cap \Delta_k)\times\{0\}$ onto itself. Then
\begin{equation}\label{eqE111j}
\tilde{\Phi}_{j,k}(x,\zeta)= (x,\alpha_{j,k}(x)\cdot \zeta),\; (x,\zeta)\in (\Delta_j\cap \Delta_k)\times \mathbb{C}\,
\end{equation}
for a nowhere vanishing holomorphic function $\alpha_{j,k}$ on $(\Delta_j\cap \Delta_k)$.
Moreover, the mapping $\tilde{\Phi}_{j,k}$
takes each $\{x\}\times {\Lambda}_k(x)\subset \{x\}\times \mathbb{C}$ onto $\{x\}\times {\Lambda}_j(x)\subset \{x\}\times \mathbb{C}$.

The $\alpha_{j,k}$ form a Cousin II cocycle (see \cite{H1}) on $X$. Since $X$ is a finite open Riemann surface, the Cousin II problem is solvable (see \cite{Fo}), i.e. there exist holomorphic nowhere vanishing functions $\alpha_j$ on $\Delta_j$, such that on $(\Delta_j\cap \Delta_k)$ the equality $\alpha_{j,k}=\alpha_k\cdot \alpha_j^{-1}$ holds.
Then for each $j,k$ with $\Delta_j\cap \Delta_k\neq \emptyset$ the equality $\alpha_j(x)\Lambda_j(x)=\alpha_k(x)\Lambda_k(x),\; x\in \Delta_j\cap \Delta_k$ holds.
Hence, we obtain a well-defined holomorphic family of lattices $\Lambda$ on $X$, $\Lambda(x)=\alpha_j(x)\Lambda_j(x),\, x \in \Delta_j$.
For this family of lattices we consider the bundle
$\big((X\times\mathbb{C})\diagup \Lambda, \mathcal{P}_{\Lambda}, \mathbold{s}_{\Lambda},X\big)$
of the form \eqref{eqE111e}.

Let $c_j:(\Delta_j\times \mathbb{C})\diagup {\Lambda}_j\to (\Delta_j\times \mathbb{C})\diagup {\Lambda}$ be the biholomorphic mapping,
that lifts to the biholomorphic self-mapping $(x,\zeta)\to(x, \alpha_j(x)\,\zeta)$ of $\Delta_j\times \mathbb{C}$, which maps each $\{x\}\times\Lambda_j(x)$ onto $\{x\}\times\alpha_j(x)\Lambda_j(x)=\{x\}\times\Lambda(x)$.
For each non-empty $\Delta_j\cap \Delta_k$ we let $c_{j,k}$ be the biholomorphic mapping
from $((\Delta_j\cap \Delta_k)\times \mathbb{C})\diagup \Lambda_k)$ onto $((\Delta_j\cap \Delta_k)\times \mathbb{C})\diagup \Lambda_j)$ that lifts to the mapping \eqref{eqE111j}.

Put $\Phi'_j= c_j \Phi_j$ on $\mathcal{X}_{\Delta_j}$. Then on each non-empty intersection $(\Delta_j\cap \Delta_k)$ the equality $\Phi'_j (\Phi'_k)^{-1}=c_j\Phi_j (c_k\Phi_k)^{-1}= c_j (\Phi_j (\Phi_k)^{-1})c_k^{-1}= c_j c_{j,k} c_k^{-1}$ holds. Consider the lifts of the mappings $c_j,c_k,c_{j,k}$ to (complex linear) self-mappings of $\mathbb{C}$. Since $\alpha_{j,k}=\alpha_k (\alpha_j)^{-1}$ on $(\Delta_j\cap \Delta_k)$, we obtain $\Phi'_j=\Phi'_k$ on $\mathcal{X}_{\Delta_j}\cap \mathcal{X}_{\Delta_k}$. The mapping $\Phi':\mathcal{X}\to (X\times \mathbb{C})\diagup \Lambda$, for which $\Phi'(x)=\Phi'_j(x)$ for $x\in \mathcal{X}_{\Delta_j}$, is well defined and determines a holomorphic isomorphism from the original bundle to a bundle of the form \eqref{eqE111e}. Proposition \ref{prop1} is proved. \hfill $\Box$

\smallskip

Similar arguments as used in the proof of Proposition \ref{prop1} give the following statement.\\
\smallskip

\noindent {\bf Proposition $1'$.} {\it Each holomorphic elliptic fiber bundle over a finite open Riemann surface $X$ is holomorphically isomorphic to a holomorphic bundle that admits a holomorphic section.}\\

\smallskip

\noindent {\bf Proof.} Define the holomorphic bundle isomorphism \eqref{eq7} for some chosen locally defined holomorphic sections $\mathbold{s}_{\Delta_j}$. For each $j$ we obtain a holomorphic isomorphism $\tilde{\Phi}_j$ of the total space of the bundle on the left onto the total space of the bundle on the right of \eqref{eq7}, that maps fibers of the first bundle to fibers of the second one.

For $\Delta_j\cap \Delta_k\neq \emptyset$ the mapping
$\tilde{\Phi}_{j,k}=\tilde{\Phi}_j\circ (\tilde{\Phi}_k)^{-1}$
is a fiber preserving biholomorphic mapping
of $ (\Delta_j\cap \Delta_k)\times \mathbb{C}$ onto itself. Then the equation
\begin{equation}\label{eq10}
\tilde{\Phi}_{j,k}(x,\zeta)= (x, a_{j,k}(x) + \alpha_{j,k}(x)\cdot \zeta),\; (x,\zeta)\in (\Delta_j\cap \Delta_k)\times \mathbb{C}\,
\end{equation}
holds for a holomorphic function $a_{j,k}$ and a nowhere vanishing holomorphic function $\alpha_{j,k}$ on $(\Delta_j\cap \Delta_k)$.
Since $\tilde{\Phi}_{j,k}\,\tilde{\Phi}_{k,i}\, \tilde{\Phi}_{i,j}={\rm Id}$, we obtain
\begin{equation}\label{eq11}
a_{j,k}(x) +\alpha_{j,k}(x) a_{k,i}(x) +\alpha_{j,k}(x)\alpha_{k,i}(x) a_{i,j}(x) + \alpha_{j,k}(x)\alpha_{k,i}(x)\alpha_{i,j}(x)\cdot \zeta \,=\,\zeta\,.
\end{equation}
Hence, $\alpha_{j,k}(x)\alpha_{k,i}(x)\alpha_{i,j}(x) \equiv 1$, in other words, the $\alpha_{j,k}$ form a Cousin II cocycle. Since $X$ is an open Riemann surface the Cousin problem has a solution, i.e. there exist nowhere vanishing holomorphic functions $\alpha_j$ on $\Delta_j$ such that $\alpha_{j,k}=\frac{\alpha_k}{\alpha_j}$ on $(\Delta_j\cap \Delta_k)$.
Then by equation \eqref{eq11} $\alpha_j a_{j,k} + \alpha_k a_{k,i} + \alpha_i a_{i,j}=0$, in other words, the $\alpha_j a_{j,k}$ form a Cousin I cocycle. Since $X$ is a Stein manifold the Cousin problem is solvable, i.e. there exist
holomorphic functions $a_j$ on $\Delta_j$ such that $\alpha_j a_{j,k}=a_j-a_k$ on $(\Delta_j\cap \Delta_k)$. Consider the biholomorphic mapping
$\tilde{\Phi}_j^*:\tilde{\mathcal{X}}_{\Delta_j}\to \Delta_j\times \mathbb{C}$ which is the composition $\tilde{\Phi}_j^*=A_j\circ\tilde{\Phi}_j$ of $\tilde{\Phi}_j$ with the mapping $A_j:\Delta_j\times \mathbb{C}\, \toitself\,$, $A_j(x,\zeta)=(x, \alpha_j(x) \zeta -a_j(x))$. Since $A_k^{-1}(x,\zeta)=(x, \frac{1}{\alpha_k}\zeta +\frac{a_k}{\alpha_k})$ we obtain
\begin{align}\label{eq12}
\tilde{\Phi}_j^*(\tilde{\Phi}_k^*)^{-1} (x,\zeta)  =  A_j \tilde{\Phi}_j (\tilde{\Phi}_k)^{-1}A_k^{-1}(x,\zeta)
= &  A_j \tilde{\Phi}_j (\tilde{\Phi}_k)^{-1}(x, \frac{1}{\alpha_k}\zeta +\frac{a_k}{\alpha_k})\nonumber\\
=  A_j(x, a_{j,k} + \alpha_{j,k}(\frac{1}{\alpha_k}\zeta +\frac{a_k}{\alpha_k}))
= & A_j(x, a_{j,k} + \frac{1}{\alpha_j} (\zeta + a_k)) \nonumber\\
= (x,\alpha_j( a_{j,k} + \frac{1}{\alpha_j} (\zeta + a_k)) -a_j)
 = & (x,a_j-a_k +\zeta+a_k -a_j)=(x,  \zeta)\,.
\end{align}
We obtained biholomorphic mappings $\tilde{\Phi}_j^*:\tilde{\mathcal{X}}_{\Delta_j}\to \Delta_j\times \mathbb{C}$ such that on $\Delta_j\cap \Delta_k$ the mappings $\tilde{\Phi}_j^*$ and $\tilde{\Phi}_k^*$ coincide.
${\mathcal{X}}_{\Delta_j}$ is the quotient of $\tilde{\mathcal{X}}_{\Delta_j}$ by the group of covering transformations of the covering $p_{\Delta_j}: \tilde{\mathcal{X}}_{\Delta_j}\to\mathcal{X}_{\Delta_j}   $.
Conjugating for each $j$ the group of covering transformations 
by the mapping $\tilde{\Phi}_j^*$, we obtain a group acting on 
$\Delta_j\times \mathbb{C}$
which can be identified with a holomorphic family of lattices $\Lambda_j$ on $\Delta_j$.
The mapping $\tilde{\Phi}_j^*$ descends to a biholomorphic mapping $\Phi_j:\mathcal{X}_{\Delta_j}\to (\Delta_j\times \mathbb{C})\diagup \Lambda_j$. The lattices $\Lambda_j$ and $\Lambda_k$ coincide on $\Delta_j\cap\Delta_k$, and $\Phi_j=\Phi_k$ on $\Delta_j\cap\Delta_k$.
For the lattice $\Lambda$ on $X$ that equals $\Lambda_j$ on $\Delta_j$
we obtain a bundle
$((X\times \mathbb{C})\diagup \Lambda,\mathcal{P}_{\Lambda}, \,\mathbold{s}_{\Lambda}\,,X)$
of the form \eqref{eqE111e}, and the mappings $\Phi_j$ define a holomorphic bundle isomorphism from the original bundle to $((X\times \mathbb{C})\diagup \Lambda,\mathcal{P}_{\Lambda},X)$. The proposition is proved. \hfill $\Box$

\section{Special $(0,4)$-bundles and double branched coverings}
\label{sec:9.6}
Take any element $\mathring{E}=\{z_1,z_2,z_3\}\in C_3(\mathbb{C})\diagup \mathcal{S}_3$, and put $E\stackrel{def}=\mathring{E}\cup \{\infty\}$. The set
\begin{equation}\label{eqEL31a}
\mathring{Y}_E\stackrel{def}= \Big\{(z,w)\in \mathbb{C}^2: w^2 = 4 \, (z - z_1) (z-z_2)
(z-z_3)\Big\}
\end{equation}
is a one-dimensional complex submanifold of $\mathbb{C}^2$.
Each point of $\mathring{Y}_E$ has a neighbourhood in $\mathring{Y}_E$ on which one of the functions, $z$ or $w$, defines local holomorphic coordinates.
The mapping
\begin{equation}\label{eqEl31d}
\mathring{Y}_E\ni (z,w)\to z \in \mathbb{C}
\end{equation}
is a branched holomorphic covering of $\mathbb{C}$ with branch locus $\mathring{E}$.

Consider the $1$-point compactification $Y_E$ of $\mathring{Y}_E$. The complex structure on it is obtained as
follows. Let
$\mathring{Y}_E^{\,r}$ be the subset of $\mathring{Y}_E$ where $|z|>r$ for a large positive number $r$. On this set $|\frac{z}{w}|$ is small. Put $(\tilde{z},\tilde{w})=(\frac{1}{z},
\frac{z}{w})$ on a small neighbourhood in $\mathbb{C}^2$ of
$\mathring{Y}_E^{\,r}$. In these coordinates
the equation for $\mathring{Y}_E^{\,r}$ becomes
\begin{equation}\label{eqEl36b}
\tilde{w}^2= \tilde{z} \frac{1}{4 \prod_{j=1}^3(1- z_j\tilde{z})}\,,
\end{equation}
and $\mathring{Y}_E^{\,r}$ can be identified with the subset
\begin{align}\label{eqEl36a}
\left\{(\tilde{z},\tilde{w})\in \mathbb{C}^2,\, 0<|\tilde{z}|<\frac{1}{r}:\tilde{w}^2=
\tilde{z} \frac{1}{4 \prod_{j=1}^3(1-z_j\tilde{z})}\right\}
\end{align}
of $\mathbb{C}^2$. Adding the point $(\tilde{z},\tilde{w})=(0,0)\in \mathbb{C}^2$ to
$\mathring{Y}_E^{\,r}$  we obtain a complex manifold ${Y}_E^{\,r}$.
The two manifolds ${Y}_E^{\,r}$  and $\mathring{Y}_E$
form an open cover of the desired compact complex manifold $Y_E$.
Denote the point $(0,0)$ in coordinates $(\tilde{z},\tilde{w})$ on $Y^r_{E}$ by $s^{\infty}$.
Each point of ${Y}_E$ has
a neighbourhood where one of the functions $z$, $w$, $\tilde{z}$, or
$\tilde{w}$ defines local holomorphic coordinates. The holomorphic projection $Y_E\to
\mathbb{P}^1$ is correctly defined by $(z,w)\to z,\, (z,w)\in \mathring{Y}_E,\,$
and $(\tilde{z},\tilde{w})\to  \tilde{z},\,  (\tilde{z},\tilde{w})\in {Y}_E^{\,r}$.
We obtain a double branched covering  $Y_E\to \mathbb{P}^1$
over $\mathbb{P}^1$ with
branch locus equal to $E\stackrel{def}=\mathring{E} \cup \{\infty\}$. The manifold $Y_E$ is a closed Riemann surface of genus $1$. We will consider
it as a closed Riemann surface of genus $1$ with distinguished point being the preimage $s^{\infty}$ of $\infty$ under the branched
covering.
The set $\mathring{E}$ (considered as subset of $\mathbb{C}$) will be called the finite branch locus of
the covering $Y_E\to \mathbb{P}^1$.

Let $Y$ be a closed Riemann surface of genus $1$ with distinguished point $s$ and let $Y_1$ be equal to $\mathbb{P}^1$ with
set of distinguished points $\mathring{E}\cup \{\infty\}$ for $\mathring{E}\subset C_3(\mathbb{C})\diagup
\mathcal{S}_3$. Suppose ${\sf{P}}:Y\to Y_1$ is a double branched covering with branch
locus $\mathring{E}\cup \infty$ and ${\sf{P}}(s)=\infty$. A mapping class $\mathfrak{m}\in
\mathfrak{M}(Y;s,\emptyset)$ is called a lift of a mapping class $\mathfrak{m}_1\in
\mathfrak{M}(Y_1;\infty,\mathring{E})$ if there are representing homeomorphisms $\varphi \in
\mathfrak{m}$ and $\varphi_1 \in \mathfrak{m}_1$, such that $\varphi$ lifts $\varphi_1$,
i.e. $\varphi_1({\sf{P}}(\zeta))={\sf{P}}(\varphi(\zeta))$, $\zeta\in Y$.

Let $X$ be a finite open Riemann surface.
Recall that a special holomorphic (smooth, respectively) $(0,4)$-bundle $\mathfrak{F}_1$ over $X$ is a bundle of the form
$\big(X\times \mathbb{P}^1,{\rm pr}_1,{\mathbold{E}} ,\,X\big)$ for a complex (smooth, respectively) submanifold  $\mathbold{E}$ of $X\times \mathbb{P}^1$ which intersects each fiber $\{x\}\times\mathbb{P}^1 $ along a set of distinguished points $E_x= \{x\}\times (\mathring{E}_x\cup \{\infty\})$ with  $\mathring{E}_x\subset C_3(\mathbb{C})\diagup \mathcal{S}_3$. A smooth special $(0,4)$-bundle over $X$ carries in a natural way the structure of a family of Riemann surfaces of type $(0,4)$ depending smoothly on the parameter in $X$.
Let $\mathfrak{F}=\big(\mathcal{X},\mathcal{P}, \mathbold{s},X\big)$ be a
complex analytic (smooth, respectively) family of Riemann
surfaces of type $(1,1)$.
The family $\mathfrak{F}$ is called a double branched covering of the special
holomorphic (smooth, respectively)
$(0,4)$-bundle $\big(X\times \mathbb{P}^1,{\rm pr}_1,\mathbold{E} ,\,X\big)$
if there exists
a holomorphic (smooth, respectively) mapping ${\sf{P}}: \mathcal{X} \to   X\times\mathbb{P}^1$
that maps each fiber $\mathcal{P}^{-1}(x)$ of the $(1,1)$-family $\mathfrak{F}$
onto the fiber
$\{x\}\times \mathbb{P}^1$
of the $(0,4)$-bundle over the same point $x$, such that the
restriction
${\sf{P}}: \mathcal{P}^{-1}(x) \to \{x\}\times \mathbb{P}^1$ is a holomorphic double
branched covering with branch locus being the set $\{x\}\times (\mathring{E}_x\cup \{\infty\})$ of distinguished
points in the fiber $\{x\}\times \mathbb{P}^1$, and $\sf{P}$
maps the distinguished point $s_x$  in the fiber $\mathcal{P}^{-1}(x)$ over $x$ to
the point $\{x\}\times \{\infty\}$ in the fiber $\{x\}\times \mathbb{P}^1$ of the special $(0,4)$-bundle.
We will also write
$(X\times \mathbb{P}^1  , {\rm pr}_1  ,\mathbold{E}, X)={\sf{P}}((\mathcal{X},\mathcal{P},\mathbold{s},X))$, and call
the family $(\mathcal{X},\mathcal{P},\mathbold{s}, X)$ a lift of $(X\times\mathbb{P}^1, {\rm pr}_1 ,\mathbold{E} ,X)$.

Each special holomorphic 
$(0,4)$-bundle $\big(X\times \mathbb{P}^1,{\rm pr}_1,\mathbold{E} ,\,X\big)$
over a finite open Riemann surface  $X$ admits a double
branched covering by a complex analytic family of Riemann surfaces of type $(1,1)$.
Each special smooth
$(0,4)$-bundle over a differentiable manifold $X$ has a double branched covering by a
differentiable family of Riemann surfaces of type $(1,1)$.
Notice that the statement for the smooth case is true also for families over products $X\times I$ where $X$ is a finite open surface and $I$ is an interval, in other words, it is true for isotopies of bundles.

We prove the statement for the holomorphic case.
The smooth case is treated similarly as the holomorphic case.
Assume the sets $\mathring{ E}_x,\, x\in X,$ (with $\mathbold{E}\cap (\{x\}\times \mathbb{P}^1)=\{x\}\times (\mathring{E}_x \cup \{\infty\})$) are uniformly
bounded in $C_3(\mathbb{C})\diagup\mathcal{S}_3$.
Consider the set
\begin{equation}\label{eqEl37}
\mathring{\mathcal Y}_{\mathbold{E}} = \left\{ (x,z,w) \in X
\times {\mathbb C}^2 : w^2 = 4 \prod_{z_j \in \mathring{E}_x} (z
- z_j) \right\} \,,
\end{equation}
equipped with the structure of an embedded complex hypersurface  in $X\times\mathbb{C}^2$.
In a neighbourhood of a point $(x_0 , z_0 , w_0)$ on
$\mathring{\mathcal Y}_{\mathbold{E}}$ with $z_0 \notin \mathring{E}_{x_0}$ the pair $(x,z)$
defines holomorphic coordinates. If $z_0 \in \mathring{E}_{x_0}$ the pair
$(x,w)$ defines holomorphic coordinates in a neighbourhood of $(x_0
, z_0 , w_0)$ on $\mathring{\mathcal Y}_{\mathbold{E}}$.
The projection ${\mathcal P} :\mathring{\mathcal Y}_{\mathbold{E}} \to X$, ${\mathcal P}
(x,z,w) = x$ is holomorphic.

For some large positive number $r$ we may define the set
\begin{align}\label{eqEl36c}
{\mathcal Y}_{\mathbold{E}}^{\,r}\stackrel{def}=\left\{(x,\tilde{z},\tilde{w})\in X\times
\mathbb{C}^2,\, |\tilde{z}|<\frac{1}{r}:\tilde{w}^2= \tilde{z} \frac{1}{4\prod_{z_j \in \mathring{E}_x}
(1-{z}_j\tilde{z})}\right\}\,.
\end{align}
It is a complex hypersurface of complex dimension two of $X\times \mathbb{C}^2$.
Each of its points has a neighbourhood on which either $(x,\tilde{z})$ or
$(x,\tilde{w})$ defines holomorphic coordinates. The mapping $(x,\tilde{z},\tilde{w})\to x$ is
holomorphic.
The part
$\mathring{\mathcal{Y}}_{\mathbold{E}}^{\,r}\stackrel{def}=\big\{(x,\tilde{z},\tilde{w})
\in{\mathcal{Y}}_{\mathbold{E}}^{\,r}:
\tilde{z}\neq 0\big\}$ of ${\mathcal Y}_{\mathbold{E}}^{\,r}$ can be identified with the subset
$\big\{(x,z,w) \in \mathring{\mathcal{Y}}_{\mathbold{E}}: \, |z|>r\big\}$  of
$\mathring{\mathcal{Y}}_{\mathbold{E}}$ using the transition functions $(x,z)\to
(x,\tilde{z})=(x,\frac{1}{z})$, $(x,w)\to(x,\tilde{w})=(x,\frac{z}{w})$.
The sets $\mathring{\mathcal{Y}}_{\mathbold{E}}$ and ${\mathcal Y}_{\mathbold{E}}^{\,r}$
form an open cover of a complex manifold denoted by ${\mathcal Y}_{\mathbold{E}}$ equipped
with a proper holomorphic submersion $\mathcal{P}:{\mathcal Y}_{\mathbold{E}}\to X$, such that $\mathcal{P}^{-1}(x)$ is a torus for each $x\in X$.
We proved that $({\mathcal Y}_{\mathbold{E}},\mathcal{P},X)$ is a holomorphic elliptic bundle. Let $\mathbold{s}^{\infty}$ be the submanifold of $\mathcal{X}$
that intersects each fiber $\mathcal{P}^{-1}(x)$ along the
distinguished point $s_x^{\infty}\in \mathcal{Y}^r_{\mathbold{E}}$ that is written in coordinates on $\mathcal{Y}^r_{\mathbold{E}}$ as $(x,0,0)$.

Let ${\sf P}: \mathcal{Y}_{\mathbold{E}}\to X\times \mathbb{P}^1$ be the map that assigns to $(x,z,w)$ the point $(x,z)$ (and to  $(x,\tilde{z},\tilde{w})$ the point $(x,\tilde{z})$).
By the construction it is clear that the obtained $(1,1)$-bundle $({\mathcal
Y}_{\mathbold{E}},\mathcal{P},\,\mathbold{s}^{\infty},\,X)$ is a double branched covering of the $(0,4)$-bundle $\big(X\times \mathbb{P}^1,{\rm pr}_1,\mathbold{E} ,\,X\big)$ whose set of finite distinguished points in
the fiber over $x$ is equal to the finite branch locus $\{x\}\times \mathring{E}_x$. The
double branched covering map ${\sf{P}}: {\mathcal Y}_{\mathbold{E}}\to X\times \mathbb{P}^1$
maps the distinguished point $s_x^{\infty} \in \mathcal{P}^{-1}(x)$ to the point $\{x\}\times \infty$.

In the general case (i.e. without the assumption that the sets $\mathring{E}_x,\,x \in X$, are
uniformly bounded) the statement is proved by considering an exhaustion of $X$ by
relatively compact open sets.

Let $\big(X\times \mathbb{P}_1, {\rm pr}_1, \mathbold{E}_j,X\big)$, $j=0,1,$ be
two special holomorphic (smooth, respectively) $\,(0,4)$-bundles that are isotopic through smooth special $\,(0,4)$-bundles.
Then the
families $({\mathcal Y}_{\mathbold{E}_j},\mathcal{P}_j, \mathbold{s}_j^{\infty},\,X),\, j=0,1,$ of Riemann surfaces of type $(1,1)$ are isotopic.
Indeed, let $I$ be an open interval containing $[0,1]$.
The isotopy of the special $(0,4)$-bundles is given by a bundle $\big((I\times X)\times \mathbb{P}^1, {\rm pr}_1, {\mathbold{E}},I\times X\big)$ over $I\times X$. Here ${{\mathbold{E}}}$ is a smooth submanifold of $(I\times X)\times \mathbb{P}^1$ such that for each $(t,x)\in I\times X$ the intersection of the set ${\mathbold{E}}$ with the fiber over $(t,x)$ equals $\{(t,x)\}\times {E}(t,x)$ for subsets ${E}(t,x)$ that are the union of $\infty$ with an element of $C_3(\mathbb{C}\diagup \mathcal{S}_3)$. Moreover, ${{\mathbold{E}}}\cap ((\{j\}\times X) \times \mathbb{P}^1)=\mathbold{E}_j,\,j=0,1$.
The smooth special $(0,4)$-bundle $\big((I\times X)\times \mathbb{P}^1, {\rm pr}_1, {\mathbold{E}},I\times X\big)$ has a double branched covering by a smooth family of Riemann surfaces of type $(1,1)$ which defines the required isotopy for the families $({\mathcal Y}_{\mathbold{E}_j},\mathcal{P}_j, \mathbold{s}_j^{\infty},\,X),\, j=0,1,$ of Riemann surfaces of type $(1,1)$.

Notice that for a non-contractible finite open smooth surface $X$ and a
smooth special $(0,4)$-bundle $\big(X\times \mathbb{P}^1, {\rm pr}_1, \mathbold{E},X\big)$ over $X$ the obtained double branched covering
$({\mathcal Y}_{\mathbold{E}},\mathcal{P},\,\mathbold{s}^{\infty},\,X)$ is not the only double branched covering of the $(0,4)$-bundle (see Section \ref{sec:9.7}). We will call the bundle $({\mathcal Y}_{\mathbold{E}},\mathcal{P},\,\mathbold{s}^{\infty},\,X)$ the
canonical double branched covering of the special $(0,4)$-bundle
 $\big(X\times \mathbb{P}^1, {\rm pr}_1, \mathbold{E},X\big)$.

\section{Lattices and double branched coverings}
\label{sec:9.7}

Consider the quotient $\mathbb{C}\diagup \Lambda$ of the complex plane by a lattice $\Lambda$. We want to associate to the quotient
a double branched covering over the Riemann sphere with covering space being conformally
equivalent to $\mathbb{C}\diagup \Lambda$. The standard tool for this purpose is
the Weierstra\ss\ $\wp$-function  $\wp_{\Lambda}$.
Put $\Lambda\stackrel{def}=a {\mathbb Z} + b{\mathbb Z}$, where $a$ and $b$ are real linearly independent complex numbers.
The Weierstra\ss\
$\wp$-function,\index{Weierstra\ss \ $\wp$-function}
\begin{equation}\label{eqEl30}
\wp_{\Lambda} (\zeta) = \frac1{\zeta^2} + \sum_{(n,m)
\in {\mathbb Z}^2 \atop (n,m) \ne (0,0)} \left( \frac1{(\zeta -a\, n -b
\,m )^2} - \frac1{(a\,n+b\,m )^2} \right) , \quad \zeta \in
{\mathbb C} \setminus \Lambda\,,
\end{equation}
is meromorphic on ${\mathbb C}$, has poles
of second order at points of $ \Lambda$ and
is holomorphic on ${\mathbb C}\setminus \Lambda$. It has
periods $a$ and $b$ and
principal part $\zeta \to \frac1{\zeta^2}$ at $0$.
The summation is over all non-zero elements of the lattice. Hence, the function depends
only on the lattice, not on the choice of the generators $a$ and $b$ of the lattice.

Write the lattice in the form $\Lambda=\alpha \big(\mathbb{Z}+ \tau \mathbb{Z}\big)$.
Consider the lattice $\Lambda_{\tau}\stackrel{def}={\mathbb Z} + \tau \,{\mathbb Z}$. If the pair $1$ and $\tau$ generates the lattice, then also the pair $1$ and $-\tau$ generates it. Hence, we may assume that ${\rm{Re}}{\tau} >0$. Put $\wp_{\tau}\stackrel{def}=
\wp_{\Lambda_{\tau}}$,
\begin{equation}\label{eqEl30a}
\wp_{\tau} (\zeta) = \frac{1}{\zeta^2} + \sum_{(n,m)
\in {\mathbb Z}^2 \atop (n,m) \ne (0,0)} \left( \frac{1}{(\zeta -\, n -
\,m \, \tau)^2} - \frac{1}{(\,n+ \tau\,m )^2} \right) , \quad \zeta \in
{\mathbb C} \setminus \Lambda_{\tau}\,.
\end{equation}

The function $\wp_{\tau}$ satisfies the differential equation
\begin{equation}\label{eqEl31}
(\wp'_{\tau})^2 (\zeta) = 4 (\wp_{\tau} (\zeta) - e_1
(\tau))(\wp_{\tau} (\zeta) - e_2 (\tau)) (\wp_{\tau} (\zeta) - e_3
(\tau)) \, ,
\end{equation}
where
\begin{equation}\label{eqEl32}
e_1 (\tau) = \wp_{\tau} \left(\frac12 \right) \, ,
\quad e_2 (\tau) = \wp_{\tau} \left(\frac\tau2 \right) \, , \quad
e_3 (\tau) = \wp_{\tau} \left(\frac{1+\tau}2 \right)
\end{equation}
are the values of $\wp_{\tau}$ at points which are contained in the
lattice $\frac12 \, {\mathbb Z} + \frac\tau2 \, {\mathbb Z}$ but are
not contained in the lattice ${\mathbb Z} + \tau \, {\mathbb Z}$.
We put $w = \wp'_{\tau}(\zeta)$, $z = \wp_{\tau}(\zeta)$, and
\begin{equation}\label{eqEL31a}
\mathring{Y}(\Lambda_{\tau})\stackrel{def}= \left\{(z,w)\in \mathbb{C}^2: w^2 = 4 \, (z -
e_1 (\tau)) (z-e_2 (\tau)) (z-e_3
(\tau))\right\}\,.
\end{equation}

Consider an arbitrary lattice $\Lambda$ and write it as $\alpha(\mathbb{Z} +\tau \mathbb{Z})$.
The function $\wp_{\Lambda}$ can be expressed in terms of $\alpha$ and $\wp_{\tau}$
by the equality
\begin{equation}\label{eqEl29c}
\wp_{\Lambda}(\zeta)=\alpha^{-2} \, \wp_{\tau} \left( \frac{\zeta}{\alpha} \right),\,
\zeta \in {\mathbb C} \backslash \Lambda\,.
\end{equation}
The differential equation for this function is obtained by replacing in  \eqref{eqEl31}
$\wp_{\tau}(\zeta)$ by $\alpha^{-2}\wp_{\tau} \big( \frac{\zeta}{\alpha} \big)$,
$\wp_{\tau}'(\zeta)$ by $\alpha^{-3}\wp_{\tau}' \big( \frac{\zeta}{\alpha} \big)$ and $e_j(\tau)$ by $\alpha^{-2}e_j(\tau)$.
By equations \eqref{eqEl32} and \eqref{eqEl29c}  the unordered triple of the $\alpha^{-2}e_j(\tau)$ is equal to the unordered triple
of the values of $\wp_{\Lambda}(\zeta)= \alpha^{-2} \wp_{\tau}(\frac{\zeta}{\alpha})$ at points which are contained in the
lattice $\frac12 \,\alpha ({\mathbb Z} + \tau \, {\mathbb Z})$ but are
not contained in the lattice $\alpha({\mathbb Z} + \tau \, {\mathbb Z})$. Hence, the unordered triple depends only on the lattice $\Lambda$. We denote it by
\begin{equation}\label{eqEl29x}
\{e_1 (\Lambda),e_2 (\Lambda),e_3 (\Lambda)\}= \{\alpha^{-2}e_1(\tau),\alpha^{-2}e_2(\tau),\alpha^{-2}e_3(\tau)\}\,.
\end{equation}
The differential equation for $\wp_{\Lambda}$ becomes
\begin{equation}\label{eqEl29b}
(\wp'_{\Lambda})^2 (\zeta) = 4 (\wp_{\Lambda} (\zeta) - e_1
(\Lambda))(\wp_{\Lambda} (\zeta) - e_2 (\Lambda)) (\wp_{\Lambda} (\zeta) - e_3
(\Lambda)) \, .
\end{equation}

Put $z (\zeta)= \wp_{\Lambda}(\zeta)= \alpha^{-2} \, \wp_{\tau}
\left(\frac\zeta\alpha\right)$ and
$w(\zeta) = \wp'_{\Lambda}(\zeta) =\alpha^{-3} \, \wp'_{\tau}
\left(\frac\zeta\alpha\right)$. The mapping
\begin{equation}\label{eqEl29}
{\mathbb C} \backslash \Lambda \ni \zeta \to \big(z(\zeta),w(\zeta)\big)\,=\, \big(
 \wp_{\Lambda} ( \zeta ) ,
 \wp'_{\Lambda}({\zeta})
\big) \in {\mathbb C}^2 \,
\end{equation}
descends to a conformal mapping $\omega_{\Lambda}$ from the punctured torus
$({\mathbb C} \setminus \Lambda) \diagup \Lambda$ onto the complex hypersurface $\mathring{{Y}}({\Lambda})$ of ${\mathbb C}^2$,
\begin{align}\label{eqEl34}
\omega_{\Lambda}:(\mathbb{C}\setminus \Lambda)\diagup \Lambda \to       \mathring{{Y}}({\Lambda})\stackrel{def} =
\Big\{(z,w)\in \mathbb{C}^2: w^2  = 4 \prod_{j=1}^3\,(z -  \, e_j (\Lambda))
\Big\}\,.
\end{align}
To see this we notice first that the mapping $\omega_{\Lambda}$
is one-to-one. Indeed, the mapping $\wp_{\Lambda}$ descends to a holomorphic mapping ${\sf{p}}_{\Lambda}:(\mathbb{C}\setminus \Lambda)\diagup\Lambda \to \mathbb{C}$, such that $\wp_{\Lambda} = {\sf{p}}_{\Lambda}\circ p$. Here $p$ denotes the projection $p:\mathbb{C}\setminus \Lambda \to (\mathbb{C}\setminus \Lambda)\diagup \Lambda$.
The mapping ${\sf{p}}_{\Lambda}$ is $2$ to $1$. Further,
$\wp_{\Lambda}(\zeta)=\wp_{\Lambda}(-\zeta)$.
Hence, the preimage under $\wp_{\Lambda}$ of each point in $\mathbb{C}\setminus \{0\}$
equals $(\{\zeta\}+\Lambda)\cup (\{-\zeta\}+\Lambda)$.
Since,
$\wp_{\Lambda}'(\zeta)=-\wp_{\Lambda}'(-\zeta)$, the mapping $\omega_{\Lambda}$ is one-to-one.

Moreover, the mapping \eqref{eqEl29} is locally conformal. Indeed, if $\wp_{\Lambda}(\zeta_0)\neq e_j({\Lambda}),\, j=1,2,3,$ then $\wp_{\Lambda}'(\zeta_0)\neq 0$
and the mapping $\zeta\to \wp_{\Lambda}(\zeta)$ is conformal in a neighbourhood of $\zeta_0$. Suppose
$\wp_{\Lambda}(\zeta_0)=e_j(\Lambda)$. The differential equation \eqref{eqEl29b}  implies
that $\wp_{\Lambda}' \wp_{\Lambda}''= 2\wp_{\Lambda}'(\sum_{\ell=1}^3 \prod_{k\neq \ell}  (\wp_{\Lambda}-e_k(\Lambda))$. Hence, $\wp_{\Lambda}''(\zeta_0)=2 \prod_{k\neq j}(e_j(\Lambda)- e_k(\Lambda))\neq 0$. Hence, in this case  the mapping $\zeta\to \wp_{\Lambda}'(\zeta)$ is conformal in a neighbourhood of $\zeta_0$.

The set \eqref{eqEl34} is the covering space of the double branched covering
\begin{equation}\label{eqEl34a}
\mathring{{Y}}({\Lambda})\ni (z,w)\to z \in\mathbb{C}
\end{equation}
of ${\mathbb C}$ with branch locus
\begin{equation}\label{eqEl35}
\mathring{BL}(\Lambda)\stackrel{def} =\big\{e_1(\Lambda),e_2(\Lambda),e_3(\Lambda)\big\}\,.
\end{equation}
With $BL(\Lambda)\stackrel{def}=\mathring{BL}(\Lambda)\cup\{\infty\}$ the equality
$\mathring{{Y}}({\Lambda})=\mathring{{Y}}_{BL(\Lambda)}$ holds.
\index{$\mathring{BL}(\Lambda)$}

For a family of lattices $\Lambda (z)$ depending holomorphically on a complex parameter $z$ the sets $\mathring{BL}(\Lambda (z))$ (see \eqref{eqEl35})
depend holomorphically on
$z$. Indeed, locally we can write $\Lambda (z) = \alpha (z) ({\mathbb
Z} \, + \, \tau (z) \, {\mathbb Z})$ with $\alpha$ and $\tau$
holomorphic. The branch locus $\mathring{BL}(\Lambda (z))$ of the covering map \eqref{eqEl34a}
depends only on the lattice $\Lambda $, not on the choice of $\alpha$ and $\tau \in
\mathbb{C}_+$, and
can be locally expressed in
terms of $\tau$ and $\alpha$ by equations \eqref{eqEl32} and \eqref{eqEl29x}.
Since the function  $\wp_{\tau} (\zeta)$ depends
holomorphically on $\tau$ and $\zeta$ and the points $\frac12$,
$\frac\tau2$ and $\frac{\tau+1}2$ depend holomorphically on $\tau$,
the locally defined functions $\alpha^{-2} \, e_j (\tau)$
depend holomorphically on $\tau$ and $\alpha$.

If the family of lattices is merely smooth, then the set of finite branch points in the
fiber over $x$ depends smoothly on the point $x \in X$.

The one-point compactification $Y_{BL(\Lambda)}$ of $\mathring{Y}_{BL(\Lambda)}$ (see
Section \ref{sec:9.6}) is the covering manifold of a double branched covering of the Riemann
sphere $\mathbb{P}^1$ with branch locus $BL\stackrel{def}=\mathring{BL}(\Lambda)\cup\{\infty\}$ and finite branch locus $\mathring{BL}(\Lambda)$. We associate to the
Riemann surface $\mathbb{C}\diagup \Lambda$ the conformally equivalent Riemann surface
$Y_{BL(\Lambda)}$ which we also  denote by $Y(\Lambda)$. Notice that $\omega_{\Lambda}$ extends to a conformal mapping between the closed Riemann surfaces $\mathbb{C}\diagup \Lambda$ and $Y(\Lambda)$.

We describe now the lifts of the mapping classes of $\mathbb{P}^1$ with distinguished points $BL=\mathring{BL}(\Lambda)\cup \infty$ to mapping classes of $Y(\Lambda)$.

Notice that
for each lattice $\Lambda$ the mapping ${\mathbb C} \ni \zeta \to -
\zeta$ maps $\Lambda$ onto itself. This mapping descends to an
involution of ${\mathbb C} \diagup \Lambda$, i.e. to a
self-homeomorphism $\iota_{\Lambda}$ of ${\mathbb C} \diagup \Lambda$ such
that $\iota_{\Lambda}^2 = {\rm id}$.

\smallskip

Formula (\ref{eqEl30}) implies the following equality
\begin{equation}
\label{eqEl39} \left( \wp_{\Lambda}
\left(-\zeta\right) , \wp'_{\Lambda}
\left(-\zeta\right)\right) = \left(
\wp_{\Lambda} \left(\zeta\right) ,
-\wp'_{\Lambda} \left(\zeta\right)\right) \, .
\end{equation}

Conjugate the restriction $\iota_{\Lambda}\mid (\mathbb{C}\setminus \Lambda)\diagup \Lambda$ by the conformal mapping $\omega_{\Lambda}^{-1}$.
We obtain a self-homeomorphism of $\mathring{Y}({\Lambda})$ which we denote by $\iota$. The
mapping $\iota$ satisfies the equality
\begin{equation} \label{eqEl39a}
\iota(z,w)=(z,-w)\,.
\end{equation}
The involution $\iota$ extends to an involution of ${Y}({\Lambda})$, denoted also by $\iota$.
By \eqref{eqEl39a} the involution $\iota$ fixes the projection to
${\mathbb P}^1$ of each point of the double branched covering space $Y(\Lambda)$ and interchanges the
sheets over each point. Hence, it fixes each of the three finite branch
points, it also fixes $\infty$, and it does not fix any other point.

\smallskip

Let $\widetilde\varphi$ be any real linear self-homeomorphism of
${\mathbb C}$ which maps $\Lambda$ onto itself, and let $\varphi$ be
the induced mapping on ${\mathbb C} \diagup \Lambda$. Then $\varphi$
commutes with $\iota_{\Lambda}$. Indeed, $\widetilde\varphi (-\zeta) = -
\widetilde\varphi (\zeta)$, $\zeta \in {\mathbb C}$.

\smallskip
Each mapping class ${\mathfrak
m}$ in $\mathfrak{M}({\mathbb C} \diagup \Lambda;\,0\diagup \Lambda,\emptyset)$ can be represented by a
self-homeomorphism of ${\mathbb C} \diagup \Lambda$ which commutes
with $\iota_{\Lambda}$. Indeed, $\mathfrak{M}({\mathbb C} \diagup \Lambda;\;0\diagup \Lambda,\,\emptyset)$ contains a mapping $\varphi$ that lifts
to a real linear self-map
$\widetilde\varphi$ of ${\mathbb C}$ which maps $\Lambda$ onto
itself. This mapping $\varphi$ commutes with $\iota_{\Lambda}$.

\smallskip

As a corollary, each mapping class ${\mathfrak
m}$ in $\mathfrak{M}(Y(\Lambda);\, s^{\infty} ,\,\emptyset)\cong \mathfrak{M}({\mathbb C} \diagup \Lambda;\,0\diagup \Lambda,\emptyset)$
is a lift of a mapping class $\mathfrak{m}_1\in \mathfrak{M}(\mathbb{P}^1;\{\infty\},\mathring{E})$ for the set $\mathring{E}\stackrel{def}= \mathring{BL}(\Lambda)\subset C_3(\mathbb{C})\diagup \mathcal{S}_3$. Indeed, if $\varphi \in \mathfrak{m} \in  \mathfrak{M}({\mathbb C} \diagup \Lambda;\,0\diagup \Lambda,\,\emptyset)$ commutes with  $\iota_{\Lambda}$, then $\psi=\omega_{\Lambda}\circ \varphi \circ \omega_{\Lambda}^{-1}$ commutes with $\iota$ and represents the mapping class in $\mathfrak{M}(Y(\Lambda);\, s^{\infty} ,\,\emptyset)$ that corresponds to $\mathfrak{m}$. Denote by $\mathring{\psi}$ the restriction of $\psi$ to $\mathring{Y}$.

In coordinates $(z,w)$ on $\mathring {Y}(\Lambda)$ we write
\begin{equation}\label{eqEl19a}
\mathring{\psi}(z,w) = (\mathring{\psi}_1 (z,w) ,\mathring{ \psi}_2 (z,w)) \,.
\end{equation}
Since $\mathring{\psi}$ commutes with $\iota$ we obtain by (\ref{eqEl39a})
\begin{equation}\label{eqEl39d}
(\mathring{\psi}_1 (z,-w) ,\mathring{ \psi}_2 (z,-w)) = \mathring{\psi} \circ \iota (z,w) = \iota \circ \mathring{\psi} (z,w)=
(\mathring{\psi}_1 (z,w) , -\mathring{\psi}_2 (z,w))\,.
\end{equation}
Hence $\mathring{\psi}_1(z,w)=\mathring{\psi}_1(z,-w)$. Since each point in $(z,w)\in \mathring{Y}_{\Lambda}$ is determined by $z$ and the sign of $w$, this means that $\mathring{\psi}_1(z,w)$ depends only on the coordinate $z \in {\mathbb C}$, not on the sheet (determined by $w$). Further, if $\iota$ fixes
$(z,w)$, then $w=0$, and by  \eqref{eqEl39d} for $w=0$, $\mathring{\psi}_2(z,w)=\mathring{\psi}_2 (z,-w)= -\mathring{\psi}_2(z,w)$. In other words, if $\iota$ fixes $(z,w)$, then it also fixes $\mathring{\psi}(z,w)$. Hence, $\mathring{\psi}$
maps the set of finite branch points (the preimage of $\mathring{BL}({\Lambda})$ under the branched covering map) onto itself, and its extension $\psi$ maps the preimage $s^{\infty}$ of
$\infty$ to itself.
We saw, that
${\psi}$ induces a self-homeomorphism ${\sf P}\psi$ of ${\mathbb P}^1$ in the
class ${\mathfrak M} ({\mathbb P}^1 ; \{\infty\} , \mathring{BL}({\Lambda}))$,
${\sf P}\psi |\mathbb{C}= {\sf P} \mathring{\psi} $, ${\sf P} \mathring{\psi}(z)=\mathring{\psi}_1(z,w), \, (z,w)\in \mathbb{C}^2$.
We call ${\sf{P}}(\psi)$ the projection of $\psi$. \index{$ {\sf{P}}(\psi)$}
The mapping class $\mathfrak{m}$ is a lift of the mapping class $\mathfrak{m}_1$ of $\psi_1$. Notice that the provided arguments imply also the following fact. For two mapping
classes $\mathfrak{m}^1,\mathfrak{m}^2\in\mathfrak{M}(Y(\Lambda);s^{\infty}, \emptyset)$ the equality ${\sf P}(\mathfrak{m}^1 \mathfrak{m}^2)={\sf P}(\mathfrak{m}^1){\sf P}(\mathfrak{m}^2)$ holds. Indeed,
take representing maps
$\psi^1$ and $\psi^2$, that commute with the involution $\iota$. Using the arguments given above, the equality ${\sf P}(\psi^1\circ \psi^2)={\sf P}(\psi^1)\circ {\sf P}(\psi^2)$ follows.

Let again $\mathring{E}$ be an element of $C_3(\mathbb{C})\diagup \mathcal{S}_3$ and $E=\mathring{E}\cup \{\infty\}$. Any mapping class $\mathfrak{m}_1\in \mathfrak{M}(\mathbb{P}^1;\{\infty\},\mathring{E})$ has exactly two lifts $\mathfrak{m}_{\pm} \in \mathfrak{M}(Y_E;s^\infty,\emptyset)$. Here $Y_E$ is the double branched covering of $\mathbb{P}^1$ with finite branch locus $\mathring{E}$ and $s^\infty$ is the branch point over $\infty$. Indeed, for any representative $\psi_1$ of $\mathfrak{m}_1$ there are exactly two self-homeomorphisms of the double branched covering that lift $\psi_1$.
They are obtained as follows. Cut $\mathbb{P}^1$ along two disjoint simple curves $\gamma_1$ and $\gamma_2$, that join disjoint pairs of points of $E$. $\psi_1$ maps the pair of curves $\gamma_1$ and $\gamma_2$ to another pair of curves $\gamma'_1$ and $\gamma'_2$ joining (maybe different) disjoint pairs of points of $E$. The double branched covering over $\mathbb{P}^1$ is obtained in two different ways, either gluing two sheets of $\mathbb{P}^1\setminus (\gamma_1\cup\gamma_2)$ crosswise together, or gluing two sheets of $\mathbb{P}^1\setminus (\gamma'_1\cup\gamma'_2)$ crosswise together. $\varphi_1$ maps $\mathbb{P}^1\setminus (\gamma_1\cup\gamma_2)$ homeomorphically onto $\mathbb{P}^1\setminus (\gamma'_1\cup\gamma'_2)$. Consider the mapping that takes the first sheet of
$\mathbb{P}^1\setminus (\gamma_1\cup\gamma_2)$ onto the first sheet of $\mathbb{P}^1\setminus (\gamma'_1\cup\gamma'_2)$ and the second sheet of
$\mathbb{P}^1\setminus (\gamma_1\cup\gamma_2)$ onto the second sheet of $\mathbb{P}^1\setminus (\gamma'_1\cup\gamma'_2)$ and lifts $\varphi_1\mid \mathbb{P}^1\setminus (\gamma_1\cup\gamma_2)$. This mapping extends to a self-homeomorphism $\varphi$ of the double branched covering that lifts $\psi_1$.
There is exactly one more self-homeomorphism of the double branched covering that lifts
$\psi_1$. This self-homeomorphism maps the first sheet of
$\mathbb{P}^1\setminus (\gamma_1\cup\gamma_2)$ onto the second sheet of $\mathbb{P}^1\setminus (\gamma'_1\cup\gamma'_2)$ and the second sheet of
$\mathbb{P}^1\setminus (\gamma_1\cup\gamma_2)$ onto the first sheet of $\mathbb{P}^1\setminus (\gamma'_1\cup\gamma'_2)$.
In other words, there are
two lifts of $\psi_1$ to the double branched covering of ${\mathbb P}^1$ with branch locus 
$E$, and they differ by involution.

The following proposition relates Theorem \ref{thmEl.0} to the respective theorem for
elliptic bundles that will be formulated below.
\begin{prop}\label{propEl.2} Let $X$ be a Riemann surface (smooth surface, respectively) of genus $g$ with ${m}\geq 1$ holes with base point $x_0$ and curves $\gamma_j$ representing a set of generators $e_j\in \pi_1(X,x_0)$.\\
\noindent $(1)$ Each complex analytic (differentiable, respectively) family of Riemann surfaces of type $(1,1)$ over $X$ is holomorphically (smoothly, respectively) isomorphic to the double branched covering $\mathfrak{F}$ of a special holomorphic (smooth, respectively) $(0,4)$-bundle $\mathfrak{F}_1$ over $X$. The monodromy of the bundle $\mathfrak{F}$ along each $\gamma_j$ is a lift of the respective monodromy of the bundle $\mathfrak{F}_1$.\\
\noindent $(2)$ Vice versa, for each special holomorphic (smooth, respectively) $(0,4)$-bundle over $X$ and each collection $\mathfrak{m}^j$ of lifts of the $2{ g} +{ m}-1 $ monodromy mapping classes $\mathfrak{m}_1^j$ of the bundle along the $\gamma_j$ there exists a double branched covering by
a complex analytic (differentiable, respectively) family of Riemann surfaces of type $(1,1)$ with collection of monodromy mapping classes equal to the $\mathfrak{m}^j$.
For each given holomorphic (smooth, respectively) special $(0,4)$-bundle over $X$ there are up to holomorphic (smooth, respectively) isomorphisms exactly
$2^{2{g}+ {m}-1}$ holomorphic (smooth, respectively) families of Riemann surfaces of type $(1,1)$ that lift the $(0,4)$-bundle.\\
\noindent $(3)$ A lift of a special $(0,4)$-bundle is reducible if and only if the special
 $(0,4)$-bundle is reducible.
\end{prop}

\noindent {\bf Proof.}
We start with the proof of the first statement of the proposition.
Consider a complex analytic (smooth, respectively) family of Riemann surfaces of type $(1,1)$ over $X$. By Lemma \ref{lem1} and Proposition \ref{prop1} we may assume that the family has the form 
$(\mathcal{X}_{\Lambda},\mathcal{P}_{\Lambda}, \,\mathbold{s}_{\Lambda}\,,X)\,$
for a
holomorphic (smooth, respectively) family of lattices $\Lambda(x),\, x \in X$.

We consider the complex (smooth, respectively) submanifold $\mathring{\mathbold{BL}}(\Lambda)$ of $X\times \mathbb{C}$ that intersects each fiber $\{x\}\times \mathbb{C}$ along
the set $\{x\}\times \mathring{BL}(\Lambda(x))$ (see equations \eqref{eqEl29x} and \eqref{eqEl35})
and  the complex (smooth, respectively) submanifold $\mathbold{BL}(\Lambda)$ of $X\times \mathbb{P}^1$ that intersects each fiber $\,\{x\}\times \mathbb{P}^1\,$ along the set $\,\{x\}\times BL(\Lambda(x))\,$ with $BL(\Lambda(x))\,=\, \mathring{BL}(\Lambda(x))\cup \{\infty\}$. The holomorphic (smooth, respectively) family $(\mathcal{X}_{\Lambda},\mathcal{P}_{\Lambda}, \,\mathbold{s}_{\Lambda}\,,X)\,$ is holomorphically isomorphic (isomorphic as a smooth family of Riemann surfaces, respectively) to the canonical double branched covering  $(\mathcal{Y}_{\mathbold{BL}(\Lambda)}, \mathcal{P},\mathbold{s}^{\infty},X)$ (see Section \ref{sec:9.6}) of the
special $(0,4)$-bundle $(X\times \mathbb{P}^1, {\rm pr}_1, \mathbold{BL}(\Lambda),X  )$.
To see this we recall that for each $x\in X$ the mapping $\omega_{\Lambda_x}:(\mathbb{C}\setminus \Lambda_x)\diagup \Lambda_x \to \mathring{Y}_{{BL}(\Lambda_x)}$ is a surjective conformal mapping that depends holomorphically (smoothly, respectively) on $x$ (see equation \eqref{eqEl29}). Hence,
the mapping defines a holomorphic (smooth, respectively) homeomorphism from $\mathcal{X}_{\Lambda}\setminus \mathbold{s}_{\Lambda}$ onto $\mathring{\mathcal{Y}}_{\mathbold{BL}(\Lambda)}$, that maps the fiber over $x$ of the first bundle conformally onto the fiber over $x$ of the second bundle.

The set ${\mathcal{Y}}_{\mathbold{BL}(\Lambda)}$
is obtained as in Section \ref{sec:9.6}
by ''adding a point to each fiber''.
As in Section \ref{sec:9.6} the mapping extends to a holomorphic (smooth, respectively) homeomorphism $\mathcal{X}_{\Lambda} \to  {\mathcal{Y}}_{\mathbold{BL}(\Lambda)}$,
that maps the fiber over $x$ conformally onto the fiber over $x$. The extension maps $\mathbold{s}_{\Lambda}$ to $\mathbold{s}^{\infty}$.
We proved that the bundles $(\mathcal{X}_{\Lambda},\mathcal{P}_{\Lambda}, \,\mathbold{s}_{\Lambda}\,,X)\,$ and $(\mathcal{Y}_{\mathbold{BL}(\Lambda)}, \mathcal{P},\mathbold{s}^{\infty},X)$ are holomorphically isomorphic (isomorphic as smooth families of Riemann surfaces, respectively).
Hence, each holomorphic (smooth, respectively) $(1,1)$-family over a finite open Riemann surface is holomorphically isomorphic (smoothly isomorphic), and in particular, isotopic to a double branched covering of a holomorphic (smooth, respectively) $(0,4)$-bundle.

We identify the mapping class groups in the fiber over the base point $x_0$ of the bundles
$(\mathcal{X}_{\Lambda},\mathcal{P}_{\Lambda}, \,\mathbold{s}_{\Lambda}\,,X)\,$ and $(\mathcal{Y}_{\mathbold{BL}(\Lambda)}, \mathcal{P},\mathbold{s}^{\infty},X)$ using the described isomorphism.
Having in mind this identification, we  prove  now  that  the  monodromy  mapping  class along $\gamma_j$ of  the  $(1,1)$-bundle
$(\mathcal{X}_{\Lambda},\mathcal{P}_{\Lambda}, \,\mathbold{s}_{\Lambda}\,,X)\,$
is a lift of the monodromy mapping class along $\gamma_j$ of the special $(0,4)$-bundle $(X\times \mathbb{P}^1, {\rm pr}_1, \mathbold{BL}(\Lambda),X  )$. Parameterise $\gamma_j$ by the unit interval $[0,1]$.
Write $\Lambda_t\stackrel{def}=\Lambda(\gamma_j(t))=a(t)\mathbb{Z}+b(t)\mathbb{Z},\, t\in[0,1],$ for  smooth functions $a$ and $b$ on $[0,1]$.
For $t \in [0,1]$ we denote by $\widetilde{\varphi}^t$ the real linear self-homeomorphism of $\mathbb{C}$ that maps $a(0)$ to $a(t)$ and $b(0)$ to $b(t)$. Let $\varphi^t$ be the homeomorphism from the fiber over $\gamma_j(0)$ onto the fiber over $\gamma_j(t)$ of the $(1,1)$-bundle that lifts to $\widetilde{\varphi}^t$. Then ${\varphi}^0= {\rm Id}$ and ${\varphi}^1$ represents the monodromy mapping class of the $(1,1)$-bundle $(\mathcal{X}_{\Lambda},\mathcal{P}_{\Lambda}, \,\mathbold{s}_{\Lambda}\,,X)\,$
along $\gamma_j$. Let $\mathring{\psi}^t:\mathring{Y}(\Lambda_0) \to \mathring{Y}(\Lambda_t)$ be
obtained from the commutative diagram
$$
\xymatrix{(\mathbb{C}\setminus \Lambda_0)\diagup \Lambda_0
 \ar[rr]^{\mathring{\varphi}^t}  \ar[d]^{\omega_{\Lambda_0}} &&(\mathbb{C}\setminus \Lambda_t)\diagup \Lambda_t \ar[d]^{\omega_{\Lambda_t}}\\
\ar[rr]^{\mathring{\psi}^t} {\mathring{Y}(\Lambda_0)}  && \mathring{Y}(\Lambda_t)}
$$
where $\mathring{\varphi}^t=\varphi^t|(\mathbb{C}\setminus \Lambda_0)\diagup \Lambda_0$.

We use the restriction to $\mathring{Y}(\Lambda_0)$ of the coordinates $(z,w)$ on $\mathbb{C}^2$, and let $\iota_0$ be the involution on $\mathring{Y}(\Lambda_0)$ and $\iota_t$ the involution on $\mathring{Y}(\Lambda_t)$. Choose respective coordinates on
$\mathring{Y}(\Lambda_t)$ and write $\mathring{\psi}^t(z,w)=(\mathring{\psi}^t_1(z,w), \mathring{\psi}^t_2(z,w)), \; (z,w) \in \mathring{Y}(\Lambda_0)$. Since $\widetilde{\varphi}^t(-\zeta)=-\widetilde{\varphi}^t(\zeta), \zeta \in\mathbb{C},$ we obtain
\begin{equation}\label{eqEl19d}
\mathring{\psi}^t\,\iota_0= \iota_t \, \mathring{\psi}^t\,.
\end{equation}
In coordinates $(z,w)$ this means 
$(\mathring{\psi}^t_1(z,-w), \mathring{\psi}^t_2(z,-w))=
(\mathring{\psi}^t_1(z,w), -\mathring{\psi}^t_2(z,w))$. Hence, $\mathring{\psi}^t_1(z,w)=\mathring{\psi}^t_1(z,-w)$, so that $\mathring{\psi}^t_1$ depends only on the coordinate $z \in \mathbb{C}$.
Moreover, if $\iota_0$ fixes $(z,w)$ (equivalently, if $z\in \mathring{BL}(\Lambda_0)$ ), then $w=0$, hence $\mathring{\psi}^t_2(z,w)=\mathring{\psi}^t_2(z,-w)=-\mathring{\psi}^t_2(z,w)$, i.e. $\iota_t$ fixes $\mathring{\psi}^t(z,w)$ (equivalently, $\mathring{\psi}^t_1(z,w)\in \mathring{BL}(\Lambda_t)$).
We saw that
the self-homeomorphism $\mathring{\psi}^t_1$ of $\mathbb{C}$ maps the set  $\mathring{BL}(\Lambda_0)$ onto the set $\mathring{BL}(\Lambda_t)$.
Each $\mathring{\psi}^t_1$ extends to a self-homeomorphism $\psi^t_1$ of $\mathbb{P}^1$ that maps the distinguished points $BL(\Lambda_0)\stackrel{def}= \mathring{BL}(\Lambda_0)\cup\{\infty\}$ onto the distinguished points $BL(\Lambda_t)\stackrel{def}= \mathring{BL}(\Lambda_t)\cup\{\infty\}$. Hence,  $\psi^1_1$ represents the monodromy mapping class along $\gamma_j$ of the special $(0,4)$-bundle with set of finite distinguished points $\mathring{BL}(\Lambda_t)$ in the fiber over $\gamma_j(t)$.
The mapping $\mathring{\psi}^1=(\mathring{\psi}_1^1,\mathring{\psi}_2^1)$ is equal to $\omega_{\Lambda_0} \circ \mathring{\varphi}^1 \circ \omega_{\Lambda_0}^{-1}$. Identifying the mapping class groups of $\mathbb{C}\diagup \Lambda_0$ with the mapping class group of $Y_{BL(\Lambda_0)}$ by the isomorphism induced by $\omega_{\Lambda_0}$, we may identify the monodromy of the bundle
$(\mathcal{X}_{\Lambda},\mathcal{P}_{\Lambda}, \,\mathbold{s}_{\Lambda}\,,X)\,$ along $\gamma_j$ with the mapping class of the extension $\psi^1$ of
$\mathring{\psi}^1=(\mathring{\psi}_1^1,\mathring{\psi}_2^1)$. The
mapping $\psi_1^1$ is the projection of
$\psi^1$ under the double branched covering from $\mathbb{C}\diagup\Lambda_0$ with distinguished point $0\diagup\Lambda_0$ onto $\mathbb{P}^1$ with distinguished points ${BL}(\Lambda_0)\cup\{\infty\}$. We proved that the monodromy mapping class of the bundle $(\mathcal{X}_{\Lambda}, \mathcal{P}_{\Lambda}, \mathbold{s}_{\Lambda},X)$ along each $\gamma_j$ is a lift of the respective monodromy mapping class of the bundle $(X\times \mathbb{P}^1, {\rm pr}_1, \mathbold{BL}(\Lambda),X)$.

We prove now the second statement.
Consider the special holomorphic, or smooth, respectively, $(0,4)$-bundle $\mathfrak{F}_1\stackrel{def}=\big(X\times \mathbb{P}^1, {\rm pr}_1, \mathbold{E}, X\big)$. The complex or smooth, respectively, submanifold $\mathbold{E}$ of $X\times \mathbb{P}^1$  intersects each fiber $\{x\}\times \mathbb{P}^1$ along the set $\{x\}\times E_x$, with  $E_x= (\mathring{E}_x\cup\{\infty\})$, where $\mathring{E}_x\subset C_3(\mathbb{C})\diagup \mathcal{S}_3$.
In Section \ref{sec:9.6} we obtained the canonical lift $(\mathcal{\mathcal{Y}}_{\mathbold{E}},\mathcal{P},\mathbold{s}^{\infty},X)$. Denote the monodromy mapping class along $\gamma_j$ of the canonical lift
by $\mathfrak{m}^j_+$.
Recall that each mapping class $\mathfrak{m}_1^j$ has two lifts and they differ by involution.
Let $\mathfrak{m}^j_-$ be the mapping class whose representatives differ from those of $\mathfrak{m}^j_+$ by composition with the involution $\iota$ of the fiber over the base point.

Take any subset $J$ of $\{1,2,\ldots,2{g}+{m}-1 \}$. The lift of the $(0,4)$-bundle whose monodromy mapping class along $\gamma_j$ equals $\mathfrak{m}^j_+$ if $j\notin J$, and
$\mathfrak{m}^j_-$ if $j\in J$, is obtained as follows.
Let $\tilde U$ be the subset of the universal covering
$\tilde X  \overset{{ p}}{-\!\!\!\longrightarrow} X$ that was used in Section \ref{sec:9.2}, and let $D$, and $\tilde{D}_j,\,j=0,\ldots,$ be as in Section \ref{sec:9.2}. Consider the lift of the bundle $(\mathcal{\mathcal{Y}}_{\mathbold{E}},\mathcal{P},\mathbold{s}^{\infty},X)$
to $\tilde X$ and restrict the lifted bundle to $\tilde U$. Denote the obtained bundle on $\tilde U$ by $\tilde{\mathfrak{F}}$.
For each point $x\in D$ we let $\tilde{x}_j\in \tilde{D}_j,\, j=0,\ldots,$ be the points
with $p(\tilde{x}_j)=x$.
We identify the fiber of $\tilde{\mathfrak{F}}$ over $\tilde{x}_j,\,j=0,\ldots,$ with the fiber ${Y}_x$, $x=p(\tilde{x}_j)$, of the canonical double branched covering of the
$(0,4)$-bundle. $\mathring{Y}_x$ is obtained from $Y_x$ by removing the point of $\mathbold{s}^{\infty}$ from $Y_x$.

Glue for each $x\in D$ and each $j=1,\ldots,$ the fiber of $\tilde{\mathfrak{F}}$ over $\tilde{x}_j$ to the fiber of $\tilde{\mathfrak{F}}$ over $\tilde{x}_0$ using the identity if $j\notin J$ and the mapping $\iota$
if $j\in J$.
More detailed, the gluing mapping of the punctured fibers in the case $j\in J$ equals
\begin{equation}\label{eqEl19e}
(\tilde{x}_j,z,w)\to (\tilde{x}_0,z,-w),\;\; x\in D,\; (z,w)\in \mathring{Y}_x\,.
\end{equation}
Since the gluing mappings are holomorphic (smooth, respectively), we obtain a complex analytic (smooth, respectively) family over $X$ of double branched coverings of $\mathbb{C}$.
Extend the family to a complex analytic (smooth, respectively) family of double branched coverings of $\mathbb{P}^1$ over $X$. This can be done as in Section \ref{sec:9.6} by considering
the families $\mathcal{Y}^{r_n}_{\mathbold{E}|{X}_n}$ for an exhausting sequence of relatively compact open subsets $X_n$ of $X$ and a suitable sequence of real numbers $r_n$.
We obtained
for each collection of lifts of the  $\mathfrak{m}^j_1$ a complex analytic (smooth, respectively) family of Riemann surfaces of type $(1,1)$ with given fiber over the base point and with monodromy mapping classes equal to this collection. It is easy to see that
two holomorphic (smooth, respectively) $(1,1)$-bundles that lift a given $(0,4)$-bundle over a finite open Riemann surface (smooth surface, respectively) are holomorphically isomorphic
(smoothly isomorphic, respectively), if and only if their monodromies coincide.
Hence, for each given holomorphic (smooth, respectively) special $(0,4)$-bundle over $X$ there are up to holomorphic (smooth, respectively) isomorphism exactly
$2^{2{g}+ {m}-1}$ holomorphic (smooth, respectively) families of Riemann surfaces of type $(1,1)$ that lift the $(0,4)$-bundle.

It remains to prove the third statement. Let ${\sf P}:(\mathcal{X},\mathcal{P},\mathbold{s},X)\to (X\times \mathbb{P}^1,{\rm pr}_1,\mathbold{E},X)$ be a double branched covering of a $(0,4)$-bundle by a $(1,1)$-family of Riemann surfaces.
If the special $(0,4)$-bundle is reducible, then there is a simple closed curve $\gamma$ that
divides the fiber $\mathbb{P}^1$ over the base point $x_0$ into two connected components, each of which contains two
distinguished points. The curve $\gamma$ is mapped by each monodromy mapping class of the $(0,4)$-bundle to a class of isotopic curves.
The preimage of the closed curve $\gamma$ under the covering map consists of two simple closed curves
$\tilde{\gamma}_1$ and $\tilde{\gamma}_2$, that are homotopic to each other in $\mathcal{P}^{-1}(x)\setminus\mathbold{s}$, and each of
the two curves cuts the torus into an annulus. Each monodromy mapping class of the $(1,1)$-bundle, being a lift
of the respective monodromy mapping class of the $(0,4)$-bundle, takes each $\tilde{\gamma}_j$ to a class of
isotopic curves. Hence, the $(1,1)$-bundle is also reducible.

Suppose now that a double branched covering $\mathfrak{F}=(\mathcal{X},\mathcal{P},\mathbold{s},X)$ of
a special $(0,4)$-bundle ${\sf P}(\mathfrak{F})=(X\times \mathbb{P}^1, {\rm pr}_1, \mathbold{E},X)$
is reducible and prove that the $(0,4)$-bundle is reducible. Let
$x_0\in X$ be the base point of $X$, and let the fiber of the bundle $\mathfrak{F}$ over $x_0$
be $\mathbb{C}\diagup\Lambda$ with distinguished point $0\diagup\Lambda$.
Any admissible system of curves on the torus $\mathbb{C}\diagup\Lambda$ with a distinguished point contains exactly one curve.

Let $\gamma_0$ be a simple closed curve in $(\mathbb{C}\setminus \Lambda)\diagup\Lambda$  that reduces all monodromy mapping classes of the bundle $\mathfrak{F}$.
The curve  $\gamma_0$ is free isotopic in $\mathbb{C}\diagup\Lambda$ to a curve $\gamma$ with base point $0\diagup \Lambda$ that lifts under the covering map $p:\mathbb{C}\to \mathbb{C}\diagup\Lambda$ to the straight line segment that joins $0$ with another lattice point $\lambda$ in $\Lambda$. Since $\gamma_0$ is free isotopic in $\mathbb{C}\diagup\Lambda$ to
a simple closed curve on the torus, $\lambda$ is a primitive element of the lattice. Indeed, with $\Lambda=\{na+mb:n,m\in \mathbb{Z}\}$ for real linearly independent complex numbers $a$ and $b$ we get $\lambda=n(\lambda)a+m(\lambda)b$, where $n(\lambda)$ and $m(\lambda)$ are relatively prime integer numbers. Then there is another element $\lambda'=n(\lambda')a+m(\lambda')b\in \Lambda$ such that
\[
\begin{vmatrix}
n(\lambda) & n(\lambda') \\
m(\lambda) & m(\lambda') \\
\end{vmatrix}
= 1\,.
\]
\noindent The two elements $\lambda$ and $\lambda'$ generate $\Lambda$. Multiplying $\Lambda$ by a non-zero complex number and changing perhaps $\lambda'$ to $-\lambda'$, we may assume that $\lambda=1$ and $\lambda'=\tau$ with ${\rm Im}\tau>0$. In other words, after changing $\mathbb{C}\diagup\Lambda$  to the conformally equivalent Riemann surface $ \mathbb{C}\diagup\Lambda_{\tau}$, we may assume that the fiber of $\mathfrak{F}$ over $x_0$ equals $\mathbb{C}\diagup\Lambda_{\tau}$. After a free isotopy in $(\mathbb{C}\setminus\Lambda_{\tau}) \diagup\Lambda_{\tau}$ we may assume that the reducing curve $\gamma_0$ lifts under the covering map $p:\mathbb{C}\to \mathbb{C}\diagup \Lambda_{\tau}$ to the segment $\frac{\tau}{2}+[0,1]\subset \mathbb{C}$.
Denote by $\gamma$ the curve in $\mathbb{C}\diagup \Lambda_{\tau}$ that lifts to $[0,1]$,
and by $\gamma'$ the curve on the torus that lifts to $[0,\tau]$. $\gamma$ and $\gamma'$ represent a pair of generators of the fundamental group of the closed torus with base point $0\diagup \Lambda_{\tau}$.

The complement $(\mathbb{C}\diagup\Lambda_{\tau}) \setminus \gamma_0$ of $\gamma_0$ in the fiber over $x_0$ is a topological annulus with distinguished point $0\diagup  \Lambda_{\tau}$. Since $\gamma_0$ reduces all monodromy mapping classes of the $(1,1)$-family $\mathfrak{F}$, each monodromy mapping class has a representative that fixes $\gamma_0$ pointwise. These representatives map the topological annulus $(\mathbb{C}\diagup\Lambda_{\tau}) \setminus \gamma_0$ homeomorphically onto itself, fixing the distinguished point and fixing the boundary pointwise, or fixing the boundary pointwise after an involution. Hence, each monodromy is a power of a Dehn twist about $\gamma_0$, maybe composed with an involution.

The
Dehn twist on ${\mathbb C} \diagup \Lambda_{\tau}$  about $\gamma_0$ can be represented by a
self-homeomorphism $\psi_{\tau}$ of ${\mathbb C} \diagup \Lambda_{\tau}$
which lifts under $ p$ to the real linear self-map $\widetilde{\psi}_{\tau}$ of
${\mathbb C}$ which maps $1$ to $1$ and $\tau$ to $1+\tau$. This can
be seen by looking at the action of $\widetilde\psi_{\tau}$ on the lifts of the
curves $\gamma$ and $\gamma'$. $\widetilde\psi_{\tau}$ takes $\widetilde{\gamma}=[0,1]$ to itself, and $\widetilde{\gamma}'=[0,\tau]$ to $[0,\tau+1]$.
Note that the
real $2 \times 2$ matrix corresponding to $\widetilde\psi_{\tau}$ is
$\begin{pmatrix} 1 &({\rm Im} \, \tau)^{-1} \\ 0 & 1 \end{pmatrix}$.

We consider the double branched covering ${\sf{p}}_{\Lambda_{\tau}}: ({\mathbb C}\setminus \Lambda_{\tau})
\diagup \Lambda_{\tau}\to \mathbb{C}$ with branch locus $\mathring{BL}(\Lambda_{\tau})=\{e_1(\tau), e_2(\tau), e_3(\tau)\}$ determined by the Weierstra\ss\ $\wp$-function, ${\sf{p}}_{\Lambda_{\tau}}= \wp_{\Lambda_{\tau}}\circ p$ (with $p$ being the projection $p: \mathbb{C}\setminus \Lambda_{\tau} \to (\mathbb{C}\setminus \Lambda_{\tau})\diagup \Lambda_{\tau}$) and its extension ${\sf{p}}_{\Lambda_{\tau}}^c: {\mathbb C} \diagup \Lambda_{\tau}\to \mathbb{P}^1$. The Weierstrass $\wp$-function extends to a mapping, also denoted by $\wp_{\Lambda_{\tau}}$,
$\wp_{\Lambda_{\tau}}: \mathbb{C}\to \mathbb{P}^1$ that takes $0$ to $\infty$, $\frac{1}{2}$ to $e_1$, $\frac{\tau}{2}$ to $e_2$, and $\frac{1+\tau}{2}$ to $e_3$.
The line segments $[0,\frac{1}{2}]$, $[0,\frac{\tau}{2}]$,
$[\frac{1}{2}, \frac{1+\tau}{2}]$,
and $[\frac{\tau}{2}, \frac{1+\tau}{2}]$
are mapped under ${\wp}_{\Lambda_{\tau}}$  to simple curves $\gamma_{\infty ,e_1}$, $\gamma_{\infty ,e_2}$, $\gamma_{e_1,e_3}$, and
$\gamma_{e_2,e_3}$
in $\mathbb{P}^1$, each of which joins the first mentioned point with the last mentioned point. The union of the four curves (each with suitable orientation) is a closed curve in $\mathbb{P}^1$ that divides $\mathbb{P}^1$ into two connected components $\mathcal{C}_1$ and $\mathcal{C}_2$. The Weierstrass $\wp$-function ${\wp}_{\Lambda_{\tau}}$ maps the open parallelogram $R$ in the complex plane with vertices $0$, $1$, $\frac{\tau}{2} $, $1+\frac{\tau}{2}$ conformally onto its image. Indeed, the "double" of this parallelogram, i.e. the parallelogram with vertices $0$, $1$, $\tau $, $\tau +1$, is the interior of a fundamental polygon for the covering $p:\mathbb{C}\to \mathbb{C}\diagup {\Lambda_{\tau}}$ and ${\wp}_{\Lambda_{\tau}}(\zeta)= {\wp}_{\Lambda_{\tau}}(-\zeta + n + m \tau ) ,\, n,m\in \mathbb{Z}  $.
By the last equation we get in particular
\begin{equation}\label{eqEl19d'}
{\wp}_{\Lambda_{\tau}}(\frac{\tau}{2}+t)={\wp}_{\Lambda_{\tau}}(\frac{\tau}{2}-t)\,,
\end{equation}
$t\in [0,\frac{1}{2}]$. After perhaps relabeling $\mathcal{C}_1$ and $\mathcal{C}_2$, we may assume that ${\wp}_{\Lambda_{\tau}}$ maps the open parallelogram $R_-$ with vertices $0,\,\frac{1}{2},\frac{\tau}{2},\frac{1+\tau}{2},$ (the ''left half'' of $R$) onto $\mathcal{C}_1$, and then it maps the parallelogram $R_+$ with vertices $\frac{1}{2},\,1,\, \frac{1+\tau}{2},1+\frac{\tau}{2},$ (the ''right half'' of $R$) onto $\mathcal{C}_2$.

The mapping $\tilde{\psi}_{\tau}$ fixes the segment $[0,\frac{1}{2}]$, and translates the segment $[\frac{\tau}{2}, \frac{1+\tau}{2}]$ by $\frac{1}{2}$. It maps $R_-$ to the parallelogram $\tilde{\psi}_{\tau}(R_-) \subset R$ whose horizontal sides are $[0,\frac{1}{2}]$ and $[\frac{1+\tau}{2}, \frac{2+\tau}{2}]$.
${\wp}_{\Lambda_{\tau}}$ maps
$[\frac{1}{2}+\frac{\tau}{2}, 1+ \frac{\tau}{2}]$
to $\gamma_{e_3,e_2}$ which equals $\gamma_{e_2,e_3}$ with inverted
orientation (see \eqref{eqEl19d'}).
It maps the segment $[0,\frac{1+\tau}{2}]$ to a curve $\gamma_{\infty, e_3}$ which is contained in $\mathcal{C}_1$ except of its endpoints and joins $\infty$ and $e_3$, and the
segment $[\frac{1}{2},1+\frac{\tau}{2}]$ is mapped to a curve $\gamma'_{e_1,e_2}$ with initial point $e_1$ and terminating point $e_2$,
which is contained in  $\mathcal{C}_2$ except its endpoints, so that
the union of the curves $\gamma_{\infty,e_1}$, $\gamma'_{e_1,e_2}$,
$\gamma_{e_2,e_3}$, and $\gamma_{e_3,\infty}$
is the oriented boundary of a domain $\mathcal{C}'_1$
(i.e. the domain $\mathcal{C}'_1$ is ''on the left'' when walking along the oriented curve). Here $\gamma_{e_3,\infty}$ equals $\gamma_{\infty,e_3}$
with inverted orientation.
The mapping ${\wp}_{\Lambda_{\tau}}\circ \tilde{\psi}_{\tau}$ takes $R_-$ to the domain $\mathcal{C}'_1$.
The mapping $\varphi_{\tau,1}\stackrel{def}={\wp}_{\Lambda_{\tau}}\circ \tilde{\psi}_{\tau} \circ ({\wp}_{\Lambda_{\tau}}|R_-)^{-1}: \mathcal{C}_1\to \mathcal{C}'_1 $ is homotopic on $\mathcal{C}_1$ to a mapping that fixes a (large) neighbourhood of $\gamma_{\infty,e_1}$ in $\mathcal{C}_1$, maps a neighbourhood of $\gamma_{e_2,e_3}$ in $\mathcal{C}_1$ into a neighbourhood of $\gamma_{e_2,e_3}$ in $\mathbb{C}$,
such that it
rotates (in suitable coordinates) a smaller neighbourhood of $\gamma_{e_2,e_3}$ in $\mathcal{C}_1$ by the angle $\pi$ around $e_3$, and rotates points in the remaining set by an angle between $0$ and $\pi$.

A similar argument for the parallelogram with vertices $\pm\frac{1}{2}$, $\frac{\tau}{2}\pm \frac{1}{2}$ and the domain $\mathcal{C}_2$ shows that with $R_+ -1\stackrel{def}=\{z-1:z\in R_+\}$ the mapping $\varphi_{\tau,2}\stackrel{def}={\wp}_{\Lambda_{\tau}}\circ \tilde{\psi}_{\tau} \circ ({\wp}_{\Lambda_{\tau}}|R_+ -1)^{-1} $
is homotopic on $\mathcal{C}_2$ to a mapping with the same properties as above
for the role of $\mathcal{C}_1$ and $\mathcal{C}_2$ interchanged.

The mapping $\varphi_{\tau}$ that is equal to $\varphi_{\tau,j}$ on $\mathcal{C}_j$ extends to a self-homeomorphism of $\mathbb{P}^1$ whose mapping class is conjugate to $\mathfrak{m}_{\sigma_1}$.

Since $\tilde{\psi}_{\tau}|R_{-}=(p|R)^{-1} \circ \psi_{\tau}\circ p|R_- $
and ${\wp}_{\Lambda_{\tau}}\circ \tilde{\psi}_{\tau}|R_-= \varphi_{\tau,1}\circ{\wp}_{\Lambda_{\tau}}|R_-$ we obtain
${\wp}_{\Lambda_{\tau}}\circ (p|R)^{-1}\circ \psi_{\tau} =\varphi_{\tau,1}\circ{\wp}_{\Lambda_{\tau}}\circ (p|R_-)^{-1}  $,
hence, since ${\wp}_{\Lambda_{\tau}}={\sf p}_{\Lambda_{\tau}}\circ p$, the equality
${\sf p}_{\Lambda_{\tau}} \circ \psi_{\tau}=\varphi_{\tau,1}\circ{\sf p}_{\Lambda_{\tau}}$ holds on $p(R_-)$. By the same reason the equality holds for $p(R_+)$, and since $p(R_-\cup R_+)$ is dense in $\mathbb{P}^{1}$, it holds on $\mathbb{P}^{1}$. We proved that the projection of $\psi_{\tau}$ is homotopic to a conjugate of $\mathfrak{m}_{\sigma_1}$, which is a reducible mapping class. We proved that the  $(0,4)$-bundle is reducible.
Proposition \ref{propEl.2} is proved. \hfill $\Box$

\medskip

Let ${\mathfrak F}$ be a $(\sf{g},\sf{m})$-fiber bundle over
$\partial {\mathbb D}$ and $\widehat{\mathfrak F}$ its isomorphism class. Let $\mathfrak{m}_{\mathfrak{F}}$ be the monodromy mapping class of the bundle $\mathfrak{F}$ and $\widehat{\mathfrak{m}_{\mathfrak{F}}}$ its conjugacy class.
By the bijective correspondence of Theorem \ref{thmEl1} we will write also
\begin{equation}\label{eqEl9}
{\mathcal M} (\widehat{{\mathfrak m}_{\mathfrak F}})
= {\mathcal M} (\widehat{\mathfrak F}) \,
\end{equation}
for the conformal module ${\mathcal M} (\widehat{\mathfrak F})$ of the bundle.
Corollary \ref{corEl2} below shows that for $(1,1)$-bundles this invariant descends to the conformal module of braids.

Let $\widehat{\mathfrak{F}}$ be a smooth isomorphism class of $(1,1)$-bundles over the circle,
and $\widehat{{\mathfrak m}_{\mathfrak F}}$ the associated conjugacy class of mapping classes. For any representative ${\mathfrak m}_{\mathfrak F}\in \widehat{{\mathfrak m}_{\mathfrak F}}$ and any mapping class $\mathfrak{m}$ the equality ${\sf P}(\mathfrak{m}^{-1}{\mathfrak m}_{\mathfrak F}\mathfrak{m})={\sf P}( \mathfrak{m})^{-1}       {\sf P} ({\mathfrak m}_{\mathfrak F})       {\sf P}({\mathfrak m})$ holds. Hence, $  {\sf P}(\widehat{{\mathfrak m}_{\mathfrak F}})$ is well-defined. Let ${\sf P}(\widehat{\mathfrak{F}})$ be the smooth isomorphism class of special $(0,4)$-bundles that corresponds to $  {\sf P}(\widehat{{\mathfrak m}_{\mathfrak F}})$. Recall that a representative of the latter can be obtained as the projection of
a smooth family of Riemann surfaces of type $(1,1)$ that is
in the class of $\widehat{\mathfrak{F}}$ and is a double branched covering of a special $(0,4)$-bundle.
The class ${\sf P}(\widehat{\mathfrak{F}})$ defines a conjugacy class of $3$-braids $\hat b$ that is defined up to multiplication by $\Delta_3^{2k}$ for integers $k$. For the
conformal module $\mathcal{M}(\hat{b})$ of any such conjugacy class of braids and the conformal module $\mathcal{M}(\widehat{\mathfrak{F}})$ of the smooth isomorphism class
$\widehat{\mathfrak{F}}$ the following corollary holds.
\begin{cor}\label{corEl2}
\begin{align*}
\mathcal{M}(\widehat{\mathfrak{F}})=\mathcal{M}(\hat b)\,.
\end{align*}
\end{cor}

\noindent {\bf Proof.}
$\mathcal{M}(\widehat{\mathfrak{F}})$ is the supremum of the conformal modules of annuli on which there exists a holomorphic $(1,1)$-bundle that represents $\mathfrak{F}$.
$\mathcal{M}({{\sf P}(\widehat{\mathfrak{F}})})$ is the supremum of the conformal modules of annuli on which there exists a special holomorphic $(0,4)$-bundle that represents ${\sf P}(\widehat{\mathfrak{F}})$. By \cite{Jo} we have the equality $\mathcal{M}(\hat{b})=\mathcal{M}(\widehat{b \Delta_3^{2k}})=\mathcal{M}({{\sf P}(\widehat{\mathfrak{F}})})$ for any braid $b$ that is associated to a bundle representing the class ${{\sf P}(\widehat{\mathfrak{F}})}$.
Statement (1) of Proposition \ref{propEl.2} implies the inequality $\mathcal{M}(\widehat{\mathfrak{F}})\leq \mathcal{M}({{\sf P}(\widehat{\mathfrak{F}})})$, Statement (2) implies the opposite inequality. The corollary is proved. \hfill $\Box$

\medskip

In the following Theorem \ref{thmEl.9} we consider isotopies of smooth $(1,1)$-families
over a torus with a hole to complex analytic $(1,1)$-families.
Theorem \ref{thmEl.9} is a more comprehensive version of Theorem \ref{thmEl.9a} from the introduction. (Theorem 4 states the first assertion of Theorem \ref{thmEl.9}.)

\begin{thm}\label{thmEl.9}
$(1)$ Let $X$ be a smooth surface of genus one with a hole with base point $x_0$, and with a chosen set
$\mathcal{E}=\{e_1,e_2\}$ of generators of $\pi_1(X,x_0)$. Define the set
$\mathcal{E}_0=
\{e_1,e_2, e_1 e_2^{-1},  e_1 e_2^{-2}, e_1e_2 e_1^{-1} e_2^{-1}\}$ as in Theorem {\rm \ref{thm8.0}}. Consider a smooth $(1,1)$-bundle $\mathfrak{F}$
on $X$. Suppose for each $e\in \mathcal{E}_0$ the restriction of the bundle $\mathfrak{F}$ to an annulus representing $\hat e$ has the Gromov-Oka property.
Then the bundle $\mathfrak{F}$ has the Gromov-Oka property on $X$.\\
$(2)$ If a bundle $\mathfrak{F}$ as in $(1)$ is irreducible, then it is isotopic
to an isotrivial bundle, and hence, for each conformal structure $\omega$ on $X$ the bundle
$\mathfrak{F}_{\omega}$ is isotopic
to a bundle that extends to a holomorphic $(1,1)$-bundle on
the closed torus $\omega(X)^{c}$. In particular, $\mathfrak{F}$ is isotopic to a holomorphic bundle for any conformal structure $\omega$ on $X$ (maybe, of first kind). \\
$(3)$ Any smooth reducible bundle $\mathfrak{F}$ on $X$ has a single irreducible bundle component. This
irreducible bundle component is isotopic
to an isotrivial bundle. There is a Dehn twist in the fiber over the base point such that the
$(1,1)$-bundle $\mathfrak{F}$ can be
recovered from the irreducible bundle component up to composing each monodromy by a power of the Dehn twist.
A smooth reducible bundle on $X$ is isotopic to a holomorphic $(1,1)$-bundle for each conformal structure of second kind on $X$.\\
$(4)$ A reducible holomorphic $(1,1)$-bundle over a punctured Riemann surface is locally holomorphically trivial.
\end{thm}

\noindent {\bf Proof.}
By Lemma \ref{lem1} and Proposition \ref{propEl.2} we may assume that $\mathfrak{F}$ is a double branched covering of a smooth special $(0,4)$-bundle, denoted by ${\sf P}(\mathfrak{F})$.
Suppose that for each $e\in \mathcal{E}_0,$
the restriction of the smooth $(1,1)$-bundle $\mathfrak{F}$ to an annulus $A_{\hat{e}}$ representing $\hat e$ has the Gromov-Oka property. Then there is a conformal structure $\omega_e:A_{\hat{e}}\to \omega_e(A_{\hat{e}})$ such that
$\omega_e(A_{\hat{e}})$ has conformal module bigger than $\frac{\pi}{2}\log(\frac{3+\sqrt{5}}{2})^{-1}$ and the pushed forward bundle
$(\mathfrak{F}|A_{\hat{e}})_{\omega_e}$
is isotopic
to a holomorphic bundle $\mathfrak{F}_e$  on
$\omega_e( A_{\hat{e}}  )$. By 
Proposition \ref{prop1} we may assume that
$\mathfrak{F}_e$ is of the form \eqref{eqE111e}, and by Proposition
\ref{propEl.2} we may assume that
$\mathfrak{F}_e$ is a double branched covering of a holomorphic special $(0,4)$-bundle on $\omega_e(A_{\hat{e}})$,
denoted by ${\sf P}(\mathfrak{F}_e)$.
Since the monodromy homomorphisms of $\mathfrak{F}|A_{\hat{e}}$ and $\mathfrak{F}_e$ are conjugate,
the monodromy homomorphisms of ${\sf P}(\mathfrak{F}|A_{\hat{e}})$ and ${\sf P}(\mathfrak{F}_e)$ are also conjugate. Hence, ${\sf P}(\mathfrak{F}|A_{\hat{e}})$ and ${\sf P}(\mathfrak{F}_e)$ are isotopic.
We saw that for each $e\in \mathcal{E}$ the restriction of the special $(0,4)$-bundle ${\sf P}(\mathfrak{F})$ to an annulus representing $\hat e$ is isotopic to a holomorphic
bundle for a conformal structure on the annulus of conformal module bigger than $\frac{\pi}{2}\log(\frac{3+\sqrt{5}}{2})^{-1}$.
By Theorem \ref{thmEl.0} the special $(0,4)$-bundle ${\sf P}(\mathfrak{F})$ is isotopic to a holomorphic special $(0,4)$-bundle for any conformal structure of second kind on $X$.
Hence, the canonical double branched covering of ${\sf P}(\mathfrak{F})$ is isotopic to a holomorphic $(1,1)$-bundle for each conformal structure of second kind on $X$.
Since the monodromy mapping class of any double branched covering of the $(0,4)$-bundle along each loop $\gamma$ in $X$ with base point $x_0$ differs at most by involution from the monodromy mapping class along $\gamma$ of the canonical double branched covering,
by Proposition \ref{propEl.2} also the bundle $\mathfrak{F}$ is isotopic to a holomorphic $(1,1)$-bundle for each conformal structure of second kind on $X$.
Statement (1) is proved.

Suppose a bundle $\mathfrak{F}$ as in (1) is irreducible. Then by Proposition \ref{propEl.2}
the special $(0,4)$-bundle ${\sf{P}}(\mathfrak{F})$ is also irreducible. Moreover,
the $(0,4)$-bundle ${\sf{P}}(\mathfrak{F})$
is isotopic to an isotrivial bundle, i.e. its lift to a finite covering $\hat X$ is isotopic to the trivial bundle, and therefore, for each conformal structure $\omega$ on $X$
(maybe, of first kind)
the $(0,4)$-bundle ${\sf{P}}(\mathfrak{F})_{\omega}$ is isotopic to a bundle
that extends holomorphically to the closed torus $\omega(X)^c$. Then also the canonical double branched covering of the $(0,4)$-bundle is isotopic to an isotrivial bundle.
Hence, the bundle $\mathfrak{F}$ is isotopic to an isotrivial bundle. In particular, for each conformal structure (maybe, of first kind) on $X$ the bundle $\mathfrak{F}_{\omega}$ is isotopic to a holomorphic bundle on $\omega(X)$.
This gives statement (2).

Suppose now $\mathfrak{F}$ is any smooth reducible bundle on $X$. Then ${\sf P}(\mathfrak{F})$ is reducible. By Theorem  \ref{thmEl.0} the special $(0,4)$-bundle ${\sf P}(\mathfrak{F})$ is isotopic to a holomorphic special $(0,4)$-bundle for any conformal structure of second kind on $X$. Hence, for any conformal structure of second kind on $X$ the bundle $\mathfrak{F}$  is isotopic to a holomorphic $(1,1)$-bundle which is a double branched covering.

Each admissible set of curves on a punctured torus consists of exactly one curve. Hence, there is an admissible closed curve $\gamma$ in the punctured
fiber $\mathcal{P}^{-1}(x_0)\setminus \mathbold{s}^{\infty}$ over the base point $x_0$,
that is mapped by the monodromy mapping class along each curve $\gamma_j$ representing an element $e_j\in \mathcal{E}$ to a curve that is isotopic to $\gamma$ in the punctured fiber over $x_0$. Choose for each monodromy mapping class a representative that maps the curve $\gamma$ onto itself.

The complement of $\gamma$ on the punctured torus $\mathcal{P}^{-1}(x_0)\setminus \mathbold{s}^{\infty}$ is connected and homeomorphic to the thrice punctured Riemann sphere. This implies first that there is a single irreducible component of $\mathfrak{F}$. Further, the restrictions of suitable representatives of the monodromy mapping classes to the complement of $\gamma$ on the punctured torus are conjugate to self-homeomorphisms of the thrice punctured Riemann sphere. The conjugated homeomorphisms extend to self-homeomorphisms
of the Riemann sphere with three distinguished points. The self-homeomorphisms may interchange the two points that come from different edges of $\gamma$, but each of them fixes the third distinguished point. Except the identity there is only one isotopy class of such mappings and its square is the identity.
Hence, the monodromy mapping classes of the irreducible bundle component are powers of a single periodic mapping class.
In other words, the irreducible bundle component is isotopic to an isotrivial bundle.

The monodromy mapping classes of the original bundle $\mathfrak{F}$ are powers of Dehn twists about $\gamma$,
maybe, composed with the involution, and up to powers of Dehn twists about $\gamma$ in the fiber $\mathcal{P}^{-1}(x_0)$ the bundle can be recovered from the irreducible bundle component.
Statement (3) is proved.

We prove now statement (4). Consider a reducible holomorphic (1,1)-bundle $\mathfrak{F}$
over a punctured Riemann surface of genus $1$. By Proposition \ref{propEl.2}
it is holomorphically isomorphic to the double branched covering of a holomorphic special $(0,4)$-bundle ${\sf P}(\mathfrak{F})$ over a punctured Riemann surface.
By Proposition \ref{propEl.2} the special holomorphic$(0,4)$-bundle ${\sf P}(\mathfrak{F})$ is also reducible.
By Theorem \ref{thmEl.0} the bundle ${\sf P}(\mathfrak{F})$ is holomorphically trivial.
Then the double branched covering $\mathfrak{F}$ of ${\sf P}(\mathfrak{F})$ is locally holomorphically trivial.
Theorem \ref{thmEl.9} is proved. \hfill $\Box$

\end{document}